\documentclass[10pt,a4paper]{amsart}
\usepackage{amssymb,amsmath}
\usepackage[american,french]{babel}
\usepackage{amsopn}
\usepackage{graphicx}
\usepackage{fancyhdr}
\usepackage[T1]{fontenc}
\usepackage{mathrsfs}
\usepackage{color}

\title{About quadratic differential operators}
\newcommand{\rr}{\mathbb{R}}
\newcommand{\eps}{\varepsilon}
\newcommand{\nn}{\mathbb{N}}
\newcommand{\cc}{\mathbb{C}}

\def\un{{\mathrm{1~\hspace{-1.4ex}l}}}

\def\wrtext#1{\relax\ifmmode{\leavevmode\hbox{#1}}\else{#1}\fi}

\def\begeq{\begin{equation}}
\def\endeq{\end{equation}}

\def\part#1{\frac{\partial}{\partial #1}}

\newcommand{\real}{\mathbb{R}}

\begin{document}
\newcounter{equa}
\selectlanguage{american}

\newtheorem{Corollary}{Corollary}
\newtheorem{lemma}{Lemma}
\newtheorem{definition}{Definition}
\newtheorem{proposition}{Proposition}
\newtheorem{theorem}{Theorem}

\title[HYPOELLIPTIC ESTIMATES FOR A LINEAR MODEL OF THE BOLTZMANN EQUATION]{HYPOELLIPTIC ESTIMATES FOR A LINEAR MODEL OF THE BOLTZMANN EQUATION WITHOUT ANGULAR CUTOFF}

\author{Nicolas Lerner, Yoshinori Morimoto, Karel Pravda-Starov}

\address{\noindent \textsc{Institut de Math\'ematiques de Jussieu, 
Universit\'e Pierre et Marie Curie (Paris VI),  
4 Place Jussieu, 
75252 Paris cedex 05, 
France}}
\email{lerner@math.jussieu.fr}

\address{\noindent \textsc{Graduate School of Human and Environmental Studies,
Kyoto University, Kyoto 606-8501, Japan}}
\email{morimoto@math.h.kyoto-u.ac.jp }

\address{\noindent \textsc{
Universit\'e de Cergy-Pontoise,
CNRS UMR 8088,
D\'epartement de Math\'ematiques,
95000 Cergy-Pontoise, France}}
\email{karel.pravda-starov@u-cergy.fr}


\maketitle

\noindent
\textsc{Abstract.} In this paper, we establish optimal hypoelliptic estimates for a class of kinetic equations, which are simplified linear models for the spatially inhomogeneous Boltzmann equation without angular cutoff.

\section{Introduction}

In this paper, we study the following kinetic operator 
\begin{equation}\label{yo2}
P=\partial_t+v \cdot \partial_x  + a(t,x,v)(-\tilde{\Delta}_v)^{\sigma}, \ t \in \rr, \ x,v \in \rr^{n},
\end{equation}
where $0<\sigma<1$ is a constant parameter, $x \cdot y$ stands for the standard dot-product on $\rr^n$ and $a$ denotes a $C_b^{\infty}(\rr^{2n+1})$ function satisfying 
\begin{equation}\label{eq0.5}
\exists a_0>0, \forall (t,x,v) \in \rr^{2n+1}, \ a(t,x,v) \geq a_0>0.
\end{equation}
Here the notation $C_b^{\infty}(\rr^{2n+1})$ stands for the space of $C^{\infty}(\rr^{2n+1})$ functions whose derivatives of any order are bounded on $\rr^{2n+1}$ and 
$(-\tilde{\Delta}_v)^{\sigma}$ is the Fourier multiplier with symbol
\begin{equation}\label{eq0}
F(\eta)=|\eta|^{2\sigma}w(\eta)+|\eta|^2\big(1-w(\eta)\big), \ \eta \in \rr^n,
\end{equation}
with $|\cdot|$ being the Euclidean norm, $w$ a  $C^{\infty}(\rr^n)$ function satisfying $0 \leq w \leq 1$, $w(\eta)=1$ if $|\eta| \geq 2$, $w(\eta)=0$ if $|\eta| \leq 1$; and $D_t=(2\pi i)^{-1}\partial_t$, $D_x=(2\pi i)^{-1}\partial_x$, $D_v=(2\pi i)^{-1}\partial_v$.

When $\sigma=1$, this operator reduces to the so-called Vlasov-Fokker-Planck operator, whereas  when $0<\sigma<1$, it stands for a simplified linear model of the spatially inhomogeneous Boltzmann equation without angular cutoff (see the end of this introduction together with section~\ref{kkboltz} in appendix). This is our motivation for studying the regularizing properties of this linear model and establishing hypoelliptic estimates with optimal loss of derivatives with respect to the exponent $0<\sigma<1$ of the fractional Laplacian $(-\tilde{\Delta}_v)^{\sigma}$. This linear model has the familiar structure 
$$\textrm{Transport part in the }(t,x) \textrm{ variables}\ + \textrm{ Elliptic part in the }v \textrm{ variable}$$ 
and it is easy to get the regularity in the $v$ variable. The non-commutation of the transport part (the skew-adjoint part) with the self-adjoint elliptic part $(-\tilde{\Delta}_v)^{\sigma}$ will produce the regularizing effect in the $x$ variable.

Regarding this linear model, the existence and the $C^{\infty}$ regularity for the solutions of the Cauchy problem to linear and semi-linear equations associated with the operator (\ref{yo2}) were proved in \cite{MoXu}. H.~Chen, W.-X.~Li and C.-J.~Xu have also recently studied its Gevrey hypoellipticity. More specifically, they established in~\cite{chen} (Proposition~2.1) the following hypoelliptic estimate. Let $P$ be the operator defined in (\ref{yo2}) and $K$ a compact subset of $\rr^{2n+1}$. For any $s \geq 0$, there exists a positive constant $C_{K,s}>0$  such that for any $u \in C_0^{\infty}(K)$,
\begin{equation}\label{yo4}  
\|u\|_{s+\delta} \leq C_{K,s}\big(\|Pu\|_{s}+\|u\|_{s}\big),
\end{equation}
with $\|\cdot\|_{s}$ standing for the $H^s(\rr^{2n+1})$ Sobolev norm and 
\begin{equation}\label{yo4.5}
\delta=\textrm{max}\Big(\frac{\sigma}{4},\frac{\sigma}{2}-\frac{1}{6}\Big)>0.
\end{equation}
The notation $C_0^{\infty}(K)$ stands for the set of $C_0^{\infty}(\rr^{2n+1})$ functions with support in $K$. This hypoelliptic estimate with loss of 
$$\max(2\sigma,1)-\delta>0,$$ 
derivatives is then a key instrumental ingredient for their investigation of the Gevrey hypoellipticity of the operator $P$. 
However, this hypoelliptic estimate (\ref{yo4}) is not optimal. 
In the present work, we are interested in establishing hypoelliptic estimates with optimal loss of derivatives with respect to the exponent $0<\sigma<1$ of the fractional Laplacian $(-\tilde{\Delta}_v)^{\sigma}$. More specifically, we shall show by using different microlocal techniques that the operator $P$ is hypoelliptic with a loss of 
$$\frac{\max(4\sigma^2,1)}{(2\sigma+1)}>0,$$ 
derivatives, that is, that the hypoelliptic estimates (\ref{yo4}) hold with the new positive gain 
\begin{equation}\label{yo6}
\delta=\frac{2\sigma}{2\sigma+1}>0,
\end{equation}
which improves for any $0<\sigma<1$ the gain provided by (\ref{yo4.5}),
$$\frac{2\sigma}{2\sigma+1}>\textrm{max}\Big(\frac{\sigma}{4},\frac{\sigma}{2}-\frac{1}{6}\Big).$$
Our main result reads as follows.

\bigskip

\begin{theorem}\label{TTH1}
Let $P$ be the operator defined in \emph{(\ref{yo2})}, $K$ be a compact subset of $\rr^{2n+1}$ and $s \in \rr$. Then, there exists a positive constant $C_{K,s}>0$ such that for all $u \in C_0^{\infty}(K)$,
\begin{equation}\label{yo5}
\big\|(1+|D_t|^{\frac{2\sigma}{2\sigma+1}}+|D_x|^{\frac{2\sigma}{2\sigma+1}}+|D_{v}|^{2\sigma})u\big\|_{s} \leq C_{K,s} \big(\|Pu\|_{s}+\|u\|_{s}\big),
\end{equation}
with  $\|\cdot\|_{s}$ being the $H^s(\rr^{2n+1})$ Sobolev norm.
\end{theorem}

\bigskip

The hypoelliptic estimates (\ref{yo5}) are optimal in term of the exponents of derivative terms appearing in their left-hand-side, namely, $2\sigma/(2\sigma+1)$ for the regularity in the time and space variables and $2\sigma$ for the regularity in the velocity variable. The exponent $2\sigma$ for the regularity in the velocity variable has indeed the same growth as the diffusive part of the kinetic operator (\ref{yo2}). Regarding the optimality of the exponent $2\sigma/(2\sigma+1)$ for the regularity in the time and space variables, we first notice that Theorem~\ref{TTH1} is a natural extension for the values of the parameter $0<\sigma<1$ of the well-known optimal hypoelliptic estimates with loss of $4/3$ derivatives known for the Vlasov-Fokker-Planck operator, case $\sigma=1$, (see~\cite{bouchut,chen1,perthame}),
$$\big\|(1+|D_t|^{2/3}+|D_x|^{2/3}+|D_{v}|^{2})u\big\|_{s} \leq C_{K,s} \big(\|Pu\|_{s}+\|u\|_{s}\big).$$
See also~\cite{landau} for general microlocal methods for proving optimal hypoelliptic estimates with loss of $4/3$ derivatives for certain classes of kinetic equations.

We deduce the optimality of the exponent $2\sigma/(2\sigma+1)$ for the regularity in the time and space variables by using a simple scaling argument. Indeed, let us consider the specific case when the function $a$ appearing in the definition of the kinetic operator $P$ is constant and assume that the hypoelliptic estimates (\ref{yo4}) hold for a positive gain $\delta>0$. It follows that the estimate 
$$\|(|D_t|^{\delta}+|D_x|^{\delta}+|D_{v}|^{\delta})u\|_{L^2} \leq C\big(\|iD_tu+iv \cdot D_xu  + |D_v|^{2\sigma}u\|_{L^2}+\|u\|_{L^2}\big),$$
holds for any $u \in C_0^{\infty}(\rr^{2n+1})$  with support in the closed unit Euclidean ball $B_1=\overline{B(0,1)}$ in $\rr^{2n+1}$. The symplectic invariance of the Weyl quantization (Theorem~18.5.9 in~\cite{hormander}) shows that for any $\lambda \geq 1$, 
$$T_{\lambda}^{-1}\big(iD_t+iv \cdot D_x  + |D_v|^{2\sigma}\big)T_{\lambda}=\lambda^{\frac{2\sigma}{2\sigma+1}}\big(iD_t+iv \cdot D_x  + |D_v|^{2\sigma}\big)$$
and
$$T_{\lambda}^{-1}\big(|D_t|^{\delta}+|D_x|^{\delta}+|D_{v}|^{\delta}\big)T_{\lambda}=\lambda^{\frac{2\sigma\delta}{2\sigma+1}}|D_t|^{\delta}+\lambda^{\delta}|D_x|^{\delta}+\lambda^{\frac{\delta}{2\sigma+1}}|D_{v}|^{\delta},$$
where $T_{\lambda}$ stands for the following unitary transformation on $L^{2}(\rr^{2n+1})$, 
$$T_{\lambda}u(t,x,v)=\lambda^{\frac{2\sigma+n}{2(2\sigma+1)}+\frac{n}{2}}u(\lambda^{\frac{2\sigma}{2\sigma+1}}t,\lambda x,\lambda^{\frac{1}{2\sigma+1}}v).$$
Recalling that $\lambda \geq 1$, we notice that the $C_0^{\infty}(\rr)$ function $T_{\lambda}u$ is supported in $B_1$ if the $C_0^{\infty}(\rr)$ function $u$ is supported in $B_1$.
By applying the previous a priori estimate to functions $T_{\lambda}u$, we obtain that 
$$\|(|D_t|^{\delta}+|D_x|^{\delta}+|D_{v}|^{\delta})T_{\lambda}u\|_{L^2} \leq C\big(\|(iD_t+iv \cdot D_x  + |D_v|^{2\sigma})T_{\lambda}u\|_{L^2}+\|T_{\lambda}u\|_{L^2}\big).$$
By using that $T_{\lambda}^{-1}$ is a unitary transformation on $L^{2}(\rr^{2n+1})$, we deduce that for any $\lambda \geq 1$, 
\begin{multline*}
\|(\lambda^{\frac{2\sigma\delta}{2\sigma+1}}|D_t|^{\delta}+\lambda^{\delta}|D_x|^{\delta}+\lambda^{\frac{\delta}{2\sigma+1}}|D_{v}|^{\delta})u\|_{L^2} \\\ \leq C\big(\lambda^{\frac{2\sigma}{2\sigma+1}}\|(iD_t+iv \cdot D_x  + |D_v|^{2\sigma})u\|_{L^2}+\|u\|_{L^2}\big).
\end{multline*}
This estimate may hold for any $\lambda \geq 1$ only if 
$$\delta \leq \frac{2\sigma}{2\sigma+1}.$$
This scaling argument shows that the positive gain $2\sigma/(2\sigma+1)>0$ in the hypoelliptic estimates (\ref{yo4}) is the optimal possible one.

As a consequence of these optimal hypoelliptic estimates, we obtain the following result where we write
$$f \in H^s_{\textrm{loc}, (t_0,x_0)}(\rr^{2n+1}_{t,x,v}),$$
if there exists an open neighborhood $U$ of the point $(t_0,x_0)$ in $\rr^{n+1}$ such that
$\phi(t,x) f\in H^s(\rr_{t,x,v}^{2n+1})$
for any $\phi \in C_0^\infty(U)$.

\bigskip

\begin{Corollary}\label{ev1}
Let $P$ be the operator defined in \emph{(\ref{yo2})} and $N \in \nn$. 
If $u \in H_{-N}(\rr_{t,x,v}^{2n+1})$ and $Pu \in H^s_{\emph{\textrm{loc}}, (t_0,x_0)}(\rr^{2n+1}_{t,x,v})$ with $s \ge 0$, then there exists 
an integer $k \geq 1$ such that 
\[
\frac{u}{\langle  v \rangle^k} \in H^{s+\frac{2\sigma}{2\sigma+1}}_{\emph{\textrm{loc}}, (t_0,x_0)}(\rr^{2n+1}_{t,x,v}),
\]
where $\langle v \rangle=(1+|v|^2)^{1/2}$.
In particular, if $u \in H_{-N}(\rr^{2n+1})$ and $Pu \in H^{\infty}(\rr^{2n+1})$ then $u \in C^{\infty}(\rr^{2n+1})$.
\end{Corollary}

\bigskip

Corollary~\ref{ev1} allows to recover the $C^{\infty}$ hypoellipticity proved in \cite{MoXu} (Theorem~1.2) with now optimal loss of derivatives. Notice that the equation (\ref{yo2}) is not a classical pseudodifferential equation. Indeed, the coefficient $v$ in (\ref{yo2}) is unbounded and the fractional Laplacian $(-\tilde{\Delta}_v)^{\sigma}$ is a classical pseudo-differential operator in the velocity variable $v$ but not in all the variables $t,x,v$. This accounts for parts of the difficulties encountered when studying this kinetic operator in particular when using cutoff functions in the velocity variable. This also accounts for the weight $\langle  v \rangle^{-k}$ appearing in the statement of Corollary~\ref{ev1}. Notice that the proof of this result is constructive and that one may derive an explicit (possibly not sharp) bound on the integer $k \geq 1$.

The proof of Theorem~\ref{TTH1} is relying on some microlocal techniques developed by N.~Lerner for proving energy estimates while using the Wick quantization \cite{cubo}. Let us mention that some of these techniques were already used in the work~\cite{lms}. The strategy for proving Theorem~\ref{TTH1} is the following. We want to consider this equation as an evolution equation along the characteristic curves of $\partial_t+v \cdot \partial_x$; for that purpose, we straighten this vector field and get the normal form
$$iD_t+a(t,D_{x_2},D_{x_1}+tD_{x_2})F(x_2-tx_1).$$
This normal form suggests to derive some a priori estimates for the one-dimensional first-order differential operator 
$$iD_t+a(t,\xi_2,\xi_1+t\xi_2)F(x_2-tx_1),$$
depending on the parameters $x_1,x_2,\xi_1,\xi_2 \in \rr^n$. We then deduce from those a priori estimates some a priori estimates for the operator 
$$iD_t+\big[a(t,\xi_2,\xi_1+t\xi_2)F(x_2-tx_1)\big]^{\textrm{Wick}},$$
defined by using the Wick quantization. As a last step, we need to control some remainder terms in order to come back to the standard quantization and derive a priori estimates for the original operator
$$iD_t+a(t,D_{x_2},D_{x_1}+tD_{x_2})F(x_2-tx_1).$$
For the sake of completeness and to keep the paper essentially self-contained, the definition and all the main features of the Wick quantization are recalled in appendix (Section~\ref{appendix}). Next section is devoted to the proof of a key hypoelliptic estimate (Proposition~\ref{th1}) which is the core of the present work. This key estimate is then the main ingredient in Section~\ref{section3} for proving Theorem~\ref{TTH1}. Corollary~\ref{ev1} is established in Section~\ref{section4}.

Before ending this introduction, we give some references and comments about the hypoelliptic properties of the non-cutoff Boltzmann equation. Following the rigorous derivation of the lower bound for the non-cutoff Boltzmann collision operator in~\cite{alexandre0}, there have been many works on the regularity for the solutions of the Boltzmann equation in both 
spatially homogeneous (see \cite{alexandre1, al-saf-1, desv1,HMUY, MU,Mo2} and references therein) and
inhomogeneous cases (see \cite{amuxy3, amuxy4,  amuxy5-3}).
In all those works, it was highlighted that the Boltzmann collision operator behaves essentially as a fractional Laplacian $(-\tilde{\Delta}_v)^{\sigma}$  
under the angular singularity assumption on the collision cross-section (see \cite{alexandre0, amuxy3}). We refer the reader to Section~\ref{kkboltz} in appendix for comprehensive explanations about the relevance of this model. In the spatially homogeneous case, this diffusive structure implies a $C^{\infty}$ smoothing effect for the weak solutions to the Cauchy problem constructed in~\cite{Villani} (See \cite{amuxy7, HMUY}). Furthermore, as in the case of the heat equation, the spatially homogeneous Boltzmann equation without angular cutoff enjoys a
smoothing effect in the Gevrey class of order $\sigma$ (see \cite{ L-X, MU, Mo2}). Related to this Gevrey smoothing effect for the spatially homogeneous Boltzmann equation, an  
ultra-analytic smoothing effect was proved in~\cite{Mo1} for both non-linear homogeneous Landau equations and inhomogeneous linear Landau equations.
Regarding the study of the Boltzmann equation in Gevrey spaces, we also refer the reader to the seminal work \cite{Ukai} which establishes the existence and uniqueness in Gevrey classes of a local solution to the Cauchy problem for the Boltzmann equation in both spatially homogeneous and inhomogeneous cases.
Considering now the spatially inhomogeneous Boltzmann equation without angular cutoff, 
the $C^\infty$ hypoellipticity was established in \cite{amuxy3,amuxy4,amuxy5-3} by using the coercivity estimate for proving the regularity with respect to the velocity variable $v$ (see~\cite{alexandre0,amuxy3}) and a version of the uncertainty principle for proving the regularity in the time and space variables $t,x$ via bootstrap arguments (see~\cite{amuxy2}). However, it should be noted that those works do not provide any optimal hypoelliptic estimates for the spatially inhomogeneous Boltzmann equation without angular cutoff. As an attempt to understand further the smoothing effect induced by the Boltzmann collision operator and to relate exactly the structure of the angular singularity in the collision cross-section to this regularizing effect, the present work studies the hypoelliptic properties of the simplified linear model (\ref{yo2}) and aims at giving insights on the hypoelliptic properties what may be expected for the general spatially inhomogeneous Boltzmann equation without angular cutoff.

\section{A key hypoelliptic estimate}\label{section2}

Theorem~\ref{TTH1} will be derived from the following key hypoelliptic estimate:

\bigskip

\begin{proposition}\label{th1}
Let $P$ be the operator defined in~\emph{(\ref{yo2})} and $T>0$ be a positive constant. Then, there exists a positive constant $C_T>0$ such that for all $u \in \mathscr{S}(\rr_{t,x,v}^{2n+1})$ satisfying
\begin{equation}\label{hypo0}
\emph{\textrm{supp }} u(\cdot,x,v) \subset [-T,T], \ (x,v) \in \rr^{2n},
\end{equation}
we have
\begin{equation}\label{hypo}
\big\|(1+|D_x|^{\frac{2\sigma}{2\sigma+1}}+|D_{v}|^{2\sigma})u\big\|_{L^2(\rr^{2n+1})} \leq C_T \big(\|Pu\|_{L^2(\rr^{2n+1})}+\|u\|_{L^2(\rr^{2n+1})}\big).
\end{equation}
\end{proposition}

\bigskip

\noindent

In the following, we use standard notations for symbol classes, see~\cite{hormander} (Chapter~18) or \cite{birkhauser}.
The symbol class $S(m,\Gamma)$ associated to the order function $m$ and metric
$$\Gamma=\frac{dx^2}{\varphi(x,\xi)^2}+\frac{d\xi^2}{\Phi(x,\xi)^2},$$
with $\varphi$, $\Phi$ given positive functions,
stands for the set of functions $a\in C^{\infty}(\real_{x,\xi}^{2n},\cc)$ satisfying for all $\alpha \in \nn^{n}$ and $\beta \in \nn^n$, 
\begin{equation}\label{ro1}
\exists C_{\alpha, \beta}>0, \forall (x,\xi) \in \rr^{2n}, \  |\partial_{x}^{\alpha}\partial_{\xi}^{\beta} a(x,\xi)| \leq C_{\alpha, \beta} m(x,\xi) \varphi(x,\xi)^{-|\alpha|}\Phi(x,\xi)^{-|\beta|}.
\end{equation}

\subsection{Some symbol reductions} We begin by few symbol reductions in order to reduce the symbol of the operator to a convenient normal form. For convenience only, we shall use the Weyl quantization rather than the standard one.
Notice that it will be sufficient to prove the hypoelliptic estimate (\ref{hypo}) for the operator $p^w$ defined by the Weyl quantization of the symbol 
\begin{equation}\label{eq1}
p(t,x,v;\tau,\xi,\eta)=2\pi i\tau+2\pi iv\cdot \xi+a(t,x,v)F(2\pi\eta),
\end{equation}
with $\tau$, $\xi$, $\eta$ respectively standing for the dual variables of the variables $t$, $x$, $v$ and $F$ being the function defined in (\ref{eq0}),
while using the following normalization for the Weyl quantization
$$(a^wu)(x)=\int_{\rr^{2n}}e^{2i\pi(x-y)\cdot \xi}a\big(\frac{x+y}{2},\xi\big)u(y)dyd\xi.$$
Indeed, symbolic calculus shows that the operator
$$R=P-p^w+\frac{1}{2i}\nabla_va(t,x,v) \cdot (\nabla F)(2\pi D_v),$$
is bounded on $L^2$. This is a direct consequence of the $L^2$ continuity theorem in the class $S_{00}^0$ (case $m=\varphi=\Phi=1$ in (\ref{ro1})) after noticing that the Weyl symbol of the operator $R$ together with all its derivatives of any order are bounded on $\rr^{4n+2}$, since all the derivatives of order greater or equal to 2 of the symbol $F$ are bounded given the choice of the positive parameter $0 < \sigma <1$. We refer the reader to formula (2.1.26) in \cite{birkhauser} for explicit constants in the composition formula with the normalization of the Weyl quantization chosen here. It then remains to notice that for any $\eps>0$, there is a positive constant $C_{\eps}>0$ such that 
\begin{multline*}
\big\|\nabla_va(t,x,v) \cdot (\nabla F)(2\pi D_v)u\big\|_{L^2(\rr^{2n+1})} \lesssim \big\|(\nabla F)(2\pi D_v)u\big\|_{L^2(\rr^{2n+1})}\\ \lesssim \|(1+|D_v|^{2\sigma-1})u\|_{L^2(\rr^{2n+1})}  \leq \eps  \|(1+|D_v|^{2\sigma})u\|_{L^2(\rr^{2n+1})} +C_{\eps}\|u\|_{L^2(\rr^{2n+1})},
\end{multline*}   
since the function $\nabla_va$ is bounded on $\rr^{2n+1}$. When studying the operator $p^w$, it is convenient to work on the Fourier side in the variables $x,v$. This is equivalent to study the operator defined by the Weyl quantization of the symbol
$$2\pi i\tau-2\pi i x \cdot \eta +a(t,-\xi,-\eta)F(2\pi v).$$
After relabeling the variables and simplifying the notations, we are reduced to study the operator $P=p^w$ defined by the Weyl quantization of the symbol
\begin{equation}\label{eq2}
p(t,x,v;\tau,\xi,\eta)= i\tau- iy \cdot \xi+a(t,\xi,\eta)F(x),
\end{equation}
where $a$ stands for a $C_b^{\infty}(\rr^{2n+1})$ function satisfying 
\begin{equation}\label{eq00.5}
\exists a_0>0, \ \forall (t,\xi,\eta) \in \rr^{2n+1}, \ a(t,\xi,\eta) \geq a_0>0,
\end{equation}
and $F$ is the function
\begin{equation}\label{eq00}
F(x)=|x|^{2\sigma}w(x)+|x|^2\big(1-w(x)\big),
\end{equation}
with $w$ a new function having experienced a small homothetic transformation compared to the function appearing in (\ref{eq0}).
For the sake of simplicity, we shall keep the definition given in (\ref{eq0}) and assume that $w \in C^{\infty}(\rr^n)$, $0 \leq w \leq 1$, $w(\eta)=1$ if $|\eta| \geq 2$, and $w(\eta)=0$ if $|\eta| \leq 1$. 
Using these notations, Proposition~\ref{th1} is equivalent to the proof of the following a priori estimate. For any $T>0$, there exists a positive constant $C_T>0$ such that for all $u \in \mathscr{S}(\rr_{t,x,y}^{2n+1})$ satisfying
\begin{equation}\label{eq3}
\textrm{supp } u(\cdot,x,y) \subset [-T,T], \ (x,y) \in \rr^{2n},
\end{equation}
we have
\begin{equation}\label{hypo2}
\big\|(1+|x|^{2\sigma}+|y|^{\frac{2\sigma}{2\sigma+1}})u\big\|_{L^2(\rr^{2n+1})} \leq C_T \big(\|Pu\|_{L^2(\rr^{2n+1})}+\|u\|_{L^2(\rr^{2n+1})}\big).
\end{equation}
In order to do so, we consider the new variables 
$$(x_1,x_2)=(y,x+ty),$$ 
with $t \in \rr$ fixed. We define $A(x,y)=(y,x+ty)$. Associated to this linear change of variables are the real linear symplectic transformation
$$(x_1,x_2;\xi_1,\xi_2)=\chi(x,y;\xi,\eta)=(A(x,y);(A^{-1})^T(\xi,\eta))=(y,x+ty;\eta-t\xi,\xi),$$
and the two unitary operators on $L^2(\rr^{2n})$
$$(M_tu)(x_1,x_2)=u(x_2-tx_1,x_1); \ (M_t^{-1}u)(x,y)=u(y,x+ty).$$  
Keeping on considering the $t$-variable as a parameter and defining the symbol
$$b_t(x,y;\xi,\eta)=- iy \cdot \xi+a(t,\xi,\eta)F(x),$$
we have
$$(b_t \circ \chi^{-1})(x_1,x_2;\xi_1,\xi_2)=-ix_1 \cdot \xi_2+a(t,\xi_2,\xi_1+t\xi_2)F(x_2-tx_1),$$
and the symplectic invariance of the Weyl quantization (Theorem~18.5.9 in~\cite{hormander}) or a simple direct calculation
$$M_t^{-1}a^wM_t=(a \circ \chi)^w,$$ 
show that 
$$M_t b_t^w(x,y,D_x,D_y)M_t^{-1}=-ix_1 \cdot D_{x_2}+\big[a(t,\xi_2,\xi_1+t\xi_2)F(x_2-tx_1)\big]^w.$$
We then consider the two unitary operators acting on $L^2(\rr_{t,x,y}^{2n+1})$,
$$(Mu)(t,x_1,x_2)=u(t,x_2-tx_1,x_1), \ (M^{-1}u)(t,x,y)=u(t,y,x+ty),$$ 
and notice that 
$$(MiD_t M^{-1}u)(t,x_1,x_2)=iD_t u(t,x_1,x_2)+ix_1 \cdot D_{x_2}u(t,x_1,x_2).$$
It follows that 
$$iD_t+\big[a(t,\xi_2,\xi_1+t\xi_2)F(x_2-tx_1)\big]^w=MPM^{-1}.$$
Notice also that 
$$1+|x_2-tx_1|^{2\sigma}+|x_1|^{\frac{2\sigma}{2\sigma+1}}=M(1+|x|^{2\sigma}+|y|^{\frac{2\sigma}{2\sigma+1}})M^{-1}.$$
We can therefore reduce the proof of Proposition~\ref{th1} to the proof of the following a priori estimate. For any $T>0$, there exists a positive constant $C_T>0$ such that for all $u \in \mathscr{S}(\rr_{t,x_1,x_2}^{2n+1})$ satisfying
\begin{equation}\label{eq4}
\textrm{supp } u(\cdot,x_1,x_2) \subset [-T,T], \ (x_1,x_2) \in \rr^{2n},
\end{equation}
we have
\begin{equation}\label{eq5}
\big\|(1+|x_2-tx_1|^{2\sigma}+|x_1|^{\frac{2\sigma}{2\sigma+1}})u\big\|_{L^2(\rr^{2n+1})} \leq C_T \big(\|Pu\|_{L^2(\rr^{2n+1})}+\|u\|_{L^2(\rr^{2n+1})}\big),
\end{equation}
for the operator 
\begin{equation}\label{eq6}
P=iD_t+\big[a(t,\xi_2,\xi_1+t\xi_2)F(x_2-tx_1)\big]^w.
\end{equation}

\bigskip

\subsection{Energy estimates via the Wick quantization}
In order to establish the a priori estimate (\ref{eq5}), we shall use some techniques developed in~\cite{cubo} for proving energy estimates via the Wick quantization. We recall in appendix (see Section~\ref{appendix}) the definition and all the main features of the Wick quantization which will be used here.

A first step in adapting this approach is to study the first-order differential operator on $L^2(\rr_t)$,
$$\tilde{P}=iD_t+a(t,\xi_2,\xi_1+t\xi_2)F(x_2-tx_1),$$
where the variables $x_1,x_2,\xi_1,\xi_2$ are considered as parameters, and prove some a priori estimates with respect to those parameters.  
We begin by noticing from (\ref{eq00.5}) that for any $u \in \mathscr{S}(\rr_t)$, 
$$\textrm{Re}(\tilde{P}u,u)_{L^2(\rr_t)}=\int_{\rr}a(t,\xi_2,\xi_1+t\xi_2)F(x_2-tx_1)|u(t)|^2dt 
\geq a_0 \big\|F(x_2-tx_1)^{1/2}u\big\|_{L^2(\rr_t)}^2,$$
since the operator $iD_t$ is skew-adjoint.
It follows from the Cauchy-Schwarz inequality that
\begin{equation}\label{eq7}
a_0 \big\|F(x_2-tx_1)^{1/2}u\big\|_{L^2(\rr_t)}^2 \leq \|\tilde{P}u\|_{L^2(\rr_t)}\|u\|_{L^2(\rr_t)}.
\end{equation}
By denoting $H=\un_{\rr_+}$ the Heaviside function, we may write for any $T \in \rr$,
\begin{multline}\label{eq8}
 2\textrm{Re}(\tilde{P}u,-H(t-T)u)_{L^2(\rr_t)} = 2\textrm{Re}(D_tu,iH(t-T)u)_{L^2(\rr_t)} \\
 -2\int_{\rr}H(t-T)a(t,\xi_2,\xi_1+t\xi_2)F(x_2-tx_1)|u(t)|^2dt,
\end{multline}
and obtain by a simple integration by parts that 
$$2\textrm{Re}(D_tu,iH(t-T)u)_{L^2(\rr_t)}=([D_t,iH(t-T)]u,u)_{L^2(\rr_t)}=\frac{1}{2\pi}|u(T)|^2.$$
Recalling that $a \in L^{\infty}(\rr^{2n+1})$, one may find a positive constant $C_0>0$ such that 
$$\Big|2\int_{\rr}H(t-T)a(t,\xi_2,\xi_1+t\xi_2)F(x_2-tx_1)|u(t)|^2dt\Big| \leq C_0\big\|F(x_2-tx_1)^{1/2}u\big\|_{L^2(\rr_t)}^2.$$
It follows from (\ref{eq7}), (\ref{eq8}) and the Cauchy-Schwarz inequality that there exists a constant $C_1>0$ such that for all $u \in \mathscr{S}(\rr_t)$ and $(x_1,x_2,\xi_1,\xi_2) \in \rr^{4n}$,
\begin{equation}\label{eq9}
\|u\|_{L^{\infty}(\rr_t)}^2 \leq C_1\|\tilde{P}u\|_{L^2(\rr_t)}\|u\|_{L^2(\rr_t)}.
\end{equation}
Following~\cite{cubo}, we split-up the $L^2(\rr_t)$-norm of $u$,
\begin{multline}\label{eq10}
 |x_1|^{\frac{2\sigma}{2\sigma+1}}\|u\|_{L^2(\rr_t)}^2
=  \int_{\{t \in \rr :\ |x_1|^{\frac{2\sigma}{2\sigma+1}} > |x_2-tx_1|^{2\sigma}\}}|x_1|^{\frac{2\sigma}{2\sigma+1}}|u(t)|^2dt \\
+\int_{\{t \in \rr :\ |x_1|^{\frac{2\sigma}{2\sigma+1}} \leq |x_2-tx_1|^{2\sigma}\}}|x_1|^{\frac{2\sigma}{2\sigma+1}} |u(t)|^2dt ,
\end{multline}
and estimate from above the first integral as
\begin{multline*}
\int_{\{t \in \rr :\ |x_1|^{\frac{2\sigma}{2\sigma+1}} > |x_2-tx_1|^{2\sigma}\}}|x_1|^{\frac{2\sigma}{2\sigma+1}}|u(t)|^2dt \\
\leq |x_1|^{\frac{2\sigma}{2\sigma+1}} m\big(\{t \in \rr :\ |x_1|^{\frac{2\sigma}{2\sigma+1}} > |x_2-tx_1|^{2\sigma}\}\big) \|u\|_{L^{\infty}(\rr_t)}^2 \leq 2\|u\|_{L^{\infty}(\rr_t)}^2,
\end{multline*}
with $m$ standing for the Lebesgue measure on $\rr$.
It follows from (\ref{eq9}) that this first integral can be estimated from above as
\begin{equation}\label{eq10.5}
\int_{\{t \in \rr :\ |x_1|^{\frac{2\sigma}{2\sigma+1}} > |x_2-tx_1|^{2\sigma}\}}|x_1|^{\frac{2\sigma}{2\sigma+1}}|u(t)|^2dt \leq 2C_1\|\tilde{P}u\|_{L^2(\rr_t)}\|u\|_{L^2(\rr_t)}.
\end{equation}
While estimating from above the second integral in the right-hand-side of (\ref{eq10}), we may first write that 
\begin{multline*}
\int_{\{t \in \rr :\ |x_1|^{\frac{2\sigma}{2\sigma+1}} \leq |x_2-tx_1|^{2\sigma}\}}|x_1|^{\frac{2\sigma}{2\sigma+1}} |u(t)|^2dt\\ \leq \int_{\{t \in \rr :\ |x_1|^{\frac{2\sigma}{2\sigma+1}} \leq |x_2-tx_1|^{2\sigma}\}}|x_2-tx_1|^{2\sigma} |u(t)|^2dt
\leq \int_{\rr}|x_2-tx_1|^{2\sigma} |u(t)|^2dt,
\end{multline*}
and then notice from (\ref{eq00}) that there exists a positive constant $C_2>0$ such that for all $(t,x_1,x_2,\xi_1,\xi_2) \in \rr^{4n+1}$,
$$|x_2-tx_1|^{2\sigma} \leq F(x_2-tx_1)+C_2.$$
It follows from (\ref{eq7}) that 
\begin{multline}\label{eq11}
\int_{\{t \in \rr :\ |x_1|^{\frac{2\sigma}{2\sigma+1}} \leq |x_2-tx_1|^{2\sigma}\}}|x_1|^{\frac{2\sigma}{2\sigma+1}} |u(t)|^2dt \\ \leq \|F(x_2-tx_1)^{1/2}u\|_{L^2(\rr_t)}^2+C_2\|u\|_{L^2(\rr_t)}^2 \leq
a_0^{-1} \|\tilde{P}u\|_{L^2(\rr_t)}\|u\|_{L^2(\rr_t)}+C_2\|u\|_{L^2(\rr_t)}^2.
\end{multline}
We deduce from (\ref{eq10}), (\ref{eq10.5}) and (\ref{eq11}) that there exist some positive constants $C_3>0$ and $C_4>0$ such that for all $u \in \mathscr{S}(\rr_t)$ and $(x_1,x_2,\xi_1,\xi_2) \in \rr^{4n}$,
$$ |x_1|^{\frac{2\sigma}{2\sigma+1}}\|u\|_{L^2(\rr_t)}^2 \leq C_3 \|\tilde{P}u\|_{L^2(\rr_t)}\|u\|_{L^2(\rr_t)}+C_3\|u\|_{L^2(\rr_t)}^2,$$
that is
$$|x_1|^{\frac{2\sigma}{2\sigma+1}}\|u\|_{L^2(\rr_t)} \leq C_3 \|\tilde{P}u\|_{L^2(\rr_t)}+C_3\|u\|_{L^2(\rr_t)},$$
which implies that 
\begin{equation}\label{eq12}
 \|\langle x_1\rangle^{\frac{2\sigma}{2\sigma+1}}u\|_{L^2(\rr_t)}^2 \leq C_4 \|\tilde{P}u\|_{L^2(\rr_t)}^2+C_4\|u\|_{L^2(\rr_t)}^2,
\end{equation}
with $\langle x \rangle=(1+|x|^2)^{1/2}$.
We shall now prove that there exists a positive constant $C_5>0$ such that for all $u \in \mathscr{S}(\rr_t)$ and $(x_1,x_2,\xi_1,\xi_2) \in \rr^{4n}$,
\begin{equation}\label{eq13}
\|\langle x_1\rangle^{\frac{2\sigma}{2\sigma+1}}u\|_{L^2(\rr_t)}^2+ \||x_2-tx_1|^{2\sigma}u\|_{L^2(\rr_t)}^2 \leq C_5 \|\tilde{P}u\|_{L^2(\rr_t)}^2+C_5\|u\|_{L^2(\rr_t)}^2.
\end{equation}
In order to do so, let $\eps_0 \in \{\pm1\}$ and write the components of the variables $x_1,x_2$ as
$$x_1=(x_{1,1},x_{1,2},...,x_{1,n}) \textrm{ and } x_2=(x_{2,1},x_{2,2},...,x_{2,n}).$$ 
Let $j \in \{1,...,n\}$. We shall first study the case when the real parameter 
\begin{equation}\label{oyo1}
x_{1,j} \neq 0.
\end{equation}
While expanding the following $L^2(\rr_t)$ dot-product where $H$ still denotes the Heaviside function
\begin{align*}
& \ 2\textrm{Re}\big(\tilde{P}u,|x_{2,j}-tx_{1,j}|^{2\sigma}H(\eps_0(x_{2,j}-tx_{1,j})-2\langle x_1 \rangle^{\frac{1}{2\sigma+1}})u\big)_{L^2(\rr_t)} \\
=& \ 2\textrm{Re}\big(D_tu,-i |x_{2,j}-tx_{1,j}|^{2\sigma}H(\eps_0(x_{2,j}-tx_{1,j})-2\langle x_1 \rangle^{\frac{1}{2\sigma+1}})u\big)_{L^2(\rr_t)} \\ 
+ & \  2 \int_{\rr}a(t,\xi_2,\xi_1+t\xi_2)F(x_2-tx_1)|x_{2,j}-t x_{1,j}|^{2\sigma}H(\eps_0(x_{2,j}-tx_{1,j})-2\langle x_1 \rangle^{\frac{1}{2\sigma+1}})|u(t)|^2dt,
\end{align*}
we notice that 
\begin{multline*}
 2\textrm{Re}\big(D_tu,-i |x_{2,j}-tx_{1,j}|^{2\sigma}H(\eps_0(x_{2,j}-tx_{1,j})-2\langle x_1 \rangle^{\frac{1}{2\sigma+1}})u\big)_{L^2(\rr_t)}\\
= \big([D_t,-i |x_{2,j}-tx_{1,j}|^{2\sigma}H(\eps_0(x_{2,j}-tx_{1,j})-2\langle x_1 \rangle^{\frac{1}{2\sigma+1}})]u,u\big)_{L^2(\rr_t)}. 
\end{multline*}
A direct computation gives that 
\begin{align*}
& \ \big([D_t,-i |x_{2,j}-tx_{1,j}|^{2\sigma}H(\eps_0(x_{2,j}-tx_{1,j})-2\langle x_1 \rangle^{\frac{1}{2\sigma+1}})]u,u\big)_{L^2(\rr_t)}\\
=& \ \frac{\eps_0 \sigma}{\pi}\int_{\rr} x_{1,j} |x_{2,j}-tx_{1,j}|^{2\sigma-1}H(\eps_0(x_{2,j}-tx_{1,j})-2\langle x_1 \rangle^{\frac{1}{2\sigma+1}})|u(t)|^2dt\\
& \ + \eps_0\frac{2^{2\sigma-1} }{\pi} \frac{x_{1,j}}{|x_{1,j}|} \langle x_1 \rangle^{\frac{2\sigma}{2\sigma+1}}|u(T)|^2,
\end{align*}
with
$$T=x_{2,j}x_{1,j}^{-1}-2\eps_0 x_{1,j}^{-1}\langle x_1\rangle^{\frac{1}{2\sigma+1}}$$
and we deduce from (\ref{eq9}) that 
\begin{align*}
& \ \big|\big([D_t,-i |x_{2,j}-tx_{1,j}|^{2\sigma}H(\eps_0(x_{2,j}-tx_{1,j})-2\langle x_1 \rangle^{\frac{1}{2\sigma+1}})]u,u\big)_{L^2(\rr_t)}\big|\\
\leq & \ \frac{2^{2\sigma-1}}{\pi}C_1 \|\tilde{P}u\|_{L^2(\rr_t)}\|\langle x_1 \rangle^{\frac{2\sigma}{2\sigma+1}}u\|_{L^2(\rr_t)}\\
+& \ \frac{\sigma}{\pi}\int_{\rr} |x_{1,j}| |x_{2,j}-tx_{1,j}|^{2\sigma-1}H(\eps_0(x_{2,j}-tx_{1,j})-2\langle x_1 \rangle^{\frac{1}{2\sigma+1}})|u(t)|^2dt.
\end{align*}
Since from (\ref{eq00.5}),  
\begin{multline*}
2a_0\|F(x_2-tx_1)^{1/2}|x_{2,j}-t x_{1,j}|^{\sigma}H(\eps_0(x_{2,j}-tx_{1,j})-2\langle x_1 \rangle^{\frac{1}{2\sigma+1}})u\|_{L^2(\rr_t)}^2 \\ \leq 2 \int_{\rr}a(t,\xi_2,\xi_1+t\xi_2)F(x_2-tx_1)|x_{2,j}-t x_{1,j}|^{2\sigma}H(\eps_0(x_{2,j}-tx_{1,j})-2\langle x_1 \rangle^{\frac{1}{2\sigma+1}})|u(t)|^2dt,
\end{multline*}
we deduce from the Cauchy-Schwarz inequality; and all the previous identities and estimates obtained after (\ref{oyo1}) that 
\begin{align*}
& \ 2a_0\|F(x_2-tx_1)^{1/2}|x_{2,j}-t x_{1,j}|^{\sigma}H(\eps_0(x_{2,j}-tx_{1,j})-2\langle x_1 \rangle^{\frac{1}{2\sigma+1}})u\|_{L^2(\rr_t)}^2\\
\leq & \ \frac{2^{2\sigma-1}}{\pi}C_1 \|\tilde{P}u\|_{L^2(\rr_t)}\|\langle x_1 \rangle^{\frac{2\sigma}{2\sigma+1}}u\|_{L^2(\rr_t)}\\
+ & \ \frac{\sigma}{\pi}\int_{\rr} |x_{1,j}| |x_{2,j}-tx_{1,j}|^{2\sigma-1}H(\eps_0(x_{2,j}-tx_{1,j})-2\langle x_1 \rangle^{\frac{1}{2\sigma+1}})|u(t)|^2dt\\
+& \ 2\|\tilde{P}u\|_{L^2(\rr_t)}\||x_{2,j}-tx_{1,j}|^{2\sigma}H(\eps_0(x_{2,j}-tx_{1,j})-2\langle x_1 \rangle^{\frac{1}{2\sigma+1}})u\|_{L^2(\rr_t)}.
\end{align*}
Since from (\ref{eq00}), we have
$$F(x_2-tx_1)=|x_2-t x_1|^{2\sigma} \geq |x_{2,j}-t x_{1,j}|^{2\sigma},$$
on the support of the function $H(\eps_0(x_{2,j}-tx_{1,j})-2\langle x_1 \rangle^{\frac{1}{2\sigma+1}})$, we then deduce from the previous estimate that there exists a positive constant $C_6>0$ such that for all $u \in \mathscr{S}(\rr_t)$ and $(x_1,x_2,\xi_1,\xi_2) \in \rr^{4n}$, $x_{1,j} \neq 0$, 
\begin{multline}\label{eq14}
 C_6^{-1}\||x_{2,j}-t x_{1,j}|^{2\sigma}H(\eps_0(x_{2,j}-tx_{1,j})-2\langle x_1 \rangle^{\frac{1}{2\sigma+1}})u\|_{L^2(\rr_t)}^2
\leq   \|\tilde{P}u\|_{L^2(\rr_t)}^2\\
+ \|\langle x_1 \rangle^{\frac{2\sigma}{2\sigma+1}}u\|_{L^2(\rr_t)}^2
+  \int_{\rr} |x_{1,j}| |x_{2,j}-tx_{1,j}|^{2\sigma-1}H(\eps_0(x_{2,j}-tx_{1,j})-2\langle x_1 \rangle^{\frac{1}{2\sigma+1}})|u(t)|^2dt.
\end{multline}
By using that 
$$|x_{2,j}-tx_{1,j}|^{-1} \leq \frac{1}{2} \langle x_1 \rangle^{-\frac{1}{2\sigma+1}},$$
on the support of the function $H(\eps_0(x_{2,j}-tx_{1,j})-2\langle x_1 \rangle^{\frac{1}{2\sigma+1}})$, one can estimate from above the following integral as
\begin{align*}
& \ \int_{\rr} |x_{1,j}| |x_{2,j}-tx_{1,j}|^{2\sigma-1}H(\eps_0(x_{2,j}-tx_{1,j})-2\langle x_1 \rangle^{\frac{1}{2\sigma+1}})|u(t)|^2dt \\ 
\leq & \ 2^{-1} \int_{\rr} \langle x_1 \rangle^{\frac{2\sigma}{2\sigma+1}}|x_{2,j}-tx_{1,j}|^{2\sigma}H(\eps_0(x_{2,j}-tx_{1,j})-2\langle x_1 \rangle^{\frac{1}{2\sigma+1}})|u(t)|^2dt\\
\leq & \ 2^{-1} \|\langle x_1 \rangle^{\frac{2\sigma}{2\sigma+1}}u\|_{L^2(\rr_t)}\||x_{2,j}-tx_{1,j}|^{2\sigma}H(\eps_0(x_{2,j}-tx_{1,j})-2\langle x_1 \rangle^{\frac{1}{2\sigma+1}})u\|_{L^2(\rr_t)}.
\end{align*}
According to (\ref{eq14}), this implies  that there exists a positive constant $C_7>0$ such that for all $u \in \mathscr{S}(\rr_t)$ and $(x_1,x_2,\xi_1,\xi_2) \in\rr^{4n}$, $x_{1,j} \neq 0$, 
\begin{multline}\label{eq15}
 C_7^{-1}\||x_{2,j}-t x_{1,j}|^{2\sigma}H(\eps_0(x_{2,j}-tx_{1,j})-2\langle x_1 \rangle^{\frac{1}{2\sigma+1}})u\|_{L^2(\rr_t)}^2\\
\leq   \|\tilde{P}u\|_{L^2(\rr_t)}^2
+ \|\langle x_1 \rangle^{\frac{2\sigma}{2\sigma+1}}u\|_{L^2(\rr_t)}^2.
\end{multline}
Finally, since
\begin{multline*}
|x_{2,j}-t x_{1,j}|^{2\sigma} \leq |x_{2,j}-t x_{1,j}|^{2\sigma}H(x_{2,j}-tx_{1,j}-2\langle x_1 \rangle^{\frac{1}{2\sigma+1}})\\
+ |x_{2,j}-t x_{1,j}|^{2\sigma}H(-x_{2,j}+tx_{1,j}-2\langle x_1 \rangle^{\frac{1}{2\sigma+1}})+2^{2\sigma}\langle x_1 \rangle^{\frac{2\sigma}{2\sigma+1}},
\end{multline*}
we deduce from (\ref{eq15}) that there exists a positive constant $C_8>0$ such that for all $u \in \mathscr{S}(\rr_t)$ and $(x_1,x_2,\xi_1,\xi_2) \in \rr^{4n}$, $x_{1,j} \neq 0$, 
\begin{equation}\label{eq16}
\||x_{2,j}-t x_{1,j}|^{2\sigma}u\|_{L^2(\rr_t)}^2\\
\leq  C_8 \|\tilde{P}u\|_{L^2(\rr_t)}^2 + C_8\|\langle x_1 \rangle^{\frac{2\sigma}{2\sigma+1}}u\|_{L^2(\rr_t)}^2,
\end{equation}
which together with (\ref{eq12}) proves the a priori estimate
\begin{equation}\label{eq13.11}
\|\langle x_1\rangle^{\frac{2\sigma}{2\sigma+1}}u\|_{L^2(\rr_t)}^2+ \||x_{2,j}-tx_{1,j}|^{2\sigma}u\|_{L^2(\rr_t)}^2 \leq C_9 \|\tilde{P}u\|_{L^2(\rr_t)}^2+C_9\|u\|_{L^2(\rr_t)}^2,
\end{equation}
with $C_9>0$ a positive constant; in the case when $x_{1,j} \neq 0$. Assume now that $x_{1,j}=0$. In this case, a direct computation using (\ref{eq00.5}) shows that 
\begin{multline*}
 \textrm{Re}(\tilde{P}u,|x_{2,j}|^{2\sigma}u)_{L^2(\rr_t)}=\int_{\rr}a(t,\xi_2,\xi_1+t\xi_2)F(x_2-t x_1)|x_{2,j}|^{2\sigma}|u(t)|^2dt \\
\geq a_0\|F(x_2-t x_1)^{1/2}|x_{2,j}|^{\sigma}u\|_{L^2(\rr_t)}^2,
\end{multline*}
because $iD_t$ is a skew-adjoint operator.
Since 
$$F(x_2-tx_1)=|x_2-t x_1|^{2\sigma} \geq |x_{2,j}|^{2\sigma},$$ 
when $|x_{2,j}| \geq 2$, we first deduce from the Cauchy-Schwarz inequality that for all $u \in \mathscr{S}(\rr_t)$ and $(x_1,x_2,\xi_1,\xi_2) \in \rr^{4n}$, $x_{1,j}=0$, $|x_{2,j}| \geq 2$,
$$ a_0\||x_{2,j}|^{2\sigma}u\|_{L^2(\rr_t)} \leq \|\tilde{P}u\|_{L^2(\rr_t)},$$
which also implies that  
$$ a_0^2\||x_{2,j}|^{2\sigma}u\|_{L^2(\rr_t)}^2 \leq \|\tilde{P}u\|_{L^2(\rr_t)}^2+2^{4\sigma}a_0^2\|u\|_{L^2(\rr_t)}^2.$$
Notice that this second estimate also holds when $|x_{2,j}| \leq 2$. One can then deduce from another use of (\ref{eq12}) that the estimate (\ref{eq13.11}) also holds when $x_{1,j}=0$. This proves the following lemma.

\bigskip

\begin{lemma}\label{prop0}
Consider the first-order differential operator   
$$\tilde{P}=iD_t+a(t,\xi_2,\xi_1+t\xi_2)F(x_2-tx_1),$$
with parameters $(x_1,x_2,\xi_1,\xi_2) \in \rr^{4n}$, and $a$ and $F$ standing for the functions defined in \emph{(\ref{eq00.5})} and \emph{(\ref{eq00})}. Then, there exists a positive constant $C>0$ such that for all $u \in \mathscr{S}(\rr_t)$ and $(x_1,x_2,\xi_1,\xi_2) \in \rr^{4n}$,
\begin{equation}\label{eq17}
\|(1+\langle x_1 \rangle^{\frac{2\sigma}{2\sigma+1}}+\langle x_2-tx_1 \rangle^{2\sigma})u\|_{L^2(\rr_t)} \leq C\big(\|\tilde{P}u\|_{L^2(\rr_t)}+\|u\|_{L^2(\rr_t)}\big).
\end{equation}
\end{lemma}

\bigskip

\subsection{From a priori estimates for the one-dimensional operator with parameters to a priori estimates in Wick quantization}
By applying Lemma~\ref{prop0}
$$\|(1+\langle x_1 \rangle^{\frac{2\sigma}{2\sigma+1}}+\langle x_2-tx_1 \rangle^{2\sigma})u\|_{L^2(\rr_t)}^2 \lesssim \|\tilde{P}u\|_{L^2(\rr_t)}^2+\|u\|_{L^2(\rr_t)}^2,$$
to a function $\Phi(t,x_1,x_2,\xi_1,\xi_2) \in \mathscr{S}(\rr^{4n+1})$ and integrating this a priori estimate with respect to the variables $(x_1,x_2,\xi_1,\xi_2) \in \rr^{4n}$, we obtain that we may find a positive constant $C_1>0$ such that for all $\Phi \in \mathscr{S}(\rr^{4n+1})$,
\begin{equation}\label{eq01}
\|(1+\langle x_1 \rangle^{\frac{2\sigma}{2\sigma+1}}+\langle x_2-tx_1 \rangle^{2\sigma})\Phi\|_{L^2(\rr^{4n+1})}^2 \leq C_1\big(\|\tilde{P}\Phi\|_{L^2(\rr^{4n+1})}^2+\|\Phi\|_{L^2(\rr^{4n+1})}^2\big).
\end{equation}
Let $T>0$ be a positive constant. For any $u \in \mathscr{S}(\rr_{t,x_1,x_2}^{2n+1})$ satisfying
$$\textrm{supp } u(\cdot,x_1,x_2) \subset [-T,T], \ (x_1,x_2) \in \rr^{2n},$$
we shall consider its wave-packets transform in the variables $x_1,x_2$ with parameter $0<\lambda \leq 1$ defined as
$$(W_{\lambda}u)(t,Y)=\big(u(t,\cdot),\varphi_{Y}^{\lambda}\big)_{L^2(\rr_{x_1,x_2}^{2n})} \in \mathscr{S}(\rr^{4n+1}), \ Y=(y_1,y_2;\eta_1,\eta_2) \in \rr^{4n},$$
with 
$$\varphi_{Y}^{\lambda}(x_1,x_2)=(2\lambda)^{\frac{n}{2}}e^{-\pi \lambda (|x_1-y_1|^2+|x_2-y_2|^2)}e^{2i \pi((x_1-y_1)\cdot \eta_1+(x_2-y_2)\cdot \eta_2)}.$$
We apply the a priori estimate 
$$\|(1+\langle x_1 \rangle^{\frac{2\sigma}{2\sigma+1}}+\langle x_2-tx_1 \rangle^{2\sigma})\Phi\|_{L^2(\rr^{4n+1})}^2 \leq C_1\big(\|\tilde{P}\Phi\|_{L^2(\rr^{4n+1})}^2+\|\Phi\|_{L^2(\rr^{4n+1})}^2\big),$$
to the function $$\Phi(t,X)=(W_{\lambda}u)(t,X),$$ 
with $X=(x_1,x_2;\xi_1,\xi_2) \in \rr^{4n}$. By using the fact that the wave-packets transform is an isometric mapping from $L^2(\rr_{x_1,x_2}^{2n})$ to $L^2(\rr_{x_1,x_2,\xi_1,\xi_2}^{4n})$, see appendix (Section~\ref{appendix}), we obtain that 
\begin{multline}\label{eq02}
\|(1+\langle x_1 \rangle^{\frac{2\sigma}{2\sigma+1}}+\langle x_2-tx_1 \rangle^{2\sigma})W_{\lambda}u\|_{L^2(\rr^{4n+1})}^2 \\ \leq C_1\big(\|\tilde{P}W_{\lambda}u\|_{L^2(\rr^{4n+1})}^2+\|u\|_{L^2(\rr^{2n+1})}^2\big).
\end{multline}
We recall from the appendix (Section~\ref{appendix}) that the Wick quantization with parameter $0<\lambda \leq 1$ of a symbol $a$ is formally given by
\begin{equation}\label{eq02.5}
a^{\textrm{Wick}(\lambda)}=W_{\lambda}^*a^{\mu}W_{\lambda},\ 1^{\textrm{Wick}(\lambda)}=W_{\lambda}^*W_{\lambda}=\textrm{Id}_{L^2},
\end{equation}
where $a^{\mu}$ stands for the multiplication operator by the function $a$ on $L^2$, and that the operator 
$$\pi_{\lambda}=W_{\lambda}W_{\lambda}^*,$$ 
is an orthogonal projection on a closed proper subspace of $L^2$. It follows that 
\begin{align*}
& \ \|\pi_{\lambda}(1+\langle x_1 \rangle^{\frac{2\sigma}{2\sigma+1}}
+ \langle x_2-tx_1 \rangle^{2\sigma})W_{\lambda}u\|_{L^2(\rr^{4n+1})}^2 \\
\leq & \  \|(1+\langle x_1 \rangle^{\frac{2\sigma}{2\sigma+1}}+\langle x_2-tx_1 \rangle^{2\sigma})W_{\lambda}u\|_{L^2(\rr^{4n+1})}^2.
\end{align*}
Moreover, since the wave-packets transform is an isometric mapping from $L^2(\rr_{x_1,x_2}^{2n})$ to $L^2(\rr_{x_1,x_2,\xi_1,\xi_2}^{4n})$, we may write that 
\begin{align*}
& \ \|\pi_{\lambda}(1+\langle x_1 \rangle^{\frac{2\sigma}{2\sigma+1}}+\langle x_2-tx_1 \rangle^{2\sigma})W_{\lambda}u\|_{L^2(\rr^{4n+1})}\\
= & \ \|W_{\lambda}W_{\lambda}^*(1+\langle x_1 \rangle^{\frac{2\sigma}{2\sigma+1}}+\langle x_2-tx_1 \rangle^{2\sigma})W_{\lambda}u\|_{L^2(\rr^{4n+1})}\\
= & \ \|(1+\langle x_1 \rangle^{\frac{2\sigma}{2\sigma+1}}+\langle x_2-tx_1 \rangle^{2\sigma})^{\textrm{Wick}(\lambda)}u\|_{L^2(\rr^{2n+1})}.
\end{align*}
It follows from (\ref{eq02}) that 
\begin{align*}
& \ \|(1+\langle x_1 \rangle^{\frac{2\sigma}{2\sigma+1}}+\langle x_2-tx_1 \rangle^{2\sigma})^{\textrm{Wick}(\lambda)}u\|_{L^2(\rr^{2n+1})}^2 \\
& \ +\|(1+\langle x_1 \rangle^{\frac{2\sigma}{2\sigma+1}}+\langle x_2-tx_1 \rangle^{2\sigma})W_{\lambda}u\|_{L^2(\rr^{4n+1})}^2\\
\leq & \  2C_1\big(\|\pi_{\lambda}\tilde{P}W_{\lambda}u\|_{L^2(\rr^{4n+1})}^2+\|(1-\pi_{\lambda})\tilde{P}W_{\lambda}u\|_{L^2(\rr^{4n+1})}^2+\|u\|_{L^2(\rr^{2n+1})}^2\big).
\end{align*}
By using similar arguments as previously, namely (\ref{eq02.5}) and the fact that the wave-packets transform is an isometric mapping from $L^2(\rr_{x_1,x_2}^{2n})$ to $L^2(\rr_{x_1,x_2,\xi_1,\xi_2}^{4n})$, together with recalling that we only consider here the Wick quantization the variables $x_1,x_2$, and not in the $t$-variable, we obtain that 
\begin{align*}
& \ \|\pi_{\lambda}\tilde{P}W_{\lambda}u\|_{L^2(\rr^{4n+1})}\\
= & \ \|W_{\lambda}W_{\lambda}^*\tilde{P}W_{\lambda}u\|_{L^2(\rr^{4n+1})}\\
= & \ \|W_{\lambda}W_{\lambda}^*[iD_t+a(t,\xi_2,\xi_1+t\xi_2)F(x_2-tx_1)]W_{\lambda}u\|_{L^2(\rr^{4n+1})}\\
= & \ \|W_{\lambda}(iD_t+[a(t,\xi_2,\xi_1+t\xi_2)F(x_2-tx_1)]^{\textrm{Wick}(\lambda)})u\|_{L^2(\rr^{4n+1})}\\
= & \ \|iD_tu+[a(t,\xi_2,\xi_1+t\xi_2)F(x_2-tx_1)]^{\textrm{Wick}(\lambda)}u\|_{L^2(\rr^{2n+1})}.
\end{align*}
On the other hand, notice that 
\begin{multline*}
(1-\pi_{\lambda})\tilde{P}W_{\lambda}=(1-\pi_{\lambda})[iD_t+a(t,\xi_2,\xi_1+t\xi_2)F(x_2-tx_1)]W_{\lambda}\\
=(1-\pi_{\lambda})a(t,\xi_2,\xi_1+t\xi_2)F(x_2-tx_1)W_{\lambda},
\end{multline*}
because 
$$(1-\pi_{\lambda})iD_tW_{\lambda}=(1-W_{\lambda}W_{\lambda}^*)W_{\lambda}iD_t=W_{\lambda}(1-W_{\lambda}^*W_{\lambda})iD_t=0,$$
since $D_t$ commutes with the wave-packets transform in the variables $x_1,x_2$, and that $W_{\lambda}^*W_{\lambda}=\textrm{Id}_{L^2}$.
We may then write that 
\begin{multline*}
(1-\pi_{\lambda})a(t,\xi_2,\xi_1+t\xi_2)F(x_2-tx_1)W_{\lambda}=a(t,\xi_2,\xi_1+t\xi_2)F(x_2-tx_1)W_{\lambda}\\
-[\pi_{\lambda},a(t,\xi_2,\xi_1+t\xi_2)F(x_2-tx_1)]W_{\lambda}-a(t,\xi_2,\xi_1+t\xi_2)F(x_2-tx_1)\pi_{\lambda}W_{\lambda}.
\end{multline*}
Since from (\ref{eq02.5}), 
$$\pi_{\lambda}W_{\lambda}=W_{\lambda}W_{\lambda}^*W_{\lambda}=W_{\lambda},$$
it follows that 
$$(1-\pi_{\lambda})a(t,\xi_2,\xi_1+t\xi_2)F(x_2-tx_1)W_{\lambda}=-[\pi_{\lambda},a(t,\xi_2,\xi_1+t\xi_2)F(x_2-tx_1)]W_{\lambda}.$$
We then deduce that for all $u \in \mathscr{S}(\rr_{t,x_1,x_2}^{2n+1})$ satisfying
$$\textrm{supp } u(\cdot,x_1,x_2) \subset [-T,T], \ (x_1,x_2) \in \rr^{2n},$$
we have
\begin{multline*}
 (2C_1)^{-1}\|(1+\langle x_1 \rangle^{\frac{2\sigma}{2\sigma+1}}+\langle x_2-tx_1 \rangle^{2\sigma})^{\textrm{Wick}(\lambda)}u\|_{L^2(\rr^{2n+1})}^2 \\
  + (2C_1)^{-1}\|(1+\langle x_1 \rangle^{\frac{2\sigma}{2\sigma+1}}+\langle x_2-tx_1 \rangle^{2\sigma})W_{\lambda}u\|_{L^2(\rr^{4n+1})}^2\\
\leq  \|iD_tu+[a(t,\xi_2,\xi_1+t\xi_2)F(x_2-tx_1)]^{\textrm{Wick}(\lambda)}u\|_{L^2(\rr^{2n+1})}^2\\
 + \|[\pi_{\lambda},a(t,\xi_2,\xi_1+t\xi_2)F(x_2-tx_1)]W_{\lambda}u\|_{L^2(\rr^{4n+1})}^2+\|u\|_{L^2(\rr^{2n+1})}^2.
\end{multline*}
We now need to study the commutator term
\begin{multline*}
[\pi_{\lambda},a(t,\xi_2,\xi_1+t\xi_2)F(x_2-tx_1)]=[\pi_{\lambda},a(t,\xi_2,\xi_1+t\xi_2)]F(x_2-tx_1)\\ +a(t,\xi_2,\xi_1+t\xi_2)[\pi_{\lambda},F(x_2-tx_1)].
\end{multline*}
In order to do so, we recall from (\ref{yo1}) in appendix that the kernel of the orthogonal projection $\pi_{\lambda}$ is given by 
\begin{equation}\label{EQ1}
e^{-\frac{\pi}{2}\Gamma_{\lambda}(X-Y)}e^{i\pi (x-y)\cdot(\xi+ \eta)},
\end{equation} 
with 
\begin{equation}\label{EQ2}
\Gamma_{\lambda}(X)=\lambda|x|^2+\frac{|\xi|^2}{\lambda}, \ X=(x,\xi) \in \rr^{4n}; \ x=(x_1,x_2),\ \xi=(\xi_1,\xi_2) \in \rr^{2n}.
\end{equation}

\bigskip

\begin{lemma}\label{LEM1}
Let $T>0$ be a positive constant. Then, there exists a positive constant $C>0$ such that for all $0<\lambda \leq 1$ and $u \in \mathscr{S}(\rr_{t,x_1,x_2}^{2n+1})$ satisfying 
$$\emph{\textrm{supp }} u(\cdot,x_1,x_2) \subset [-T,T], \ (x_1,x_2) \in \rr^{2n},$$
we have
\begin{multline*}
\|[\pi_{\lambda},a(t,\xi_2,\xi_1+t\xi_2)]F(x_2-tx_1)W_{\lambda}u\|_{L^2(\rr^{4n+1})} \\
\leq C\lambda^{1/2}\|\langle x_2-tx_1\rangle^{2\sigma}W_{\lambda}u\|_{L^2(\rr^{4n+1})}+C\|u\|_{L^2(\rr^{2n+1})}.
\end{multline*}
\end{lemma}

\bigskip

\noindent
\textit{Proof of Lemma~\ref{LEM1}}. Notice from (\ref{EQ1}) that the kernel of the commutator 
$$[\pi_{\lambda},a(t,\xi_2,\xi_1+t\xi_2)],$$ 
is given by
$$K_{t,\lambda}(X,Y)=e^{-\frac{\pi}{2}\Gamma_{\lambda}(X-Y)}e^{i\pi (x-y)\cdot(\xi+ \eta)}\big(a(t,\eta_2,\eta_1+t\eta_2)-a(t,\xi_2,\xi_1+t\xi_2)\big).$$
Recalling that $a$ is a $C_b^{\infty}(\rr^{2n+1})$ function and therefore a Lipschitz function, we may therefore find a positive constant $C_2>0$ such that for all $t \in [-T,T]$, $0<\lambda \leq 1$ and $X,Y \in \rr^{4n}$,
\begin{equation}\label{E0}
|K_{t,\lambda}(X,Y)| \leq C_2e^{-\frac{\pi}{2}\Gamma_{\lambda}(X-Y)}|\eta-\xi|.
\end{equation}
Notice that 
\begin{multline*}
\int_{\rr^{4n}}|K_{t,\lambda}(X,Y)|dY \leq C_2\int_{\rr^{4n}}e^{-\frac{\pi}{2}\Gamma_{\lambda}(X-Y)}|\eta-\xi|dY=C_2\int_{\rr^{4n}}e^{-\frac{\pi}{2}\Gamma_{\lambda}(Y)}|\eta|dY\\
=C_2\lambda^{1/2}\int_{\rr^{4n}}e^{-\frac{\pi}{2}|Y|^2}|\eta|dY=C_3\lambda^{1/2},
\end{multline*}
with $C_3>0$ a positive constant. 
By symmetry of the estimate (\ref{E0}), we also have 
$$\int_{\rr^{4n}}|K_{t,\lambda}(X,Y)|dX \leq C_3\lambda^{1/2}.$$
Schur test for integral operators together with (\ref{eq00}) show that there exists a positive constant $C_4>0$ such that for all $0<\lambda \leq 1$ and $u \in \mathscr{S}(\rr_{t,x_1,x_2}^{2n+1})$ satisfying
$$\textrm{supp } u(\cdot,x_1,x_2) \subset [-T,T], \ (x_1,x_2) \in \rr^{2n},$$
we have
\begin{multline*}
\|[\pi_{\lambda},a(t,\xi_2,\xi_1+t\xi_2)]F(x_2-tx_1)W_{\lambda}u\|_{L^2(\rr^{4n+1})} \leq C_3\lambda^{1/2}\|F(x_2-tx_1)W_{\lambda}u\|_{L^2(\rr^{4n+1})}\\
\leq C_4\lambda^{1/2}\|\langle x_2-tx_1\rangle^{2\sigma}W_{\lambda}u\|_{L^2(\rr^{4n+1})}+C_4\|u\|_{L^2(\rr^{2n+1})},
\end{multline*}
since we recall that for any $0<\lambda \leq 1$,
$$\|W_{\lambda}u\|_{L^2(\rr^{4n+1})}=\|u\|_{L^2(\rr^{2n+1})}. \ \Box$$

\bigskip

\begin{lemma}\label{LEM2}
Let $T>0$ be a positive constant. Then, there exists a positive constant $C>0$ such that for all $0<\lambda \leq 1$ and $u \in \mathscr{S}(\rr_{t,x_1,x_2}^{2n+1})$ satisfying 
$$\emph{\textrm{supp }} u(\cdot,x_1,x_2) \subset [-T,T], \ (x_1,x_2) \in \rr^{2n},$$
we have
\begin{multline*}
\|a(t,\xi_2,\xi_1+t\xi_2)[\pi_{\lambda},F(x_2-tx_1)]W_{\lambda}u\|_{L^2(\rr^{4n+1})} \\
\leq C \lambda^{-\frac{1+(2\sigma-1)_+}{2}} \|\langle x_2-tx_1\rangle^{(2\sigma-1)_+}W_{\lambda}u\|_{L^2(\rr^{4n+1})},
\end{multline*}
with $(2\sigma-1)_+=\emph{\textrm{max}}(2\sigma-1,0)$.
\end{lemma}

\bigskip

\noindent
\textit{Proof of Lemma~\ref{LEM2}}. We first notice that 
$$\|a(t,\xi_2,\xi_1+t\xi_2)[\pi_{\lambda},F(x_2-tx_1)]W_{\lambda}u\|_{L^2(\rr^{4n+1})} \lesssim \|[\pi_{\lambda},F(x_2-tx_1)]W_{\lambda}u\|_{L^2(\rr^{4n+1})},$$
since $a$ is a $L^{\infty}(\rr^{2n+1})$ function. Arguing as in previous lemma, we notice that the kernel of the commutator $[\pi_{\lambda},F(x_2-tx_1)]$ 
is given by
$$K_{t,\lambda}(X,Y)=e^{-\frac{\pi}{2}\Gamma_{\lambda}(X-Y)}e^{i\pi (x-y)\cdot(\xi+ \eta)}\big(F(y_2-ty_1)-F(x_2-tx_1)\big).$$
Recall from (\ref{eq00}) that there exists a positive constant $C_5>0$ such that for all $x \in \rr^n$,
\begin{equation}\label{eq271}
|\nabla F(x)| \leq C_5 \langle x \rangle^{(2\sigma-1)_+}
\end{equation}
with $(2\sigma-1)_+=\textrm{max}(2\sigma -1,0)$. Writing 
\begin{multline*}
F(y_2-ty_1)-F(x_2-tx_1)\\ =\big(y_2-x_2-t(y_1-x_1)\big) \cdot \int_0^1\nabla F\big((1-\theta)(x_2-y_2-t(x_1-y_1))+y_2-ty_1)\big)d\theta,
\end{multline*}
we have for all $(t,x_1,x_2,y_1,y_2) \in [-T,T] \times \rr^{4n}$,
\begin{align*}
& \ |F(y_2-ty_1)-F(x_2-tx_1)|\\
 \lesssim & \ |x-y| \int_0^1\langle (1-\theta)(x_2-y_2-t(x_1-y_1))+y_2-ty_1)\rangle^{(2\sigma-1)_+} d\theta\\
\lesssim & \ |x-y| \int_0^1  \langle (1-\theta)(x_2-y_2-t(x_1-y_1))\rangle^{(2\sigma-1)_+}  \langle y_2-ty_1\rangle^{(2\sigma-1)_+} d\theta\\
\lesssim & \ |x-y| \langle y_2-ty_1\rangle^{(2\sigma-1)_+} \langle x_2-y_2-t(x_1-y_1)\rangle^{(2\sigma-1)_+}.
\end{align*}
Setting
$$\tilde{K}_{t,\lambda}(X,Y)=\langle y_2-ty_1\rangle^{-(2\sigma-1)_+}K_{t,\lambda}(X,Y),$$
we have 
\begin{equation}\label{E2}
|\tilde{K}_{t,\lambda}(X,Y)| \lesssim e^{-\frac{\pi}{2}\Gamma_{\lambda}(X-Y)} |x-y|\langle x_2-y_2-t(x_1-y_1)\rangle^{(2\sigma-1)_+}.
\end{equation}
While using a change of variables, we notice that for all $t \in [-T,T]$ and $0<\lambda \leq 1$,
\begin{align*}
\int_{\rr^{4n}}|\tilde{K}_{t,\lambda}(X,Y)|dY \lesssim & \ \int_{\rr^{4n}}e^{-\frac{\pi}{2}\Gamma_{\lambda}(X-Y)}|x-y|\langle x_2-y_2-t(x_1-y_1)\rangle^{(2\sigma-1)_+}dY\\
\lesssim & \ \int_{\rr^{4n}}e^{-\frac{\pi}{2}\Gamma_{\lambda}(Y)}|y|\langle y_2-ty_1\rangle^{(2\sigma-1)_+}dY \\
\lesssim & \ \int_{\rr^{4n}}e^{-\frac{\pi}{2}|Y|^2}|\lambda^{-1/2}y|\langle \lambda^{-1/2}(y_2-ty_1)\rangle^{(2\sigma-1)_+}dY \lesssim \lambda^{-\frac{1+(2\sigma-1)_+}{2}},
\end{align*}
since we have $\langle \mu x \rangle \leq \mu \langle x \rangle$, when $\mu \geq 1$. By symmetry of the estimate (\ref{E2}), we also have 
$$\int_{\rr^{4n}}|\tilde{K}_{t,\lambda}(X,Y)|dX \lesssim \lambda^{-\frac{1+(2\sigma-1)_+}{2}}.$$
Schur test for integral operators then shows that for all $0<\lambda \leq 1$ and $u \in \mathscr{S}(\rr_{t,x_1,x_2}^{2n+1})$ satisfying 
$$\textrm{supp } u(\cdot,x_1,x_2) \subset [-T,T], \ (x_1,x_2) \in \rr^{2n},$$
we have
$$\|[\pi_{\lambda},F(x_2-tx_1)]W_{\lambda}u\|_{L^2(\rr^{4n+1})} \lesssim \lambda^{-\frac{1+(2\sigma-1)_+}{2}} \|\langle x_2-tx_1\rangle^{(2\sigma-1)_+}W_{\lambda}u\|_{L^2(\rr^{4n+1})},$$
which proves Lemma~\ref{LEM2}.~$\Box$

\bigskip

\noindent
We then deduce from Lemmas~\ref{LEM1} and \ref{LEM2} that there exists a positive constant $C_6>0$ such that for all $0<\lambda \leq 1$ and $u \in \mathscr{S}(\rr_{t,x_1,x_2}^{2n+1})$ satisfying
$$\textrm{supp } u(\cdot,x_1,x_2) \subset [-T,T], \ (x_1,x_2) \in \rr^{2n},$$
we have
\begin{align*}
& \ C_6^{-1}\|(1+\langle x_1 \rangle^{\frac{2\sigma}{2\sigma+1}}+\langle x_2-tx_1 \rangle^{2\sigma})^{\textrm{Wick}(\lambda)}u\|_{L^2(\rr^{2n+1})}^2 \\
+ & \ C_6^{-1}\|(1+\langle x_1 \rangle^{\frac{2\sigma}{2\sigma+1}}+\langle x_2-tx_1 \rangle^{2\sigma})W_{\lambda}u\|_{L^2(\rr^{4n+1})}^2\\
\leq & \  \|iD_tu+[a(t,\xi_2,\xi_1+t\xi_2)F(x_2-tx_1)]^{\textrm{Wick}(\lambda)}u\|_{L^2(\rr^{2n+1})}^2+\|u\|_{L^2(\rr^{2n+1})}^2\\
+& \ \lambda^{-1-(2\sigma-1)_+} \|\langle x_2-tx_1\rangle^{(2\sigma-1)_+}W_{\lambda}u\|_{L^2(\rr^{4n+1})}^2+\lambda\|\langle x_2-tx_1\rangle^{2\sigma}W_{\lambda}u\|_{L^2(\rr^{4n+1})}^2.
\end{align*}
By choosing a positive constant $0<\lambda_0 \leq 1$ such that 
$$\lambda_0 \leq \frac{1}{2C_6},$$
we notice that one may estimate from above the following term as 
$$\lambda\|\langle x_2-tx_1\rangle^{2\sigma}W_{\lambda}u\|_{L^2(\rr^{4n+1})}^2 \leq (2C_6)^{-1}\|(1+\langle x_1 \rangle^{\frac{2\sigma}{2\sigma+1}}+\langle x_2-tx_1 \rangle^{2\sigma})W_{\lambda}u\|_{L^2(\rr^{4n+1})}^2,$$
for all $0<\lambda \leq \lambda_0$. We then notice that for any $\eps>0$, there exists a positive constant $C_{\eps}>0$ such that for all $u \in \mathscr{S}(\rr_{t,x_1,x_2}^{2n+1})$ satisfying
$$\textrm{supp } u(\cdot,x_1,x_2) \subset [-T,T], \ (x_1,x_2) \in \rr^{2n},$$
we have
$$ \|\langle x_2-tx_1\rangle^{(2\sigma-1)_+}W_{\lambda}u\|_{L^2(\rr^{4n+1})} \leq \eps\|\langle x_2-tx_1\rangle^{2\sigma}W_{\lambda}u\|_{L^2(\rr^{4n+1})}+C_{\eps}\|W_{\lambda}u\|_{L^2(\rr^{4n+1})}.$$ By recalling that 
$$\|W_{\lambda}u\|_{L^2(\rr^{4n+1})}=\|u\|_{L^2(\rr^{2n+1})},$$ 
we finally deduce from these estimates that for any $T>0$, there exists a positive constant $c_T>0$ independent of the parameter $0<\lambda \leq \lambda_0$; and a second positive constant $C_T(\lambda)>0$, which may depend on $\lambda$ such that for all $0<\lambda \leq \lambda_0$ and $u \in \mathscr{S}(\rr_{t,x_1,x_2}^{2n+1})$ satisfying 
$$\textrm{supp } u(\cdot,x_1,x_2) \subset [-T,T], \ (x_1,x_2) \in \rr^{2n},$$
we have
\begin{multline}\label{eq18}
\|[\langle x_1\rangle^{\frac{2\sigma}{2\sigma+1}}]^{\textrm{Wick}(\lambda)}u\|_{L^2(\rr^{2n+1})} +\|[\langle x_2-tx_1\rangle^{2\sigma}]^{\textrm{Wick}(\lambda)}u\|_{L^2(\rr^{2n+1})} \\ \leq c_T\|iD_tu+[a(t,\xi_2,\xi_1+t\xi_2)F(x_2-tx_1)]^{\textrm{Wick}(\lambda)}u\|_{L^2(\rr^{2n+1})}+C_T(\lambda)\|u\|_{L^2(\rr^{2n+1})}.
\end{multline}

\bigskip

\subsection{From a priori estimates in Wick quantization to a priori estimates in Weyl quantization}
In the previous section, we established an a priori estimate with symbols quantized in the Wick quantization with parameter $0<\lambda \leq \lambda_0 \leq 1$. We shall need the following lemma to estimate error terms incurred by coming back from Wick quantization with parameter $0<\lambda \leq \lambda_0$ to the Weyl quantization for symbols appearing in the left-hand-side of (\ref{eq18}).

\bigskip

\begin{lemma}\label{lem1}
Let $a \in C^{\infty}(\rr^{2n})$ be a symbol whose derivatives of order $ \geq 2$ are bounded on $\rr^{2n}$. Then, there exists a positive constant $C>0$ depending only on the $L^{\infty}$-norms of a finite number of derivatives of order greater or equal to 2 of the symbol $a$; such that for all $\lambda > 0$ and $u \in \mathscr{S}(\rr^n)$,
$$\|a^{\emph{Wick}(\lambda)}u-a^wu\|_{L^2(\rr^n)} \leq C\Big(\lambda+\frac{1}{\lambda}\Big)\|u\|_{L^2(\rr^n)}.$$
\end{lemma}

\bigskip

\noindent
\textit{Proof of Lemma~\ref{lem1}.} We notice from (\ref{lay1}) and (\ref{lay2}) in appendix that one may write that
$$a^{\textrm{Wick}(\lambda)}=a^w+R_{\lambda}^w,$$
with 
$$R_{\lambda}(X)=\int_0^1\int_{\rr^{2n}}(1-\theta)a''(X+\theta Y)Y^2e^{-2\pi\Gamma_{\lambda}(Y)}2^ndYd\theta, \ X=(x,\xi) \in \rr^{2n},$$
where 
$$\Gamma_{\lambda}(Y)=\lambda |y|^2+\frac{|\eta|^2}{\lambda}, \ Y=(y,\eta) \in \rr^{2n}.$$
It follows that any derivative of the symbol $R_{\lambda}$, 
$$\partial_{X}^{\alpha}R_{\lambda}(X)=\int_0^1\int_{\rr^{2n}}(1-\theta)\partial_X^{\alpha}a''(X+\theta Y)Y^2e^{-2\pi\Gamma_{\lambda}(Y)}2^ndYd\theta,$$
can be estimated from above as
$$\|\partial_{X}^{\alpha}R_{\lambda}\|_{L^{\infty}(\rr^{2n})} \leq \|\partial_X^{\alpha}a''\|_{L^{\infty}(\rr^{2n})}\int_{\rr^{2n}}(|y|^2+|\eta|^2)e^{-2\pi\lambda |y|^2}e^{-2\pi\frac{|\eta|^2}{\lambda}}2^ndyd\eta.$$
Since 
\begin{multline*}
\int_{\rr^{2n}}(|y|^2+|\eta|^2)e^{-2\pi\lambda |y|^2}e^{-2\pi\frac{|\eta|^2}{\lambda}}2^ndyd\eta\\
=\int_{\rr^{2n}}\Big(\frac{|y|^2}{\lambda}+\lambda|\eta|^2\Big)e^{-2\pi (|y|^2+|\eta|^2)}2^ndyd\eta 
=\mathcal{O}\Big(\lambda+\frac{1}{\lambda}\Big),
\end{multline*}
Lemma~\ref{lem1} is then a direct consequence of the $L^2$ continuity theorem in the class $S_{00}^0$.~$\Box$

\bigskip

\noindent
Notice that we may apply Lemma~\ref{lem1} to the two symbols seen as functions of the variables $(x_1,x_2,\xi_1,\xi_2) \in \rr^{4n}$, 
$$\langle x_1\rangle^{\frac{2\sigma}{2\sigma+1}} \textrm{ and } \langle x_2-tx_1\rangle^{2\sigma}\chi_0(t),$$
with $\chi_0 \in C_0^{\infty}(\rr)$, $\chi_0=1$ on $[-T,T]$, since $0<\sigma<1$. We recall that the $t$-variable is seen as a parameter and that  
we only consider here the Wick and Weyl quantizations in the variables $x_1,x_2$. It follows from Lemma~\ref{lem1} and (\ref{eq18}) that for any fixed $T>0$, there exist a positive constant $c_T>0$ independent of the parameter $0<\lambda \leq \lambda_0$, and a second positive constant $C_T(\lambda)>0$, which may depend on $\lambda$ such that for all $0<\lambda \leq \lambda_0$ and $u \in \mathscr{S}(\rr_{t,x_1,x_2}^{2n+1})$ satisfying 
$$\textrm{supp } u(\cdot,x_1,x_2) \subset [-T,T], \ (x_1,x_2) \in \rr^{2n},$$
we have
\begin{multline}\label{eq19}
\big\|(1+|x_2-tx_1|^{2\sigma}+|x_1|^{\frac{2\sigma}{2\sigma+1}})u\big\|_{L^2(\rr^{2n+1})} \\ \leq c_T\|iD_tu+[a(t,\xi_2,\xi_1+t\xi_2)F(x_2-tx_1)]^{\textrm{Wick}(\lambda)}u\|_{L^2(\rr^{2n+1})}+C_T(\lambda)\|u\|_{L^2(\rr^{2n+1})}.
\end{multline}
We shall now establish a priori estimates for error terms incurred by coming back from Wick quantization with parameter $0< \lambda \leq \lambda_0$ to the Weyl quantization for the symbol $$a(t,\xi_2,\xi_1+t\xi_2)F(x_2-tx_1).$$

\bigskip

\begin{lemma}\label{prop1}
Let $T>0$ be a positive constant. Then, there exists a positive constant $C>0$ such that for all $0< \lambda \leq \lambda_0$ and $u \in \mathscr{S}(\rr_{t,x_1,x_2}^{2n+1})$ satisfying 
$$\emph{\textrm{supp }} u(\cdot,x_1,x_2) \subset [-T,T], \ (x_1,x_2) \in \rr^{2n},$$
we have
\begin{multline*}
C^{-1}\|[a(t,\xi_2,\xi_1+t\xi_2)F(x_2-tx_1)]^{\emph{\textrm{Wick}}(\lambda)}u-[a(t,\xi_2,\xi_1+t\xi_2)F(x_2-tx_1)]^wu\|_{L^2(\rr^{2n+1})} \\
\leq  \lambda^{1-\sigma} \|\langle x_2-tx_1\rangle^{2\sigma}u\|_{L^2(\rr^{2n+1})}+ \lambda^{-\sigma} \|\langle x_2-tx_1\rangle^{(2\sigma-1)_+}u\|_{L^2(\rr^{2n+1})}+ \lambda^{-1}\|u\|_{L^2(\rr^{2n+1})},
\end{multline*}
with $(2\sigma-1)_+=\emph{\textrm{max}}(2\sigma-1,0)$.
\end{lemma}

\bigskip

\noindent
\textit{Proof of Lemma~\ref{prop1}.} Let $\chi_0$ be a $C_0^{\infty}(\rr)$ function satisfying $\chi_0=1$ on $[-T,T]$. By using again (\ref{lay1}) and (\ref{lay2}) in appendix, one may write that
\begin{multline}\label{eq20}
\chi_0(t)[a(t,\xi_2,\xi_1+t\xi_2)F(x_2-tx_1)]^{\textrm{Wick}(\lambda)}\\ =\chi_0(t)[a(t,\xi_2,\xi_1+t\xi_2)F(x_2-tx_1)]^w +R_{t,\lambda}^w,
\end{multline}
with 
$$R_{t,\lambda}(X)=\int_0^1\int_{\rr^{4n}}(1-\theta)\chi_0(t)r_t''(X+\theta Y)Y^2e^{-2\pi\Gamma_{\lambda}(Y)}2^ndYd\theta, $$
where 
$$r_t(x_1,x_2;\xi_1,\xi_2)=a(t,\xi_2,\xi_1+t\xi_2)F(x_2-tx_1), \  X=(x_1,x_2;\xi_1,\xi_2) \in \rr^{4n},$$
and
$$\Gamma_{\lambda}(Y)=\lambda (|y_1|^2+|y_2|^2)+\frac{1}{\lambda}(|\eta_1|^2+|\eta_2|^2), \ Y=(y_1,y_2;\eta_1,\eta_2) \in \rr^{4n}.$$
Define
\begin{equation}\label{eq23}
\tilde{r}_{1,t}(X,Y,\theta)=\chi_0(t)a\big(t,\xi_2+\theta \eta_2,\xi_1+t\xi_2+\theta(\eta_1+t \eta_2)\big) (\nabla_x^2B_t)(x+\theta y).y^2,
\end{equation}
\begin{equation}\label{eq22}
\tilde{r}_{2,t}(X,Y,\theta)=\chi_0(t)(\nabla_{\xi}A_t)(\xi+\theta \eta).\eta \ (\nabla_xB_t)(x+\theta y).y,
\end{equation}
\begin{equation}\label{eq21}
\tilde{r}_{3,t}(X,Y,\theta)=\chi_0(t)(\nabla_{\xi}^2A_t)(\xi+\theta \eta).\eta^2 \ F(x_2-tx_1+\theta(y_2-ty_1)),
\end{equation}
with
\begin{equation}\label{eq5678}
A_t(\xi)=a(t,\xi_2,\xi_1+t\xi_2), \ \xi=(\xi_1,\xi_2) \textrm{ and } B_t(x)=F(x_2-tx_1), \ x=(x_1,x_2).
\end{equation}
We also define
\begin{equation}\label{eq24}
R_{1,t,\lambda}(X)=\int_0^1\int_{\rr^{4n}}(1-\theta)\tilde{r}_{1,t}(X,Y,\theta)e^{-2\pi\Gamma_{\lambda}(Y)}2^ndYd\theta,
\end{equation}
\begin{equation}\label{eq25}
R_{2,t,\lambda}(X)=\int_0^1\int_{\rr^{4n}}(1-\theta)\tilde{r}_{2,t}(X,Y,\theta)e^{-2\pi\Gamma_{\lambda}(Y)}2^ndYd\theta,
\end{equation}
\begin{equation}\label{eq26}
R_{3,t,\lambda}(X)=\int_0^1\int_{\rr^{4n}}(1-\theta)\tilde{r}_{3,t}(X,Y,\theta)e^{-2\pi\Gamma_{\lambda}(Y)}2^ndYd\theta.
\end{equation}
Recall from (\ref{eq00}) that there exists a positive constant $C_1>0$ such that for all $x \in \rr^n$,
\begin{equation}\label{eq27}
|F(x)| \leq C_1 \langle x \rangle^{2\sigma}, \ \ |F'(x)| \leq C_1 \langle x \rangle^{(2\sigma-1)_+}, \ \ \forall k \geq 2, \ F^{(k)} \in L^{\infty}(\rr^n),
\end{equation}
with $(2\sigma-1)_+=\textrm{max}(2\sigma-1,0)$, since $0<\sigma <1$.

\bigskip

\begin{lemma}\label{lem2}
There exists a positive constant $C>0$ such that 
$$\forall \ 0<\lambda \leq \lambda_0, \forall u \in \mathscr{S}(\rr_{t,x_1,x_2}^{2n+1}),\ \|R_{1,t,\lambda}^wu\|_{L^2(\rr^{2n+1})} \leq C \lambda^{-1}\|u\|_{L^2(\rr^{2n+1})}.$$
\end{lemma}

\bigskip

\noindent
\textit{Proof of Lemma~\ref{lem2}}. According to the  $L^2$ continuity theorem in the class $S_{00}^0$, it is sufficient to prove that for all $\alpha \in \nn^{4n}$, there exists a positive constant $C_{2,\alpha}>0$ such that 
$$\forall t \in \rr, \forall \ 0< \lambda \leq \lambda_0,\ \|\partial_X^{\alpha}R_{1,t,\lambda}\|_{L^{\infty}(\rr_X^{4n})} \leq C_{2,\alpha}\lambda^{-1}.$$
Recalling that $a$ is a $C_b^{\infty}(\rr^{2n+1})$ function, we deduce from (\ref{eq23}), (\ref{eq5678}), (\ref{eq24}), (\ref{eq27}) and a change of variables that for any $\alpha \in \nn^{4n}$, there exist some positive constants $C_{3,\alpha},C_{4,\alpha}>0$ such that
\begin{multline*}
\forall t \in \rr, \forall \ 0<\lambda \leq \lambda_0,\ \|\partial_X^{\alpha}R_{1,t,\lambda}\|_{L^{\infty}(\rr_X^{4n})} \\ \leq C_{3,\alpha}\int_{\rr^{4n}}\big(|y_1|^2+|y_2|^2\big)
e^{-2\pi(\lambda (|y_1|^2+|y_2|^2)+\lambda^{-1}(|\eta_1|^2+|\eta_2|^2))}dy_1dy_2d\eta_1d\eta_2=C_{4,\alpha}\lambda^{-1}.
\end{multline*}

\bigskip

\begin{lemma}\label{lem3}
Let $T>0$ be a positive constant. Then, there exists a positive constant $C>0$ such that for all $0<\lambda \leq \lambda_0$ and $u \in \mathscr{S}(\rr_{t,x_1,x_2}^{2n+1})$ satisfying 
$$\emph{\textrm{supp }} u(\cdot,x_1,x_2) \subset [-T,T], \ (x_1,x_2) \in \rr^{2n},$$
we have
$$\|R_{2,t,\lambda}^wu\|_{L^2(\rr^{2n+1})} \leq C \lambda^{-\sigma} \|\langle x_2-tx_1\rangle^{(2\sigma-1)_+}u\|_{L^2(\rr^{2n+1})},$$
with $(2\sigma-1)_+=\emph{\textrm{max}}(2\sigma-1,0) \geq 0.$ 
\end{lemma}

\bigskip

\noindent
\textit{Proof of Lemma~\ref{lem3}}. We notice from (\ref{eq5678}) and (\ref{eq27}) that 
\begin{multline*}
\big| (\nabla_xB_t)(x+\theta y)\big| \lesssim \langle x_2-tx_1+\theta(y_2-ty_1)\rangle^{(2\sigma -1)_+}\\ \lesssim \langle x_2-tx_1\rangle^{(2\sigma -1)_+}
 \langle \theta(y_2-ty_1)\rangle^{(2\sigma -1)_+}
 \lesssim \langle x_2-tx_1\rangle^{(2\sigma -1)_+} \langle y_2-ty_1\rangle^{(2\sigma -1)_+},
 \end{multline*}
when $0 \leq \theta \leq 1$ and $t \in \textrm{supp }\chi_0$. 
Recalling that $a$ is a $C_b^{\infty}(\rr^{2n+1})$ function, it follows from (\ref{eq22}), (\ref{eq5678}) and (\ref{eq25}) that for all $0<\lambda \leq \lambda_0$,
\begin{align*}
|R_{2,t,\lambda}(X)| \lesssim & \ \chi_0(t) \langle x_2-tx_1\rangle^{(2\sigma -1)_+}\int_{\rr^{4n}}|\eta||y|\langle y_2-ty_1\rangle^{(2\sigma -1)_+}e^{-2\pi\Gamma_{\lambda}(Y)}dY \\
\lesssim & \  \chi_0(t)\langle x_2-tx_1\rangle^{(2\sigma -1)_+}\int_{\rr^{4n}}|\eta||y|\langle \lambda^{-\frac{1}{2}}(y_2-ty_1)\rangle^{(2\sigma -1)_+}e^{-2\pi|Y|^2}dY\\
\lesssim & \  \lambda^{-\frac{(2\sigma -1)_+}{2}}\chi_0(t)\langle x_2-tx_1\rangle^{(2\sigma -1)_+},
\end{align*}
since we have $\langle \mu x \rangle \leq \mu \langle x \rangle$, when $\mu \geq 1$. 
After differentiating (\ref{eq22}) with respect to the variable $X=(x,\xi)$, and using the same kind of estimates, we obtain from (\ref{eq27}) that for all $0<\lambda \leq \lambda_0$,
\begin{align*}
& \ |\partial_{x}^{\alpha}\partial_{\xi}^{\beta}R_{2,t,\lambda}(X)| \\
\lesssim & \ \chi_0(t) \langle x_2-tx_1\rangle^{(2\sigma -|\alpha|-1)_+}\int_{\rr^{4n}}|\eta||y|\langle y_2-ty_1\rangle^{(2\sigma -|\alpha|-1)_+}e^{-2\pi\Gamma_{\lambda}(Y)}dY \\
\lesssim & \  \chi_0(t)\langle x_2-tx_1\rangle^{(2\sigma -|\alpha|-1)_+}\int_{\rr^{4n}}|\eta||y|\langle \lambda^{-\frac{1}{2}}(y_2-ty_1)\rangle^{(2\sigma -|\alpha|-1)_+}e^{-2\pi|Y|^2}dY\\
\lesssim & \  \lambda^{-\frac{(2\sigma -|\alpha|-1)_+}{2}}\chi_0(t)\langle x_2-tx_1\rangle^{(2\sigma -|\alpha|-1)_+},
\end{align*}
with 
$$(2\sigma -|\alpha|-1)_+=\textrm{max}(2\sigma -|\alpha|-1,0).$$ 
It follows from those estimates that the Weyl symbol
$$\tilde{a}_{\lambda}(t,X)=\lambda^{\sigma}R_{2,t,\lambda}(X) \sharp[\chi_0(t) \langle x_2-tx_1\rangle^{-(2\sigma-1)_+}],$$ 
of the operator obtained by composition
$$\lambda^{\sigma}R_{2,t,\lambda}^w [\chi_0(t) \langle x_2-tx_1\rangle^{-(2\sigma-1)_+}]^w,$$
belongs to the class $S(1,dt^2+d\tau^2+dX^2)$ uniformly with respect to the parameter $0<\lambda \leq \lambda_0$. We may then deduce from the  $L^2$ continuity theorem in the class $S_{00}^0$ that for all $0<\lambda \leq \lambda_0$ and $u \in \mathscr{S}(\rr_{t,x_1,x_2}^{2n+1})$ satisfying 
$$\textrm{supp } u(\cdot,x_1,x_2) \subset [-T,T], \ (x_1,x_2) \in \rr^{2n},$$
we have
\begin{align*}
\|R_{2,t,\lambda}^wu\|_{L^2(\rr^{2n+1})} = & \ \|\lambda^{-\sigma}\tilde{a}_{\lambda}^w\langle x_2-tx_1\rangle^{(2\sigma-1)_+}u\|_{L^2(\rr^{2n+1})} \\ \lesssim & \  
\lambda^{-\sigma} \|\langle x_2-tx_1\rangle^{(2\sigma-1)_+}u\|_{L^2(\rr^{2n+1})}.
\end{align*}
This ends the proof of Lemma~\ref{lem3}.~$\Box$

\bigskip

\begin{lemma}\label{lem4}
Let $T>0$ be a positive constant. Then, there exists a positive constant $C>0$ such that for all $0<\lambda \leq \lambda_0$ and $u \in \mathscr{S}(\rr_{t,x_1,x_2}^{2n+1})$ satisfying 
$$\emph{\textrm{supp }} u(\cdot,x_1,x_2) \subset [-T,T], \ (x_1,x_2) \in \rr^{2n},$$
we have
$$\|R_{3,t,\lambda}^wu\|_{L^2(\rr^{2n+1})} \leq C \lambda^{1-\sigma} \|\langle x_2-tx_1\rangle^{2\sigma}u\|_{L^2(\rr^{2n+1})}.$$ 
\end{lemma}

\bigskip

\noindent
\textit{Proof of Lemma~\ref{lem4}}. By using similar kind of estimates as in the previous lemma together with the fact that 
\begin{multline*}
|F(x_2-tx_1+\theta(y_2-ty_1))| \lesssim \langle x_2-tx_1+\theta(y_2-ty_1)\rangle^{2\sigma}\\ \lesssim \langle x_2-tx_1\rangle^{2\sigma}
 \langle \theta(y_2-ty_1)\rangle^{2\sigma}
 \lesssim \langle x_2-tx_1\rangle^{2\sigma} \langle y_2-ty_1\rangle^{2\sigma},
 \end{multline*}
when $0 \leq \theta \leq 1$;
it follows from (\ref{eq21}), (\ref{eq5678}) and (\ref{eq26}) that for all $0<\lambda \leq \lambda_0$,
\begin{align*}
|R_{3,t,\lambda}(X)| \lesssim & \ \chi_0(t) \langle x_2-tx_1\rangle^{2\sigma}\int_{\rr^{4n}}|\eta|^2\langle y_2-ty_1\rangle^{2\sigma}e^{-2\pi\Gamma_{\lambda}(Y)}dY \\
\lesssim & \  \chi_0(t)\langle x_2-tx_1\rangle^{2\sigma}\int_{\rr^{4n}}\lambda|\eta|^2\langle \lambda^{-\frac{1}{2}}(y_2-ty_1)\rangle^{2\sigma}e^{-2\pi|Y|^2}dY\\
\lesssim & \  \lambda^{1-\sigma}\chi_0(t)\langle x_2-tx_1\rangle^{2\sigma}.
\end{align*}
After differentiating (\ref{eq21}) and (\ref{eq26}) with respect to the variable $X=(x,\xi)$ and using similar types of estimates, we obtain from (\ref{eq27}) that for all $0<\lambda \leq \lambda_0$,
\begin{align*}
& \ |\partial_{x}^{\alpha}\partial_{\xi}^{\beta}R_{3,t,\lambda}(X)| \\
\lesssim & \ \chi_0(t) \langle x_2-tx_1\rangle^{(2\sigma -|\alpha|)_+}\int_{\rr^{4n}}|\eta|^2\langle y_2-ty_1\rangle^{(2\sigma -|\alpha|)_+}e^{-2\pi\Gamma_{\lambda}(Y)}dY \\
\lesssim & \  \chi_0(t)\langle x_2-tx_1\rangle^{(2\sigma -|\alpha|)_+}\int_{\rr^{4n}}\lambda|\eta|^2\langle \lambda^{-\frac{1}{2}}(y_2-ty_1)\rangle^{(2\sigma -|\alpha|)_+}e^{-2\pi|Y|^2}dY\\
\lesssim & \  \lambda^{1-\frac{(2\sigma -|\alpha|)_+}{2}}\chi_0(t)\langle x_2-tx_1\rangle^{(2\sigma -|\alpha|)_+}.
\end{align*} 
It follows from those estimates that the Weyl symbol
$$\tilde{a}_{\lambda}(t,X)=\lambda^{\sigma-1}R_{3,t,\lambda}(X) \sharp[\chi_0(t) \langle x_2-tx_1\rangle^{-2\sigma}],$$ 
of the operator obtained by composition
$$\lambda^{\sigma-1}R_{3,t,\lambda}^w [\chi_0(t) \langle x_2-tx_1\rangle^{-2\sigma}]^w,$$
belongs to the class $S(1,dt^2+d\tau^2+dX^2)$ uniformly with respect to the parameter $0<\lambda \leq \lambda_0$. We may then deduce from the  $L^2$ continuity theorem in the class $S_{00}^0$ that for all $0<\lambda \leq \lambda_0$ and $u \in \mathscr{S}(\rr_{t,x_1,x_2}^{2n+1})$ satisfying 
$$\textrm{supp } u(\cdot,x_1,x_2) \subset [-T,T], \ (x_1,x_2) \in \rr^{2n},$$
we have
$$\|R_{3,t,\lambda}^wu\|_{L^2(\rr^{2n+1})} =\|\lambda^{1-\sigma}\tilde{a}_{\lambda}^w\langle x_2-tx_1\rangle^{2\sigma}u\|_{L^2(\rr^{2n+1})} \lesssim 
\lambda^{1-\sigma} \|\langle x_2-tx_1\rangle^{2\sigma}u\|_{L^2(\rr^{2n+1})}.$$
This ends the proof of Lemma~\ref{lem4}.~$\Box$

\bigskip

\noindent
By coming back to (\ref{eq20}), Lemma~\ref{prop1} is then a direct consequence of Lemmas~\ref{lem2}, \ref{lem3} and~\ref{lem4}.~$\Box$

\bigskip

Notice that the power $1-\sigma$ appearing in the right-hand-side of the estimate given by Lemma~\ref{prop1} is positive by assumption, since $0<\sigma<1$.
We deduce from Lemma~\ref{prop1} and (\ref{eq19}) that for any fixed $T>0$, one may choose a new positive parameter $0< \lambda_0 \leq 1$ indexing the Wick quantization sufficiently small so that the term 
$$\lambda^{1-\sigma} \|\langle x_2-tx_1\rangle^{2\sigma}u\|_{L^2(\rr^{2n+1})},$$
appearing in the right-hand-side of the estimate given by Lemma~\ref{prop1}
can be controlled by one half times the left-hand-side terms appearing in the a priori estimate (\ref{eq19}). For any fixed $T>0$,   
there therefore exists a positive constant $c_T>0$ such that for all $u \in \mathscr{S}(\rr_{t,x_1,x_2}^{2n+1})$ satisfying 
$$\textrm{supp } u(\cdot,x_1,x_2) \subset [-T,T], \ (x_1,x_2) \in \rr^{2n},$$
we have
\begin{multline*}
c_T^{-1}\big\|(1+|x_2-tx_1|^{2\sigma}+|x_1|^{\frac{2\sigma}{2\sigma+1}})u\big\|_{L^2(\rr^{2n+1})}  \leq \|\langle x_2-tx_1\rangle^{(2\sigma-1)_+}u\|_{L^2(\rr^{2n+1})}\\
+\|iD_tu+[a(t,\xi_2,\xi_1+t\xi_2)F(x_2-tx_1)]^{w}u\|_{L^2(\rr^{2n+1})}+\|u\|_{L^2(\rr^{2n+1})}.
\end{multline*}
By noticing that, for any $\eps>0$, there exists a positive constant $C_{\eps}>0$ such that for all $(t,x_1,x_2) \in [-T,T] \times \rr^{2n}$,
$$\langle x_2-tx_1\rangle^{(2\sigma-1)_+} \leq \eps \langle x_2-tx_1\rangle^{2\sigma}+C_{\eps},$$ 
we finally obtain that for any fixed $T>0$, there exists a positive constant $c_T>0$ such that for all $u \in \mathscr{S}(\rr_{t,x_1,x_2}^{2n+1})$ satisfying 
$$\textrm{supp } u(\cdot,x_1,x_2) \subset [-T,T], \ (x_1,x_2) \in \rr^{2n},$$
we have
\begin{multline*}
c_T^{-1}\big\|(1+|x_2-tx_1|^{2\sigma}+|x_1|^{\frac{2\sigma}{2\sigma+1}})u\big\|_{L^2(\rr^{2n+1})}  \\ 
\leq \|iD_tu+[a(t,\xi_2,\xi_1+t\xi_2)F(x_2-tx_1)]^wu\|_{L^2(\rr^{2n+1})}+\|u\|_{L^2(\rr^{2n+1})}.
\end{multline*}
This proves the a priori estimate (\ref{eq5}) and ends the proof of Proposition~\ref{th1}.~$\Box$

\section{Proof of Theorem~\ref{TTH1}}\label{section3}

\noindent
The next step in the proof of Theorem~\ref{TTH1} is given by establishing the following hypoelliptic estimate:

\bigskip

\begin{proposition}\label{cor1}
Let $P$ be the operator defined in \emph{(\ref{yo2})} and $K$ a compact subset of~$\rr^{2n+1}$. Then, there exists a positive constant $C_K>0$ such that for any $u \in C_0^{\infty}(K)$,
$$\big\|(1+|D_t|^{\frac{2\sigma}{2\sigma+1}}+|D_x|^{\frac{2\sigma}{2\sigma+1}}+|D_{v}|^{2\sigma})u\big\|_{L^2(\rr^{2n+1})} \leq C_K \big(\|Pu\|_{L^2(\rr^{2n+1})}+\|u\|_{L^2(\rr^{2n+1})}\big).$$
\end{proposition}

\bigskip

\noindent
\textit{Proof of Proposition~\ref{cor1}}. 
Let $K$ be a compact subset of $\rr^{2n+1}$ and denote by $\mathcal{F}_{t,x}$ the partial Fourier transform with respect to the $t,x$ variables. Then, for any $u \in C_0^{\infty}(K)$, we may write that
$$\||D_t|^{\frac{2\sigma}{2\sigma+1}}u\|_{L^2(\rr^{2n+1})} \lesssim \||\tau+v\cdot \xi|^{\frac{2\sigma}{2\sigma+1}}\mathcal{F}_{t,x}u\|_{L^2(\rr^{2n+1})}+\||v\cdot \xi|^{\frac{2\sigma}{2\sigma+1}}\mathcal{F}_{t,x}u\|_{L^2(\rr^{2n+1})}.$$
Notice that there exists a positive constant $C_{K}>0$ such that for all $u \in C_0^{\infty}(K)$,
$$\||v\cdot \xi|^{\frac{2\sigma}{2\sigma+1}}\mathcal{F}_{t,x}u\|_{L^2(\rr^{2n+1})} \leq C_{K}  \||D_x|^{\frac{2\sigma}{2\sigma+1}}u\|_{L^2(\rr^{2n+1})}$$
and that 
\begin{multline*}
\||\tau+v\cdot \xi|^{\frac{2\sigma}{2\sigma+1}}\mathcal{F}_{t,x}u\|_{L^2(\rr^{2n+1})} \lesssim \|(\tau+v\cdot \xi)\mathcal{F}_{t,x}u\|_{L^2(\rr^{2n+1})}+\|u\|_{L^2(\rr^{2n+1})}\\
\lesssim \|Pu\|_{L^2(\rr^{2n+1})}+\|a(t,x,v)(-\tilde{\Delta}_v)^{\sigma}u\|_{L^2(\rr^{2n+1})}+\|u\|_{L^2(\rr^{2n+1})}.
\end{multline*}
Recalling that $a \in L^{\infty}(\rr^{2n+1})$, we obtain that for all $u \in C_0^{\infty}(K)$,
\begin{multline*}
\||D_t|^{\frac{2\sigma}{2\sigma+1}}u\|_{L^2(\rr^{2n+1})} \lesssim  \|Pu\|_{L^2(\rr^{2n+1})}+ \||D_x|^{\frac{2\sigma}{2\sigma+1}}u\|_{L^2(\rr^{2n+1})}\\ +\|(-\tilde{\Delta}_v)^{\sigma}u\|_{L^2(\rr^{2n+1})}+\|u\|_{L^2(\rr^{2n+1})}.
\end{multline*}
Proposition~\ref{cor1} is then a direct consequence of Proposition~\ref{th1}.~$\Box$

\bigskip

By using Proposition~\ref{cor1}, one can now complete the proof of Theorem~\ref{TTH1}. Let $K$ be a compact subset of $\rr^{2n+1}$, $s \in \rr$ and $\chi$ be a  $C_0^{\infty}(\rr^{2n+1})$ function satisfying $\chi=1$ on $K$. Setting
$$\Lambda=(1+|D_t|^2+|D_x|^2+|D_v|^2)^{\frac{1}{2}},$$
we apply Proposition~\ref{cor1} to the function $\chi \Lambda^su$ with $u \in C_0^{\infty}(K)$,
$$\big\|(1+\langle D_t\rangle^{\frac{2\sigma}{2\sigma+1}}+\langle D_x\rangle^{\frac{2\sigma}{2\sigma+1}}+\langle D_{v}\rangle ^{2\sigma})\chi \Lambda^su\big\|_{L^2} \leq C_K \big(\|P\chi \Lambda^su\|_{L^2}+\|\chi \Lambda^su\|_{L^2}\big).$$
Notice that 
$$\chi \Lambda^su= \Lambda^s(\chi u)+[\chi,\Lambda^s]u=\Lambda^su+[\chi,\Lambda^s]u,$$
since $\chi=1$ on $K$ and $u \in C_0^{\infty}(K)$. 
We have 
$$\|\chi \Lambda^su\|_{L^2} \leq \|u\|_{s}.$$
Symbolic calculus shows that the Weyl symbol of the operator $[\chi,\Lambda^s]$ belongs to the class
$S(\lambda^{s-1},\tilde{\Gamma})$ with
$$\lambda=(1+|\tau|^2+|\xi|^2+|\eta|^2)^{\frac{1}{2}} \textrm{ and } \tilde{\Gamma}=dt^2+dx^2+dy^2+\lambda^{-2}(d\tau^2+d\xi^2+d\eta^2).$$
It therefore follows that 
\begin{multline*}
\big|\big\|(1+\langle D_t \rangle^{\frac{2\sigma}{2\sigma+1}}+\langle D_x\rangle^{\frac{2\sigma}{2\sigma+1}}+\langle D_{v}\rangle^{2\sigma})\chi \Lambda^su\big\|_{L^2} 
\\ -\big\|(1+\langle D_t \rangle^{\frac{2\sigma}{2\sigma+1}}+\langle D_x\rangle^{\frac{2\sigma}{2\sigma+1}}+\langle D_{v}\rangle^{2\sigma})u\big\|_{s}\big| 
\lesssim  \|\langle D_{v}\rangle ^{(2\sigma-1)_+}u\|_{s},
\end{multline*}
with $(2\sigma-1)_+=\textrm{max}(2\sigma-1,0)$.
On the other hand, we have  
\begin{align*}
\|P\chi \Lambda^su\|_{L^2} \leq & \ \|P \Lambda^su\|_{L^2}+\|P[\chi,\Lambda^s]u\|_{L^2}\\ \leq & \ \|Pu\|_{s}+ \|[P,\Lambda^s]u\|_{L^2} +\|[P,[\chi,\Lambda^s]]u\|_{L^2}+\|[\chi,\Lambda^s]Pu\|_{L^2}.
\end{align*}
Notice that 
$$[P,\Lambda^s]=[v,\Lambda^s] \cdot \partial_x+[a(t,x,v),\Lambda^s](-\tilde{\Delta}_v)^{\sigma}.$$
Symbolic calculus directly shows that 
$$ \|[P,\Lambda^s]u\|_{L^2} \lesssim \|\langle D_{v}\rangle ^{(2\sigma-1)_+}u\|_{s}$$
and 
$$\|[\chi,\Lambda^s]Pu\|_{L^2}=\|[\chi,\Lambda^s]\Lambda^{-s}\Lambda^sPu\|_{L^2} \lesssim \|Pu\|_{s}.$$
It follows that 
\begin{align*}
& \ \big\|(1+|D_t|^{\frac{2\sigma}{2\sigma+1}}+|D_x|^{\frac{2\sigma}{2\sigma+1}}+|D_{v}|^{2\sigma})u\big\|_{s} \\ \lesssim & \
 \|Pu\|_{s}+\|[P,[\chi,\Lambda^s]]u\|_{L^2}
+\|\langle D_{v}\rangle ^{(2\sigma-1)_+}u\|_{s}.
\end{align*}
It remains to study the double commutator 
\begin{multline*}
[P,[\chi,\Lambda^s]]=[\partial_t,[\chi,\Lambda^s]]+[v \cdot \partial_x,[\chi,\Lambda^s]]+[a(t,x,v)(-\tilde{\Delta}_v)^{\sigma},[\chi,\Lambda^s]]
=[\partial_t,[\chi,\Lambda^s]]\\ +v[\partial_x,[\chi,\Lambda^s]]+[v,[\chi,\Lambda^s]]\partial_x+a(t,x,v)[(-\tilde{\Delta}_v)^{\sigma},[\chi,\Lambda^s]]+[a(t,x,v),[\chi,\Lambda^s]](-\tilde{\Delta}_v)^{\sigma}.
\end{multline*}
By using the fact that $a$ is a $C_b^{\infty}(\rr^{2n+1})$ function and $0<\sigma <1$, symbolic calculus directly shows that 
\begin{multline*}
\|[\partial_t,[\chi,\Lambda^s]]u\|_{L^2}+\|[v,[\chi,\Lambda^s]]\partial_xu\|_{L^2}+\|a(t,x,v)[(-\tilde{\Delta}_v)^{\sigma},[\chi,\Lambda^s]]u\|_{L^2}\\ +\|[a(t,x,v),[\chi,\Lambda^s]](-\tilde{\Delta}_v)^{\sigma}u\|_{L^2} \lesssim \|u\|_{s},
\end{multline*}
since we recall that the Weyl symbol of the operator $[\chi,\Lambda^s]$ belongs to the class $S(\lambda^{s-1},\tilde{\Gamma})$. It follows that 
\begin{multline*}
\|[P,[\chi,\Lambda^s]]u\|_{L^2} \lesssim \|[v,[\partial_x,[\chi,\Lambda^s]]]u\|_{L^2}+\|[\partial_x,[\chi,\Lambda^s]]]vu\|_{L^2}+\|u\|_{s}\\
\lesssim \|vu\|_{s}+\|u\|_{s}\lesssim \|\chi vu\|_{s}+\|u\|_{s} \lesssim \|u\|_{s},
\end{multline*}
since $\chi=1$ on $K$ and $u \in C_0^{\infty}(K)$.
By using that for any $\eps>0$, there exists a positive constant $C_{\eps}>0$ such that
$$\|\langle D_{v}\rangle ^{(2\sigma-1)_+}u\|_{s} \leq \eps \||D_{v}|^{2\sigma}u\|_{s}+C_{\eps}\|u\|_{s},$$
we finally obtain that there exists a positive constant $C_K>0$ such that for any $u \in C_0^{\infty}(K)$,
$$\big\|(1+|D_t|^{\frac{2\sigma}{2\sigma+1}}+|D_x|^{\frac{2\sigma}{2\sigma+1}}+|D_{v}|^{2\sigma})u\big\|_{s} \leq C_K \big(\|Pu\|_{s}+\|u\|_{s}\big).$$
This ends the proof of Theorem~\ref{TTH1}.~$\Box$

\section{Proof of Corollary~\ref{ev1}}\label{section4}

This section is devoted to the proof of Corollary~\ref{ev1}. Let $P$ be the operator defined in (\ref{yo2}) and $N_0 \in \nn$. For simplicity, we shall assume that $(t_0,x_0) =(0,0)$ and consider a function $u \in H^{-N_0}(\rr^{2n+1})$ satisfying 
\begin{equation}\label{sev20}
Pu \in H^s_{\textrm{loc}, (t_0,x_0)}(\rr^{2n+1}_{t,x,v}),
\end{equation}
with $s \geq 0$. Let $\varphi$ be a $C_0^\infty(\rr_{t,x}^{n+1})$ function supported in the set 
$$\{(t,x) \in \rr^{n+1} : |(t,x)|<2\}$$
satisfying $\varphi=1$ in a neighborhood of $(t_0,x_0) =(0,0)$ and 
\begin{equation}\label{ts9}
\varphi(t,x) Pu \in H^s(\rr_{t,x,v}^{2n+1}).
\end{equation}
For convenience only, we shall assume that 
\begin{equation}\label{ts8}
 \varphi(t,x)=1,
\end{equation}
for all $|(t,x)| \leq 3/2$.
We consider $\chi$ a $C^\infty(\rr^{n+1}_{t,x})$ function satisfying $-N_0-2 \le \chi \le  s$, 
$\chi = s$ on  $\{|(t,x)| \leq  1\}$, $\chi = -N_0-2$ when $|(t,x)| \ge 3/2$, and define
\begin{equation}\label{sev6}
M_\delta(t,x, \zeta) =\frac{\langle \zeta \rangle^{\chi(t,x)}}{
(1+ \delta \langle \zeta \rangle)^{N_0+s+2}}, \ \zeta=(\tau,\xi,\eta) \in \rr^{2n+1},
\end{equation}
with $\langle \zeta \rangle = (1+|\zeta|^2)^{1/2}$ and $0 <\delta \leq 1$. We recall that $\tau$, $\xi$ and $\eta$ stand respectively for the dual variables of $t$, $x$ and $v$. It directly follows from this definition that for any $\alpha \in \nn^{n+1}$, $\beta \in \nn^{2n+1}$, there exist some positive constants $C_{\alpha, \beta} >0$ such that for all $0< \delta \leq 1$, $(t,x) \in \rr^{n+1}$ and $\zeta \in \rr^{2n+1}$,
\begin{align}\label{symbol-est}
\big| \partial^\alpha_{t,x} \partial_{\zeta}^\beta
M_\delta(t,x, \zeta) \big|\le C_{\alpha, \beta} 
\big(\log \langle \zeta \rangle\big)^{|\alpha|} \langle \zeta \rangle^{-|\beta|}
M_\delta(t,x, \zeta).
\end{align}
Notice that 
$$M_{\delta}(t,x,\zeta) \leq \langle \zeta \rangle^s.$$
Let $\eps>0$ be a positive constant and $k \geq 0$ a nonnegative integer. The symbol $M_{\delta}$ therefore belongs to the symbol class $S^{s+\eps}$ uniformly with respect to the parameter $0<\delta \leq 1$. We recall that the notation $S^m$, with $m \in \rr$, stands for the symbol class
$$S^{m}=S(\langle \zeta \rangle^m,\Gamma), \quad  \Gamma=dz^2+\frac{d\zeta^2}{\langle \zeta \rangle^{2}},$$
with $z=(t,x,v) \in \rr^{2n+1}$. The symbol $M_{\delta}$ also belongs to the symbol class $\tilde{S}^{s}$,
\begin{equation}\label{tt1}
\tilde{S}^{s}=S(\langle \zeta \rangle^s,\tilde{\Gamma}), \quad  \tilde{\Gamma}=\langle \zeta \rangle dz^2+\frac{d\zeta^2}{\langle \zeta \rangle^{2}},
\end{equation}
uniformly with respect to the parameter $0<\delta \leq 1$.
Furthermore, notice that 
$$M_{\delta}(t,x,\zeta) \leq \delta^{-N_0-s-2} \langle \zeta \rangle^{-N_0-2},$$
when $0<\delta \leq 1$. It follows that for each fixed $0<\delta \leq 1$, the symbol $M_\delta$ also belongs to the class
$S^{-N_0-2+\eps}$. 
More specifically, we deduce from (\ref{sev6}) and the Leibniz formula that for any $\alpha \in \nn^{2n+1}$, $\beta \in \nn^{2n+1}$, we may write 
\begin{equation}\label{symbol-est1}
\partial^\alpha_{z} \partial_{\zeta}^\beta\big[\langle v \rangle^{-k} M_\delta(t,x, \zeta)\big]=a_{\alpha,\beta, \delta}(z,\zeta)\langle v \rangle^{-k}M_\delta(t,x, \zeta),
\end{equation}
with $a_{\alpha,\beta, \delta}$ a symbol belonging to the class $S^{\eps-|\beta|}$
uniformly with respect to the parameter $0<\delta \leq 1$.  
We begin by establishing the following lemma:

\bigskip

\begin{lemma}\label{llop1}
Let $0<\eps<1$, $k \geq 0$, $N \in \nn$ and $A \in S^m$, $m \in \rr$. Then, there exists a symbol $B_{\delta}$ belonging to the class $S^{m-(1-\eps)}$ uniformly with respect to the parameter $0<\delta \leq 1$ such that 
$$\emph{\textrm{Op}}(A\langle v \rangle^{-k}M_{\delta})-A(z,D_z)\langle v \rangle^{-k}M_{\delta}(t,x,D_z)-\emph{\textrm{Op}}(B_{\delta}\langle v \rangle^{-k}M_{\delta}) \in \emph{\textrm{Op}}(S^{-N}),$$
uniformly with respect to the parameter $0<\delta \leq 1$.

\end{lemma}

\bigskip

\noindent
Here $A(z,D_z)$ or $\textrm{Op}(A)$ stands for the standard quantization of the symbol $A(z,\zeta)$,
$$A(z,D_z)u(z)=\int_{\rr^{2n+1}}e^{2\pi i(z-\tilde{z}).\zeta}A(z,\zeta)u(\tilde{z})d\zeta d\tilde{z}, \ z=(t,x,v) \in \rr^{2n+1}.$$

\bigskip

\noindent
\textit{Proof of Lemma~\ref{llop1}}. Since $A \in S^m$ and $M_{\delta} \in S^{s+\eps}$ uniformly with respect to the parameter $0<\delta \leq 1$, symbolic calculus shows that $\langle v \rangle^{-k}M_{\delta} \in S^{s+\eps}$ uniformly with respect to the parameter $0<\delta \leq 1$ and
$$A(z,D_z)\langle v \rangle^{-k}M_{\delta}(t,x,D_z) -\sum_{0\le|\alpha|\leq [m+s+\eps]+N}\frac{1}{\alpha !}\textrm{Op}\big((\partial^\alpha_\zeta A) D_z^\alpha\big(\langle v \rangle^{-k}M_{\delta}\big)\big) \in \textrm{Op}(S^{-N}),$$
with $[m+s+\eps]$ being the integer part of $m+s+\eps \in \rr$. We may therefore write from (\ref{symbol-est1}) that 
$$\textrm{Op}(A\langle v \rangle^{-k}M_{\delta})-A(z,D_z)\langle v \rangle^{-k}M_{\delta}(t,x,D_z)-\textrm{Op}(B_{\delta}\langle v \rangle^{-k}M_{\delta}) \in \textrm{Op}(S^{-N}),$$
with 
$$B_{\delta}(z,\zeta)=-\sum_{1\le|\alpha|\leq [m+s+\eps]+N}\frac{1}{\alpha !}(\partial^\alpha_\zeta A)(z,\zeta)a_{\alpha,0,\delta}(z,\zeta).$$
Since $\partial^\alpha_\zeta A \in S^{m-|\alpha|}$ and $a_{\alpha,0,\delta} \in S^{\eps}$ uniformly with respect to $0<\delta \leq 1$, the symbol 
$B_{\delta}$ belongs to the class $S^{m-(1-\eps)}$ uniformly with respect to $0<\delta \leq 1$. This proves Lemma~\ref{llop1}.~$\Box$

\bigskip
\noindent

\bigskip

\begin{lemma}\label{llop2}
Let $k \geq 0$, $N \in \nn$ and $A \in S^m$, $m \in \rr$. Then, there exists a symbol $A_{\delta}$ belonging to the class $S^{m}$ uniformly with respect to the parameter $0<\delta \leq 1$ such that 
$$\emph{\textrm{Op}}(A\langle v \rangle^{-k}M_{\delta})-A_{\delta}(z,D_z)\langle v \rangle^{-k}M_{\delta}(t,x,D_z) \in \emph{\textrm{Op}}(S^{-N}),$$
uniformly with respect to the parameter $0<\delta \leq 1$.

\end{lemma}

\bigskip

\noindent
\textit{Proof of Lemma~\ref{llop2}}. Let $0<\eps<1$. Recalling that the symbol $\langle v \rangle^{-k}M_{\delta}$  belongs to the class $S^{s+\eps}$ uniformly with respect to the parameter $0<\delta \leq 1$, Lemma~\ref{llop2} follows from $k_0$ iterations of Lemma~\ref{llop1} with $m+s+\eps-k_0(1-\eps)<-N$ on the successive remainder terms $\textrm{Op}(B_{\delta}\langle v \rangle^{-k}M_{\delta})$.~$\Box$

\bigskip

\noindent
We shall now use the following commutator argument:

\bigskip

\begin{lemma}\label{lop1}
Let $0<\eps<1$, $k \geq 1$ and $u$ the function defined in \emph{(\ref{sev20})}. Then, there exists a positive constant $C_{k,\eps}>0$ such that for all $0<\delta \leq 1$,  
$$\|[\partial_t + v\cdot \nabla_x \, , \langle v \rangle^{-k} M_\delta(t,x,D_z)] u\|_{L^2} \leq C_{k,\eps}  \| \langle v \rangle^{-k+1} 
M_\delta(t,x,D_z) u \|_{\eps} +C_{k,\eps}\|u\|_{-N_0}.$$
\end{lemma}

\bigskip

\noindent
\textit{Proof of Lemma~\ref{lop1}}. Let $k \geq 1$ be a positive integer and $0<\eps<1$. We consider the commutator 
$$[\partial_t + v\cdot \nabla_x \, , \langle v \rangle^{-k} M_\delta(t,x,D_z)].$$
Symbolic calculus shows that the symbol of the commutator 
$$ [\partial_t + v\cdot \nabla_x , \langle v \rangle^{-k}M_\delta(t,x,D_z)]=\langle v \rangle^{-k} [\partial_t + v\cdot \nabla_x , M_\delta(t,x,D_z)],$$ 
in the standard quantization is 
$$\langle v \rangle^{-k} \big((\partial_t \chi)(t,x) + v\cdot (\nabla_x\chi)(t,x)\big) M_\delta(t,x, \zeta)\log \langle \zeta \rangle-\langle v \rangle^{-k} \xi \cdot (\nabla_{\eta}M_{\delta})(t,x,\zeta).$$
This symbol may be written as 
\begin{multline}\label{sev21b}
\log \langle \zeta \rangle\Big(\langle v \rangle^{-1}(\partial_t \chi)(t,x) + \frac{v}{\langle v \rangle}\cdot (\nabla_x\chi)(t,x)\Big) \langle v \rangle^{-(k-1)}M_\delta(t,x, \zeta)\\ - \langle v \rangle^{-1}\xi \cdot (\nabla_{\eta})\big[\langle v \rangle^{-(k-1)}M_{\delta}(t,x,\zeta)\big].
\end{multline} 
By using that $\log \langle \zeta \rangle \in S^{\eps}$, it follows from (\ref{symbol-est1}) and (\ref{sev21b}) that this symbol may also be written as
 \begin{equation}\label{sev21}
A_{\delta}(z,\zeta)\langle v \rangle^{-(k-1)}M_{\delta}(t,x,\zeta),
\end{equation} 
with $A_{\delta}$ belonging to the class $S^{\eps}$ uniformly with respect to $0<\delta \leq 1$. Then, we deduce from Lemma~\ref{llop2} that 
\begin{equation}\label{sev1}
[\partial_t + v\cdot \nabla_x ,\langle v \rangle^{-k} M_\delta(t,x,D_z)]-\tilde{A}_{\delta}(z,D_z)\langle v \rangle^{-(k-1)}M_{\delta}(t,x,D_z) \in \textrm{Op}(S^{-N_0}),
\end{equation}
with $\tilde{A}_{\delta}$ a new symbol belonging to the class $S^{\eps}$ uniformly with respect to $0<\delta \leq 1$.
This implies that 
\begin{multline}\label{sev2}
\|[\partial_t + v\cdot \nabla_x , \langle v \rangle^{-k}M_\delta(t,x,D_z)]u-\tilde{A}_{\delta}(z,D_z)\langle v \rangle^{-(k-1)}M_{\delta}(t,x,D_z)u\|_{L^2} \\ \lesssim \|u\|_{-N_0}.
\end{multline}
On the other hand, we have
\begin{equation}\label{rum2}
\|\tilde{A}_{\delta}(z,D_z)\langle v \rangle^{-(k-1)}M_{\delta}(t,x,D_z)u\|_{L^2} \lesssim \|\langle v \rangle^{-(k-1)}M_{\delta}(t,x,D_z)u\|_{\eps},
\end{equation}
since $\tilde{A}_{\delta} \in S^{\eps}$ uniformly with respect to $0<\delta \leq 1$.
It follows from (\ref{sev2}), (\ref{rum2}) and the triangular inequality that the estimate 
$$\|[\partial_t + v\cdot \nabla_x \, , \langle v \rangle^{-k} M_\delta(t,x,D_z)] u\|_{L^2} \lesssim  \| \langle v \rangle^{-k+1} 
M_\delta(t,x,D_z) u \|_{\eps} +\|u\|_{-N_0},$$
holds uniformly with respect to the parameter $0<\delta \leq 1$. This proves Lemma~\ref{lop1}.~$\Box$

\bigskip

\noindent
One can now estimate from above the following second commutator:

\bigskip

\begin{lemma}\label{lop2}
Let $0<\eps<1$, $k \geq 0$ and $u$ be the function defined in \emph{(\ref{sev20})}. Then, there exists a positive constant $C_{k,\eps}>0$ such that for all $0<\delta \leq 1$,  
\begin{multline*}
\|[a(t,x,v)(-\tilde{\Delta}_v)^{\sigma}, \langle v \rangle^{-k} M_\delta(t,x,D_z)] u\|_{L^2}\\ \leq C_{k,\eps} 
 \| \langle D_v\rangle^{(2\sigma-1)_++\eps} \langle v \rangle^{-k} 
M_\delta(t,x,D_z) u \|_{L^2}+C_{k,\eps}\|u\|_{-N_0},
\end{multline*}
with $(2\sigma-1)_+=\max(2\sigma-1,0)$.
\end{lemma}

\bigskip

\noindent
\textit{Proof of Lemma~\ref{lop2}}. Notice that the symbol of the operator $a(t,x,v)(-\tilde{\Delta}_v)^{\sigma}$ in the standard quantization does not belong to the class $S^{2\sigma}$ in $\rr_{z,\zeta}^{4n+2}$. Handling this commutator term 
$$[a(t,x,v)(-\tilde{\Delta}_v)^{\sigma}, \langle v \rangle^{-k} M_\delta(t,x,D_z)]$$
needs therefore some caution. By using that 
$$[a(t,x,v),\langle v \rangle^{-k}]=0 \textrm{ and } [(-\tilde{\Delta}_v)^{\sigma},M_\delta(t,x,D_z)]=0,$$ 
we may write that 
\begin{multline}\label{ryu1}
[a(t,x,v)(-\tilde{\Delta}_v)^{\sigma}, \langle v \rangle^{-k} M_\delta(t,x,D_z)]=a(t,x,v)[(-\tilde{\Delta}_v)^{\sigma}, \langle v \rangle^{-k}]M_\delta(t,x,D_z)\\
+ \langle v \rangle^{-k} [a(t,x,v), M_\delta(t,x,D_z)](-\tilde{\Delta}_v)^{\sigma}.
\end{multline}
Let $0<\eps<1$. Regarding the first term in the right-hand-side of (\ref{ryu1}), we first notice that the symbol of the operator 
$[(-\tilde{\Delta}_v)^{\sigma}, \langle v \rangle^{-k}]\langle v \rangle^{k} \langle D_v \rangle^{-(2\sigma-1)_+-\eps}$ 
belongs to the class
$$S(1,dv^2+\langle \eta \rangle^{-2}d\eta^2).$$
Indeed, the symbol of the operator $(-\tilde{\Delta}_v)^{\sigma}$ belongs to $S(\langle \eta \rangle^{2\sigma},dv^2+\langle \eta \rangle^{-2}d\eta^2)$, the symbol of the operator 
$\langle D_v \rangle^{-(2\sigma-1)_+-\eps}$ belongs to $S(\langle \eta \rangle^{-(2\sigma-1)_+-\eps},dv^2+\langle \eta \rangle^{-2}d\eta^2)$
 and
$\langle v \rangle^{m} \in S(\langle v \rangle^{m},dv^2+\langle \eta \rangle^{-2}d\eta^2)$ when $m \in \rr$. Symbolic calculus shows that the symbol of the commutator 
$[(-\tilde{\Delta}_v)^{\sigma}, \langle v \rangle^{-k}]$ belongs to $S(\langle \eta \rangle^{2\sigma-1}\langle v \rangle^{-k},dv^2+\langle \eta \rangle^{-2}d\eta^2)$ and therefore that the symbol of the operator $[(-\tilde{\Delta}_v)^{\sigma}, \langle v \rangle^{-k}]\langle v \rangle^{k}\langle D_v \rangle^{-(2\sigma-1)_+-\eps}$ 
belongs to the class
$$S(\langle \eta \rangle^{2\sigma-1-(2\sigma-1)_+-\eps},dv^2+\langle \eta \rangle^{-2}d\eta^2) \subset S(1,dv^2+\langle \eta \rangle^{-2}d\eta^2).$$
Recalling that by assumption $a \in C_b^{\infty}(\rr_{t,x,v}^{2n+1})$, we obtain that 
\begin{align}\label{ryu2}
& \ \|a(t,x,v)[(-\tilde{\Delta}_v)^{\sigma}, \langle v \rangle^{-k}]M_\delta(t,x,D_z)u\|_{L^2} \\ \notag
\lesssim & \ \|[(-\tilde{\Delta}_v)^{\sigma}, \langle v \rangle^{-k}] \langle v \rangle^{k}\langle D_v \rangle^{-(2\sigma-1)_+-\eps}\langle D_v \rangle^{(2\sigma-1)_++\eps} \langle v \rangle^{-k}M_\delta(t,x,D_z)u\|_{L^2} \\ \notag
\lesssim & \  \|\langle D_v\rangle^{(2\sigma-1)_++\eps}  \langle v \rangle^{-k}M_\delta(t,x,D_z)u\|_{L^2}.
\end{align}
Recalling that $M_{\delta} \in S^{s+\eps}$ uniformly with respect to $0<\delta \leq 1$ and that by assumption $a \in S^0$, symbolic calculus shows that the symbol of the operator 
$$\langle v \rangle^{-k}[a(t,x,v), M_\delta(t,x,D_z)],$$ 
writes as
\begin{align}\label{ryu3}
& \ -\langle v \rangle^{-k} \sum_{1\le|\alpha|\leq [s+\eps]+N_0+2}\frac{1}{\alpha !}(\partial^\alpha_\zeta M_{\delta}) D_z^\alpha a\\ 
= & \ -\sum_{1\le|\alpha|\leq [s+\eps]+N_0+2}\frac{1}{\alpha !}(D_z^\alpha a)\partial^\alpha_\zeta\big[\langle v \rangle^{-k} M_{\delta}\big],\notag
\end{align}
up to a term in $S^{-N_0-2}$. This latter term in the class $S^{-N_0-2}$ gives rise to a term bounded by 
$$\|(-\tilde{\Delta}_v)^{\sigma}u\|_{-N_0-2} \lesssim \|u\|_{-N_0},$$
because $0<\sigma <1$,
while estimating from above the term 
$$\|\langle v \rangle^{-k} [a(t,x,v), M_\delta(t,x,D_z)](-\tilde{\Delta}_v)^{\sigma}u\|_{L^2}.$$
By coming back to (\ref{ryu3}), we deduce from (\ref{symbol-est1}) that this symbol reduces to a term $A_{\delta}(z,\zeta)\langle v \rangle^{-k} M_{\delta}$ with 
$A_{\delta}$ belonging to $S^{\eps-1}$ uniformly with respect to $0<\delta \leq 1$. By using Lemma~\ref{llop2}, we notice that 
$$\textrm{Op}(A_{\delta}\langle v \rangle^{-k}M_{\delta})-\tilde{A}_{\delta}(z,D_z)\langle v \rangle^{-k}M_{\delta}(t,x,D_z) \in \textrm{Op}(S^{-N_0-2}),$$
with $\tilde{A}_{\delta}$ a new symbol belonging to the class $S^{\eps-1}$ uniformly with respect to $0<\delta \leq 1$.
It follows that 
\begin{multline}\label{ryu4}
\|\textrm{Op}(A_{\delta}\langle v \rangle^{-k}M_{\delta})(-\tilde{\Delta}_v)^{\sigma}u\|_{L^2} \lesssim \|(-\tilde{\Delta}_v)^{\sigma}u\|_{-N_0-2}\\ + \|\tilde{A}_{\delta}(z,D_z)\langle v \rangle^{-k}M_{\delta}(t,x,D_z)(-\tilde{\Delta}_v)^{\sigma}u\|_{L^2}. 
\end{multline}
Since $[M_{\delta}(t,x,D_z),(-\tilde{\Delta}_v)^{\sigma}]=0$ and $\tilde{A}_{\delta} \in S^{\eps-1}$ uniformly with respect to $0<\delta \leq 1$, we deduce from (\ref{ryu4}) that
\begin{multline}\label{ryu5}
\|\textrm{Op}(A_{\delta}\langle v \rangle^{-k}M_{\delta})(-\tilde{\Delta}_v)^{\sigma}u\|_{L^2} \lesssim \|u\|_{-N_0} + \|\langle v \rangle^{-k}M_{\delta}(t,x,D_z)(-\tilde{\Delta}_v)^{\sigma}u\|_{\eps-1}\\
\lesssim \|u\|_{-N_0} + \|\langle v \rangle^{-k}(-\tilde{\Delta}_v)^{\sigma}\langle v \rangle^{k}\langle v \rangle^{-k}M_{\delta}(t,x,D_z)u\|_{\eps-1}. 
\end{multline}
Notice that 
\begin{multline}\label{ryu6}
\|\langle v \rangle^{-k}(-\tilde{\Delta}_v)^{\sigma}\langle v \rangle^{k}\langle v \rangle^{-k}M_{\delta}(t,x,D_z)u\|_{\eps-1} 
\leq   \|(-\tilde{\Delta}_v)^{\sigma}\langle v \rangle^{-k}M_{\delta}(t,x,D_z)u\|_{\eps-1}\\
+\|[\langle v \rangle^{-k},(-\tilde{\Delta}_v)^{\sigma}]\langle v \rangle^{k}\langle v \rangle^{-k}M_{\delta}(t,x,D_z)u\|_{\eps-1}.
\end{multline}
We may write
\begin{align}\label{ryu111}
& \ \|(-\tilde{\Delta}_v)^{\sigma}\langle v \rangle^{-k}M_{\delta}(t,x,D_z)u\|_{\eps-1}  \\ 
\lesssim & \ \|\langle D_{z}\rangle^{-(1-\eps)}(-\tilde{\Delta}_v)^{\sigma}\langle D_v\rangle^{-2\sigma+1-\eps}\langle D_v\rangle^{2\sigma-1+\eps}\langle v \rangle^{-k}M_{\delta}(t,x,D_z)u\|_{L^2} \notag \\
\lesssim & \ \|\langle D_v\rangle^{2\sigma-1+\eps}\langle v \rangle^{-k}M_{\delta}(t,x,D_z)u\|_{L^2}\lesssim  \|\langle D_v\rangle^{(2\sigma-1)_++\eps}\langle v \rangle^{-k}M_{\delta}(t,x,D_z)u\|_{L^2}, \notag
\end{align}
with $\langle D_z \rangle=(1+D_t^2+|D_x|^2+|D_v|^2)^{1/2}$, since the operator 
$$\langle D_{z}\rangle^{-(1-\eps)}(-\tilde{\Delta}_v)^{\sigma}\langle D_v\rangle^{-2\sigma+1-\eps},$$ 
is bounded on $L^2$.
Notice also that
\begin{align}\label{ryu112}
& \  \|[\langle v \rangle^{-k},(-\tilde{\Delta}_v)^{\sigma}]\langle v \rangle^{k}\langle v \rangle^{-k}M_{\delta}(t,x,D_z)u\|_{\eps-1} \\
\lesssim & \ \notag  \|\langle D_{z}\rangle^{-(1-\eps)}[\langle v \rangle^{-k},(-\tilde{\Delta}_v)^{\sigma}]\langle v \rangle^{k}\langle v \rangle^{-k}M_{\delta}(t,x,D_z)u\|_{L^2}\\
\lesssim & \ \notag  \|\langle D_{z}\rangle^{-(1-\eps)}[\langle v \rangle^{-k},(-\tilde{\Delta}_v)^{\sigma}]\langle v \rangle^{k}\langle D_v\rangle^{-2\sigma+1-\eps}\langle D_v\rangle^{2\sigma-1+\eps}\langle v \rangle^{-k}M_{\delta}(t,x,D_z)u\|_{L^2} \\
\lesssim & \ \|\langle D_v\rangle^{2\sigma-1+\eps}\langle v \rangle^{-k}M_{\delta}(t,x,D_z)u\|_{L^2}\lesssim  \|\langle D_v\rangle^{(2\sigma-1)_++\eps}\langle v \rangle^{-k}M_{\delta}(t,x,D_z)u\|_{L^2}. \notag
\end{align}
Indeed, we just need to check that the operator 
\begin{multline*}
\langle D_{z}\rangle^{-(1-\eps)}[\langle v \rangle^{-k},(-\tilde{\Delta}_v)^{\sigma}]\langle v \rangle^{k}\langle D_v\rangle^{-2\sigma+1-\eps}\\
=\langle D_{z}\rangle^{-(1-\eps)}\langle D_v\rangle^{1-\eps} \langle D_v\rangle^{-1+\eps} [\langle v \rangle^{-k},(-\tilde{\Delta}_v)^{\sigma}]\langle v \rangle^{k}\langle D_v\rangle^{-2\sigma+1-\eps},
\end{multline*}
is bounded on $L^2$. This is actually the case since the operator 
$$\langle D_{z}\rangle^{-(1-\eps)}\langle D_v\rangle^{1-\eps},$$
is obviously bounded on $L^2$ since $0<1-\eps<1$, and that symbolic calculus easily shows that the symbol of the operator 
$$\langle D_v\rangle^{-1+\eps} [\langle v \rangle^{-k},(-\tilde{\Delta}_v)^{\sigma}]\langle v \rangle^{k}\langle D_v\rangle^{-2\sigma+1-\eps},$$
belongs to the class $S(1,dv^2+\langle \eta \rangle^{-2}d\eta^2)$. It follows from (\ref{ryu5}), (\ref{ryu6}), (\ref{ryu111}) and (\ref{ryu112}) that 
\begin{align}\label{ryu7}
& \ \|\langle v \rangle^{-k} [a(t,x,v), M_\delta(t,x,D_z)](-\tilde{\Delta}_v)^{\sigma}u\|_{L^2} \\ \lesssim & \ 
\|\langle D_v\rangle^{(2\sigma-1)_++\eps} \langle v \rangle^{-k}M_{\delta}(t,x,D_z)u\|_{L^2}+ \|u\|_{-N_0}. \notag
\end{align}
Together with (\ref{ryu1}) and (\ref{ryu2}), we finally obtain the estimate
\begin{align*}
&\|[a(t,x,v)(-\tilde{\Delta}_v)^{\sigma}, \langle v \rangle^{-k} M_\delta(t,x,D_z)] u\|_{L^2} \\
&\lesssim  \| \langle D_v\rangle^{(2\sigma-1)_++\eps} \langle v \rangle^{-k} 
M_\delta(t,x,D_z) u \|_{L^2} +\|u\|_{-N_0},
\end{align*}
which proves Lemma~\ref{lop2}.~$\Box$

\bigskip

\noindent
Summing up the two previous lemmas provides the following estimate:

\bigskip

\begin{lemma}\label{lop3}
Let $0<\eps<1$, $k \geq 1$, $P$ be the operator defined in \emph{(\ref{yo2})} and $u$ the function defined in \emph{(\ref{sev20})}. Then, there exists a positive constant $C_{k,\eps}>0$ such that for all $0<\delta \leq 1$,  
\begin{multline*}
\|[P, \langle v \rangle^{-k} M_\delta(t,x,D_z)] u\|_{L^2} \leq C_{k,\eps}  \|\langle D_v\rangle^{(2\sigma-1)_++\eps}\langle v \rangle^{-k} 
M_\delta(t,x,D_z) u \|_{L^2} \\ +C_{k,\eps}  \| \langle v \rangle^{-k+1} 
M_\delta(t,x,D_z) u \|_{\eps} +C_{k,\eps}\|u\|_{-N_0},
\end{multline*}
with $(2\sigma-1)_+=\max(2\sigma-1,0)$.
\end{lemma}

\bigskip

\noindent
Resuming our proof of Corollary~\ref{ev1}, we may first deduce from Proposition~\ref{th1} that for all $w \in \mathscr{S}(\rr_{t,x,v}^{2n+1})$ satisfying
\begin{equation}\label{rie1}
\textrm{supp } w(\cdot,x,v) \subset [-3,3], \ (x,v) \in \rr^{2n},
\end{equation}
we have
\begin{equation}\label{rie2}
\|\langle D_x \rangle^{\frac{2\sigma}{2\sigma+1}}w\|_{L^2}+\|\langle D_v \rangle^{2\sigma}w\|_{L^2}  \lesssim \|Pw\|_{L^2}+\|w\|_{L^2}.
\end{equation}
These estimates are maximal hypoelliptic estimates and we therefore get an upper bound for the transport term
\begin{multline}\label{rie3}
\|(\partial_t+ v\cdot \nabla_x)w\|_{L^2} \lesssim \|P w \|_{L^2}+\|a(t,x,v)(-\tilde{\Delta}_v^\sigma)w\|_{L^2}\\ \lesssim \|P w \|_{L^2}+\|\langle D_v\rangle^{2\sigma}w\|_{L^2}
\lesssim \|Pw\|_{L^2}+\|w\|_{L^2},
\end{multline}
since $a \in C_b^{\infty}(\rr^{2n+1})$. Notice from (\ref{rie2}) and (\ref{rie3}) that for those functions $w$, 
\begin{align}\label{rie4}
& \ \|\langle v \rangle^{-1}\Lambda_{t,x}^{-\frac{1}{2\sigma +1}}D_t w\|_{L^2}\\
\lesssim & \ \|\langle v \rangle^{-1}\Lambda_{t,x}^{-\frac{1}{2\sigma +1}}(\partial_t+ v\cdot \nabla_x) w\|_{L^2}+ \|\langle v \rangle^{-1}\Lambda_{t,x}^{-\frac{1}{2\sigma +1}} v\cdot \nabla_x w\|_{L^2} \notag \\ \notag
\lesssim  & \ \|(\partial_t+ v\cdot \nabla_x) w \|_{L^2} + \|\langle D_x\rangle^{\frac{2\sigma}{2\sigma+1}} w\|_{L^2} \lesssim  
 \|Pw\|_{L^2}+\|w\|_{L^2},
\end{align}
with $\Lambda_{t,x}=(1+D_t^2+|D_x|^2)^{1/2}$. It follows from (\ref{rie2}) and (\ref{rie4}) that 
\begin{multline*}
\|\Lambda_{t,x}^{\frac{2\sigma}{2\sigma +1}}\langle v \rangle^{-1}w\|_{L^2}
\lesssim  \|\langle v \rangle^{-1}\Lambda_{t,x}^{-\frac{1}{2\sigma +1}}D_tw\|_{L^2}+\|\langle v \rangle^{-1}\Lambda_{t,x}^{-\frac{1}{2\sigma +1}}D_xw\|_{L^2}\\ \lesssim
 \|Pw\|_{L^2}+\|w\|_{L^2}+\|\langle D_x \rangle^{\frac{2\sigma}{2\sigma+1}}w\|_{L^2} \lesssim \|Pw\|_{L^2}+\|w\|_{L^2}.
\end{multline*}
Another use of (\ref{rie2}) shows that 
\begin{equation}\label{rie5}
\|\Lambda_{t,x}^{\frac{2\sigma}{2\sigma +1}}\langle v \rangle^{-1}w\|_{L^2}+\|\langle D_v \rangle^{2\sigma}w\|_{L^2}  \lesssim \|Pw\|_{L^2}+\|w\|_{L^2}.
\end{equation}
Notice furthermore that 
\begin{multline*}
\|\langle D_{z}\rangle^{\frac{2\sigma}{2\sigma +1}}\langle v \rangle^{-1}w\|_{L^2}  \lesssim \|\Lambda_{t,x}^{\frac{2\sigma}{2\sigma +1}}\langle v \rangle^{-1}w\|_{L^2}+
\|\langle D_v \rangle^{\frac{2\sigma}{2\sigma +1}}\langle v \rangle^{-1}w\|_{L^2}  
\lesssim \|\Lambda_{t,x}^{\frac{2\sigma}{2\sigma +1}}\langle v \rangle^{-1}w\|_{L^2} \\ +
\|[\langle D_v \rangle^{\frac{2\sigma}{2\sigma +1}},\langle v \rangle^{-1}]w\|_{L^2}+\|\langle v \rangle^{-1}\langle D_v \rangle^{\frac{2\sigma}{2\sigma +1}}w\|_{L^2}
\lesssim \|Pw\|_{L^2}+\|w\|_{L^2},
\end{multline*}
because 
$$\|\langle v \rangle^{-1}\langle D_v \rangle^{\frac{2\sigma}{2\sigma +1}}w\|_{L^2} \lesssim \|\langle D_v \rangle^{\frac{2\sigma}{2\sigma +1}}w\|_{L^2} \lesssim \|\langle D_v \rangle^{2\sigma}w\|_{L^2}  \lesssim \|Pw\|_{L^2}+\|w\|_{L^2}$$
and that the operator $[\langle D_v \rangle^{\frac{2\sigma}{2\sigma +1}},\langle v \rangle^{-1}]$ is bounded on $L^2$ since symbolic calculus shows that its symbol belongs to the class $S^0$.
We therefore obtain from (\ref{rie5}) the estimate
\begin{equation}\label{rie6}
\|\langle D_{z}\rangle^{\frac{2\sigma}{2\sigma +1}}\langle v \rangle^{-1}w\|_{L^2}+\|\langle D_v \rangle^{2\sigma}w\|_{L^2}  \lesssim  \|Pw\|_{L^2}+\|w\|_{L^2},
\end{equation}
which may be extended to any function $w \in H^{2-\eps}(\rr^{2n+1})$ satisfying \eqref{rie1} and 
$$v \cdot \nabla_x w \in L^2(\rr^{2n+1}),$$ 
with $\eps>0$ such that $\max(2\sigma,1) \leq 2-\eps$. This is possible since $0<\sigma <1$.
Take now a $C_0^{\infty}(\rr)$ function $\psi$ such that 
$$\textrm{supp }\psi \subset [-3,3] \textrm{ and } \psi(t) = 1 \textrm{ if }|t| < 2.$$ 
Recalling that the function $u$ defined in (\ref{sev20}) belongs to $H^{-N_0}(\rr^{2n+1})$ and that for any $\eps >0$ such that $\max(2\sigma,1) \leq 2-\eps$, the symbol $M_{\delta}$ belongs to $S^{-N_0-2+\eps}$ for each fixed $0<\delta \leq 1$, we notice that
$$M_{\delta}(t,x,D_z)u \in H^{2-\eps}(\rr^{2n+1}),$$
for any $0<\delta \leq 1$.
It follows that the estimate (\ref{rie6}) may be applied to function
\begin{equation}\label{ts2}
\psi(t)\langle v \rangle^{-k}  M_\delta(t,x,D_z) u,
\end{equation}
with an integer $k \geq 1$. We obtain that 
\begin{multline}\label{rie7}
\|\langle D_{z}\rangle^{\frac{2\sigma}{2\sigma +1}}\psi(t)\langle v \rangle^{-(k+1)}  M_\delta(t,x,D_z) u\|_{L^2}+\|\langle D_v \rangle^{2\sigma}\psi(t)\langle v \rangle^{-k}  M_\delta(t,x,D_z) u\|_{L^2}  \\ 
\lesssim  \|P\psi(t)\langle v \rangle^{-k}  M_\delta(t,x,D_z) u\|_{L^2}+\|\psi(t)\langle v \rangle^{-k}  M_\delta(t,x,D_z) u\|_{L^2}.
\end{multline}
Notice that the choice of function $\chi$ in (\ref{sev6}) implies that the symbol 
$$(1-\psi(t))\langle v \rangle^{-l}M_{\delta}(t,x,\zeta),$$ 
with $l \geq 0$, belongs to the class $S^{-N_0-2}$ uniformly with respect to the parameter $0<\delta \leq 1$. Recalling that $0<\sigma<1$, it follows that 
\begin{align}\label{rie8}
\|\langle D_{z}\rangle^{\frac{2\sigma}{2\sigma +1}}(1-\psi(t))\langle v \rangle^{-(k+1)}  M_\delta(t,x,D_z) u\|_{L^2} \\ \lesssim \|(1-\psi(t))\langle v \rangle^{-(k+1)}  M_\delta(t,x,D_z) u\|_{\frac{2\sigma}{2\sigma +1}} \lesssim \|u\|_{-N_0} \notag
\end{align}
and
\begin{align}\label{rie9}
\|\langle D_v \rangle^{2\sigma}(1-\psi(t))\langle v \rangle^{-k}  M_\delta(t,x,D_z) u\|_{L^2} \\ \lesssim \|(1-\psi(t))\langle v \rangle^{-k}  M_\delta(t,x,D_z) u\|_{2\sigma} \lesssim \|u\|_{-N_0}, \notag
\end{align}
uniformly with respect to $0<\delta \leq 1$. Moreover, since the commutator 
$$[P,\psi(t)]=\psi'(t),$$ 
is bounded on $L^2$ and that 
\begin{align*}
& \ \|P\psi(t)\langle v \rangle^{-k}  M_\delta(t,x,D_z) u\|_{L^2} \\
\leq & \ \|\psi(t)P\langle v \rangle^{-k}  M_\delta(t,x,D_z) u\|_{L^2} + \|[P,\psi(t)]\langle v \rangle^{-k}  M_\delta(t,x,D_z) u\|_{L^2} \\
\lesssim & \ \|P\langle v \rangle^{-k}  M_\delta(t,x,D_z) u\|_{L^2}+ \|\langle v \rangle^{-k}  M_\delta(t,x,D_z) u\|_{L^2},
\end{align*}
we deduce from (\ref{rie7}), (\ref{rie8}), (\ref{rie9}) and the triangular inequality that
\begin{multline}\label{rie10}
\|\langle D_{z}\rangle^{\frac{2\sigma}{2\sigma +1}}\langle v \rangle^{-(k+1)}  M_\delta(t,x,D_z) u\|_{L^2}+\|\langle D_v \rangle^{2\sigma}\langle v \rangle^{-k}  M_\delta(t,x,D_z) u\|_{L^2}  \\ 
\lesssim  \|P\langle v \rangle^{-k}  M_\delta(t,x,D_z) u\|_{L^2}+\|\langle v \rangle^{-k}  M_\delta(t,x,D_z) u\|_{L^2}+ \|u\|_{-N_0},
\end{multline}
uniformly with respect to $0<\delta \leq 1$. Let $0<\eps_0<1$ such that 
\begin{equation}\label{ts1} 
0 < \varepsilon_0 < \frac{2\sigma}{2\sigma+1} \textrm{ and } (2\sigma-1)_++\eps_0 <2\sigma.
\end{equation}
We may use Lemma~\ref{lop3} to obtain the estimate
\begin{align}
\label{rie11} & \ \|P\langle v \rangle^{-k}  M_\delta(t,x,D_z) u\|_{L^2} \\ \notag
 \leq  & \ \|\langle v \rangle^{-k}  M_\delta(t,x,D_z) Pu\|_{L^2} +  \|[P,\langle v \rangle^{-k}  M_\delta(t,x,D_z)]u\|_{L^2}\\ \notag
 \lesssim  & \   \|\langle v \rangle^{-k}  M_\delta(t,x,D_z) Pu\|_{L^2}+ \|\langle D_v\rangle^{(2\sigma-1)_++\eps_0}\langle v \rangle^{-k} 
M_\delta(t,x,D_z) u \|_{L^2}  \\ \notag
& \ \quad + \| \langle v \rangle^{-k+1} 
M_\delta(t,x,D_z) u \|_{\eps_0} +\|u\|_{-N_0},
\end{align}
which holds uniformly with respect to $0<\delta \leq 1$. It follows from (\ref{rie10}) and (\ref{rie11}) that
\begin{align}\label{rie12}
& \ \|\langle D_{z}\rangle^{\frac{2\sigma}{2\sigma +1}}\langle v \rangle^{-(k+1)}  M_\delta(t,x,D_z) u\|_{L^2}+\|\langle D_v \rangle^{2\sigma}\langle v \rangle^{-k}  M_\delta(t,x,D_z) u\|_{L^2}  \\ 
\lesssim & \ \notag \|\langle v \rangle^{-k}  M_\delta(t,x,D_z) Pu\|_{L^2}+ \|\langle D_v\rangle^{(2\sigma-1)_++\eps_0}\langle v \rangle^{-k} 
M_\delta(t,x,D_z) u \|_{L^2} \\
& \ \quad \notag + \|\langle D_z \rangle^{\eps_0} \langle v \rangle^{-k+1} 
M_\delta(t,x,D_z) u \|_{L^2} 
+ \|u\|_{-N_0},
\end{align}
uniformly with respect to $0<\delta \leq 1$. Since for any $\eps>0$, there exists $C_{\eps}>0$ such that 
\begin{multline*}
\|\langle D_v\rangle^{(2\sigma-1)_++\eps_0}\langle v \rangle^{-k} M_\delta(t,x,D_z) u \|_{L^2}\\ \leq \eps \|\langle D_v \rangle^{2\sigma}\langle v \rangle^{-k}  M_\delta(t,x,D_z) u\|_{L^2}
+C_\eps \|\langle v \rangle^{-k}  M_\delta(t,x,D_z) u\|_{L^2},
\end{multline*}
because $(2\sigma-1)_++\eps_0 <2\sigma$, and
$$ \|\langle v \rangle^{-k}  M_\delta(t,x,D_z) u\|_{L^2} \lesssim \|\langle v \rangle^{-k+1}  M_\delta(t,x,D_z) u\|_{L^2} \lesssim \|\langle D_z \rangle^{\eps_0}\langle v \rangle^{-k+1}  M_\delta(t,x,D_z) u\|_{L^2},$$
because $\eps_0>0$, we deduce from (\ref{rie12}) that 
\begin{multline}\label{rie13}
 \|\langle D_{z}\rangle^{\frac{2\sigma}{2\sigma +1}}\langle v \rangle^{-(k+1)}  M_\delta(t,x,D_z) u\|_{L^2}  \\
\lesssim  \|\langle v \rangle^{-k}  M_\delta(t,x,D_z) Pu\|_{L^2} + \|\langle D_z \rangle^{\eps_0} \langle v \rangle^{-k+1} 
M_\delta(t,x,D_z) u \|_{L^2} 
+ \|u\|_{-N_0},
\end{multline}
uniformly with respect to $0<\delta \leq 1$. We shall need the following instrumental lemma:

\bigskip

\begin{lemma}\label{lop4} 
Let $s_1<s_2$, $N \geq 0$ and an integer $l \geq 2$. Then, there exists an integer $m \geq 2$ such that
$$\forall \eps >0, \exists C_{\eps}>0,  \forall w \in \mathscr{S}(\rr^{2n+1}), \ \|\langle  v \rangle^{l} w\|_{s_1}  \leq \varepsilon \| w \|_{s_2} + C_{\eps} \|\langle v \rangle^m w\|_{-N}.$$
\end{lemma}

\bigskip

\noindent
\textit{Proof of Lemma~\ref{lop4}}. Let $m \in \rr$. Since the symbol of the operator $\langle D_z \rangle^{m}$ belongs to the class $S(\langle \zeta \rangle^{m},dz^2+\langle \zeta \rangle^{-2}d\zeta^2)$ and 
$\langle v \rangle^m \in S(\langle v \rangle^{m},dz^2+\langle \zeta \rangle^{-2}d\zeta^2),$
symbolic calculus shows that the symbol of the operator 
$$\langle D_{z}\rangle^{-s_2} \langle  v \rangle^{l}\langle D_{z}\rangle^{2s_1 }\langle  v \rangle^{-l},$$ 
belongs to the class $S(\langle \zeta \rangle^{2s_1-s_2},dz^2+\langle \zeta \rangle^{-2}d\zeta^2)$.
It follows from the Cauchy-Schwarz inequality that for any $\eps>0$, there exists $C_{\eps}>0$ such that 
\begin{align*}
 \|\langle  v \rangle^{l} w\|_{s_1}^2 = & \ \|\langle D_{z}\rangle^{s_1} \langle  v \rangle^{l} w \|_{L^2}^2 =  (\langle D_{z}\rangle^{s_2} w,\langle D_{z}\rangle^{-s_2} \langle  v \rangle^{l}\langle D_{z}\rangle^{2s_1 }\langle  v \rangle^{l} w)_{L^2}\\
= & \  (\langle D_{z}\rangle^{s_2} w,\langle D_{z}\rangle^{-s_2} \langle  v \rangle^{l}\langle D_{z}\rangle^{2s_1 }\langle  v \rangle^{-l}\langle  v \rangle^{2l} w)_{L^2}\\
\leq & \  \|\langle D_{z}\rangle^{s_2} w\|_{L^2}\|\langle D_{z}\rangle^{-s_2} \langle  v \rangle^{l}\langle D_{z}\rangle^{2s_1 }\langle  v \rangle^{-l}\langle  v \rangle^{2l} w\|_{L^2}\\
\lesssim & \  \|w\|_{s_2}\|\langle  v \rangle^{2l} w\|_{s_1-(s_2-s_1)}\\
\lesssim & \  \eps \|w\|_{s_2}^2+C_{\eps}\|\langle  v \rangle^{2l} w\|_{s_1-(s_2-s_1)}^2.
\end{align*}
Notice that by assumption $s_1-(s_2-s_1)<s_1$. Lemma~\ref{lop4} then follows by $k_0$ iterations of this procedure with $s_1-k_0(s_2-s_1)<-N$.~$\Box$

\bigskip

\noindent
Recalling (\ref{ts1}), we may apply Lemma \ref{lop4} with $l=2$, $s_1=\eps_0$, $s_2=2\sigma/(2\sigma+1)$ and $N=N_0+s+1$. There therefore exists an integer $m \geq 2$ such that for all $\eps >0$,
\begin{equation}\label{ts4}
 \exists C_{\eps}>0,  \forall w \in \mathscr{S}(\rr^{2n+1}), \ \|\langle  v \rangle^{2} w\|_{\eps_0}  \leq \varepsilon \| w \|_{\frac{2\sigma}{2\sigma+1}} + C_{\eps} \|\langle v \rangle^m w\|_{-N_0-s-1}.
\end{equation}
We then choose the integer $k=m-1 \geq 1$ in (\ref{ts2}). Recalling that the symbol $M_{\delta}$ belongs to the class $S^{s+1}$ uniformly with respect to $0<\delta \leq 1$, we obtain that 
$$ \|\langle v \rangle^m w\|_{-N_0-s-1}=\|M_\delta(t,x,D_z) u\|_{-N_0-s-1} \lesssim \|u\|_{-N_0},$$
 uniformly with respect to $0<\delta \leq 1$ with $w=\langle v \rangle^{-(k+1)} M_\delta(t,x,D_z) u.$ 
The estimate (\ref{ts4}) extends by density. It follows that 
\begin{multline}\label{ts5}
\|\langle D_z \rangle^{\eps_0} \langle v \rangle^{-k+1} M_\delta(t,x,D_z) u \|_{L^2}=\|\langle v \rangle^{2}w \|_{\eps_0}  \leq \varepsilon \| w \|_{\frac{2\sigma}{2\sigma+1}}\\ + C_{\eps} \|\langle v \rangle^m w\|_{-N_0-s-1}
\lesssim  \varepsilon \| \langle D_{z}\rangle^{\frac{2\sigma}{2\sigma +1}}\langle v \rangle^{-(k+1)} M_\delta(t,x,D_z) u \|_{L^2} + C_{\eps}  \|u\|_{-N_0}.
\end{multline} 
We deduce from (\ref{rie13}) and (\ref{ts5}) that 
\begin{equation}\label{rie14}
 \|\langle D_{z}\rangle^{\frac{2\sigma}{2\sigma +1}}\langle v \rangle^{-(k+1)}  M_\delta(t,x,D_z) u\|_{L^2}  \\
\lesssim  \|\langle v \rangle^{-k}  M_\delta(t,x,D_z) Pu\|_{L^2} + \|u\|_{-N_0},
\end{equation}
uniformly with respect to $0<\delta \leq 1$.  Notice from (\ref{ts8}) and (\ref{sev6}) that the symbol of the operator  
$$\langle v \rangle^{-k}M_{\delta}(t,x,D_z)(1-\varphi(t,x))\langle v \rangle,$$ 
with $k \geq 1$, belongs to the class $S^{-N_0-2}$ uniformly with respect to $0<\delta \leq 1$. Recalling from (\ref{tt1}) that $M_{\delta} \in \tilde{S}^s$ uniformly with respect to $0<\delta \leq 1$, we obtain that 
\begin{align}\label{ts10}
& \ \|\langle v \rangle^{-k}  M_\delta(t,x,D_z) Pu\|_{L^2} \\ \notag
\leq & \  \|\langle v \rangle^{-k}  M_\delta(t,x,D_z)\varphi Pu\|_{L^2}+\|\langle v \rangle^{-k} M_\delta(t,x,D_z)(1-\varphi)Pu\|_{L^2} 
\\ \notag
\leq & \  \|M_\delta(t,x,D_z)\varphi Pu\|_{L^2}+\|\langle v \rangle^{-k} M_\delta(t,x,D_z)(1-\varphi)\langle v \rangle\langle v \rangle^{-1}Pu\|_{L^2}
\\ \notag
\lesssim & \ \|\varphi Pu\|_{s} +\|\langle v \rangle^{-1}Pu\|_{-N_0-2}.
\end{align}
Since obviously $\|\langle v \rangle^{-1}Pu\|_{-N_0-2} \lesssim \|u\|_{-N_0}$, we deduce from (\ref{rie14}) and (\ref{ts10}) that for all $0<\delta \leq 1$,
\begin{equation}\label{tt2}
 \|\langle D_{z}\rangle^{\frac{2\sigma}{2\sigma +1}}\langle v \rangle^{-(k+1)}  M_\delta(t,x,D_z) u\|_{L^2}  
\lesssim   \|\varphi Pu\|_{s} + \|u\|_{-N_0} <+\infty,
\end{equation}
since by assumption $u \in H^{-N_0}(\rr^{2n+1})$ and $\varphi Pu \in H^{s}(\rr^{2n+1})$. Let $\phi$ be $C_0^\infty(\rr_{t,x}^{n+1})$ function such that $\textrm{supp }\phi \subset \{|(t,x)|< 1\}$. Since the commutator $[\langle D_{z}\rangle^{\frac{2\sigma}{2\sigma +1}},\phi(t,x)]$ is obviously bounded on $L^2$, we deduce from (\ref{tt2}) that for all $0<\delta \leq 1$, 
\begin{align}\label{tt3}
& \   \|\langle D_{z}\rangle^{\frac{2\sigma}{2\sigma +1}}\phi(t,x)\langle v \rangle^{-(k+1)}  M_\delta(t,x,D_z) u\|_{L^2} \\ \notag
 \leq & \  \|[\langle D_{z}\rangle^{\frac{2\sigma}{2\sigma +1}},\phi(t,x)]\langle v \rangle^{-(k+1)}  M_\delta(t,x,D_z) u\|_{L^2} \\ \notag
 & \ \quad \quad  +\|\phi(t,x)\langle D_{z}\rangle^{\frac{2\sigma}{2\sigma +1}}\langle v \rangle^{-(k+1)}  M_\delta(t,x,D_z) u\|_{L^2}\\ \notag
 \lesssim & \  \|\langle v \rangle^{-(k+1)}  M_\delta(t,x,D_z) u\|_{L^2}
 +\|\langle D_{z}\rangle^{\frac{2\sigma}{2\sigma +1}}\langle v \rangle^{-(k+1)}  M_\delta(t,x,D_z) u\|_{L^2} \\\notag
 \lesssim & \ \|\langle D_{z}\rangle^{\frac{2\sigma}{2\sigma +1}}\langle v \rangle^{-(k+1)}  M_\delta(t,x,D_z) u\|_{L^2} \lesssim   \|\varphi Pu\|_{s} + \|u\|_{-N_0} <+\infty.
\end{align}
Notice from (\ref{sev6}) that 
$$\langle D_{z}\rangle^{\frac{2\sigma}{2\sigma +1}}\frac{\phi(t,x)}{\langle v \rangle^{k+1}} M_\delta(t,x,D_z) u=\langle D_{z}\rangle^{\frac{2\sigma}{2\sigma +1}}\frac{\phi(t,x)}{\langle v \rangle^{k+1}}\frac{\langle D_z \rangle^{s} }{
(1+ \delta \langle D_z \rangle)^{N_0+s+2}}  u$$
and that 
$$\langle D_{z}\rangle^{\frac{2\sigma}{2\sigma +1}}\frac{\phi(t,x)}{\langle v \rangle^{k+1}} M_\delta(t,x,D_z) u \longrightarrow
\langle D_z \rangle^{\frac{2\sigma}{2\sigma+1}}\frac{\phi(t,x)}{\langle v \rangle^{k+1}} \langle D_z \rangle^{s}u,$$
in $\mathscr{S}'(\rr^{2n+1})$ when $\delta \to 0$. Because of the weak compactness of the unit ball in $L^2$, it follows from (\ref{tt3}) that  
\begin{equation}\label{aay10}
\langle D_z \rangle^{\frac{2\sigma}{2\sigma+1}}\frac{\phi(t,x)}{\langle v \rangle^{k+1}} \langle D_z \rangle^{s}u \in L^2(\rr^{2n+1}),
\end{equation}
for any $C_0^\infty(\rr_{t,x}^{n+1})$ function $\phi$ such that $\textrm{supp }\phi \subset U_0=\{|(t,x)|< 1\}$.

To complete the proof of Corollary~\ref{ev1}, we notice that the following two assertions 
\begin{equation}\label{aay1}
\forall \phi \in C_0^\infty(U_0), \  \frac{\phi(t,x)}{\langle v \rangle^k}  u \in H^{r+m}(\rr^{2n+1}_{t,x,v})
\end{equation}
and
\begin{equation}\label{aay2}
\forall \phi \in C_0^\infty(U_0), \ \frac{ \phi(t,x) }{\langle v \rangle^k}\langle D_z \rangle^{m}  u \in H^r(\rr_{t,x,v}^{2n+1})
\end{equation}
are equivalent if $u \in H_{-N_0}(\rr^{2n+1}_{t,x,v})$ and $m, r \in \rr$. Assume first that (\ref{aay1}) holds.  Let $\phi$ be a $C_0^{\infty}(U_0)$ function and $\tilde{\phi}$ another $C_0^{\infty}(U_0)$ function such that $\tilde{\phi}=1$ on a neighborhood of $\textrm{supp }\phi$. Symbolic calculus allows to write
\begin{align}\label{aay3}
 & \ \frac{ \phi(t,x) }{\langle v \rangle^k}\langle D_z \rangle^{m}-\langle D_z \rangle^{m} \frac{ \phi(t,x) }{\langle v \rangle^k}\\ \notag
= & \ -\sum_{0 < |\alpha| < m +N_0+[r]+1} \frac{1}{\alpha!}  \textrm{Op}\Big(D_z^{\alpha}\left(\frac{ \phi(t,x) }{\langle v \rangle^k}\right) \partial_{\zeta}^{\alpha} \big(\langle \zeta \rangle^{m}\big)\Big) +R \\ \notag
= & \ -\sum_{0 < |\alpha| < m +N_0+[r]+1} \frac{1}{\alpha!}  \textrm{Op}\Big(D_z^{\alpha}\left(\frac{ \phi(t,x) }{\langle v \rangle^k}\right) \partial_{\zeta}^{\alpha} \big(\langle \zeta \rangle^{m}\big)\Big)\tilde{\phi} +\tilde{R}
\end{align}
with $R, \tilde{R} \in \textrm{Op}(S^{-N_0-[r]-1})$ and therefore 
\begin{equation}\label{aay4}
\|\tilde{R} u\|_{r} \lesssim \|u\|_{-N_0-[r]-1+r} \lesssim \|u\|_{-N_0} <+\infty. 
\end{equation}
In order to prove that (\ref{aay2}) holds, it is sufficient to check that 
\begin{multline*}
\Big\| \textrm{Op}\Big(D_z^{\alpha}\left(\frac{ \phi(t,x) }{\langle v \rangle^k}\right) \partial_{\zeta}^{\alpha} \big(\langle \zeta \rangle^{m}\big)\Big)\tilde{\phi}u\Big\|_r \\
=\Big\| \textrm{Op}\Big(D_z^{\alpha}\left(\frac{ \phi(t,x) }{\langle v \rangle^k}\right) \partial_{\zeta}^{\alpha} \big(\langle \zeta \rangle^{m}\big)\Big)\langle v \rangle^{k}\langle D_z \rangle^{-m}\langle D_z \rangle^{m}\frac{\tilde{\phi}(t,x)}{\langle v \rangle^k}u\Big\|_r
<+\infty,
\end{multline*}
when $0 < |\alpha| < m +N_0+[r]+1$. This is actually the case since symbolic calculus shows that the symbol of the operator 
$$\textrm{Op}\Big(D_z^{\alpha}\left(\frac{ \phi(t,x) }{\langle v \rangle^k}\right) \partial_{\zeta}^{\alpha} \big(\langle \zeta \rangle^{m}\big)\Big)\langle v \rangle^{k}\langle D_z \rangle^{-m},$$
belongs to $S^0$ when $0 < |\alpha| < m +N_0+[r]+1$ and that by (\ref{aay1}),
$$\langle D_z \rangle^{m}\frac{\tilde{\phi}(t,x)}{\langle v \rangle^k}u \in H^r(\rr_{t,x,v}^{2n+1}).$$
It follows that (\ref{aay1}) implies (\ref{aay2}).
Conversely, we may write for any function $\phi$ in $C_0^{\infty}(U_0)$, 
\begin{align*}
 & \ \frac{ \phi(t,x) }{\langle v \rangle^k}\langle D_z \rangle^{m}-\langle D_z \rangle^{m} \frac{ \phi(t,x) }{\langle v \rangle^k}\\ \notag
= & \ -\sum_{0 < |\alpha| < m +N_0+[r]+1} \frac{1}{\alpha!}  \textrm{Op}\Big(D_z^{\alpha}\left(\frac{ \phi(t,x) }{\langle v \rangle^k}\right) \partial_{\zeta}^{\alpha} \big(\langle \zeta \rangle^{m}\big)\Big)\langle D_z \rangle^{-m}\langle v \rangle^{k} \langle v \rangle^{-k} \langle D_z \rangle^{m}+R \\ \notag
= & \ -\sum_{0 < |\alpha| < m +N_0+[r]+1} \frac{1}{\alpha!}  \textrm{Op}\Big(D_z^{\alpha}\left(\frac{ \phi(t,x) }{\langle v \rangle^k}\right) \partial_{\zeta}^{\alpha} \big(\langle \zeta \rangle^{m}\big)\Big)\langle D_z \rangle^{-m}\langle v \rangle^{k}\frac{\tilde{\phi}(t,x)}{\langle v \rangle^{k}} \langle D_z \rangle^{m} +\tilde{R}
\end{align*}
with $R, \tilde{R} \in \textrm{Op}(S^{-N_0-[r]-1})$ if  $\tilde{\phi}$ is a $C_0^{\infty}(U_0)$ function satisfying $\tilde{\phi}=1$ on a neighborhood of $\textrm{supp }\phi$. This implies that 
$$\|\tilde{R} u\|_{r} \lesssim \|u\|_{-N_0-[r]-1+r} \lesssim \|u\|_{-N_0} <+\infty. $$
In order to prove that (\ref{aay1}) holds, it is sufficient to check that 
$$\Big\|\textrm{Op}\Big(D_z^{\alpha}\left(\frac{ \phi(t,x) }{\langle v \rangle^k}\right) \partial_{\zeta}^{\alpha} \big(\langle \zeta \rangle^{m}\big)\Big)\langle D_z \rangle^{-m}\langle v \rangle^{k}\frac{\tilde{\phi}(t,x)}{\langle v \rangle^{k}} \langle D_z \rangle^{m}u\Big\|_r
<+\infty,$$
when $0 < |\alpha| < m +N_0+[r]+1$. This is actually the case since symbolic calculus shows that the symbol of the operator 
$$\textrm{Op}\Big(D_z^{\alpha}\left(\frac{ \phi(t,x) }{\langle v \rangle^k}\right) \partial_{\zeta}^{\alpha} \big(\langle \zeta \rangle^{m}\big)\Big)\langle D_z \rangle^{-m}\langle v \rangle^{k},$$
belongs to $S^0$ when $0 < |\alpha| < m +N_0+[r]+1$ and that by (\ref{aay2}),
$$\frac{\tilde{\phi}(t,x)}{\langle v \rangle^k}\langle D_z \rangle^{m}u \in H^r(\rr_{t,x,v}^{2n+1}).$$
It follows that (\ref{aay2}) implies (\ref{aay1}).

By finally coming back to (\ref{aay10}), we deduce that the function defined in (\ref{sev20}) satisfies 
$$\frac{\phi(t,x)}{\langle v \rangle^{k+1}} u \in H^{s+\frac{2\sigma}{2\sigma+1}}(\rr^{2n+1}),$$
for any $C_0^\infty(\rr_{t,x}^{n+1})$ function $\phi$ such that $\textrm{supp }\phi \subset U_0=\{|(t,x)|< 1\}$, that is 
$$\frac{u}{\langle  v \rangle^{k+1}} \in H^{s+\frac{2\sigma}{2\sigma+1}}_{\textrm{loc}, (t_0,x_0)}(\rr^{2n+1}_{t,x,v}),$$
with $(t_0,x_0)=(0,0)$. This ends the proof of Corollary~\ref{ev1}.~$\Box$

\section{Appendix}

\subsection{The Boltzmann equation}\label{kkboltz}
The Boltzmann equation \cite{17,19} describes the behavior of a dilute gas when the only interactions taken into account are binary collisions. It reads as the evolution equation 
\begin{equation}\label{xiaoe1}
\begin{cases}
\partial_tf+v\cdot\nabla_{x}f=Q(f,f)\\
f|_{t=0}=f_0,
\end{cases}
\end{equation}
for the density distribution of the particles $f=f(t,x,v) \geq 0$  at time $t$, having position $x \in \rr^d$ and velocity $v \in \rr^d$.
The term appearing in the right-hand-side of this equation $Q(f,f)$ is the so-called quadratic Boltzmann collision operator associated to the Boltzmann bilinear operator 
\begin{equation}\label{xiaoeq1}
Q(g, f)=\int_{\rr^d}\int_{S^{d-1}}B(v-v_{*},\sigma) \big(g'_* f'-g_*f\big)d\sigma dv_*,
\end{equation}
acting on $L^1(\rr^d) \times L^1(\rr^d)$, $d \geq 2$, 
where $f'_*=f(t,x,v'_*)$, $f'=f(t,x,v')$, $f_*=f(t,x,v_*)$, $f=f(t,x,v)$, 
$$v'=\frac{v+v_*}{2}+\frac{|v-v_*|}{2}\sigma,\quad   v_*'=\frac{v+v_*}{2}-\frac{|v-v_*|}{2}\sigma,$$
for $\sigma$ belonging to the unit sphere $S^{d-1}$. Those relations between pre and post collisional velocities follow from the conservations of momentum and kinetic energy
$$\quad v+v_{\ast}=v'+v_{\ast}', \quad |v|^2+|v_{\ast}|^2=|v'|^2+|v_{\ast}'|^2,$$
in the binary collisions.
The Boltzmann operator has the fundamental properties of conserving mass, momentum and energy
$$\int_{\rr^d}Q(f,f)\phi(v)dv=0, \quad \phi(v)=1,v,|v|^2,$$
and to satisfy to the so-called Boltzmann's H theorem
$$-\frac{d}{dt} \int_{\rr^d}f \log f dv=-\int_{\rr^d}Q(f,f)\log f dv \geq 0.$$
The functional $-\int f \log f$ is the entropy of the solution and the Boltzmann's H theorem implies that at some point $x \in \rr^d$, any equilibrium distribution function, i.e., any function maximizing the entropy, has the form of a locally Maxwellian distribution
$$M(\rho,u,T)(v)=\frac{\rho}{(2\pi T)^{d/2}}e^{-\frac{|u-v|^2}{2T}},$$
where $\rho, u, T$ are respectively the density, mean velocity and temperature of the gas at the point $x$, defined by
$$\rho=\int_{\rr^d}f(v)dv, \quad u=\frac{1}{\rho}\int_{\rr^d}vf(v)dv, \quad T=\frac{1}{N\rho}\int_{\rr^d}|u-v|^2f(v)dv.$$ 
For further details on the physical background and derivation of the Boltzmann equation, we refer the reader to the extensive expositions \cite{17,19,villani2}.

For monatomic gas, the non-negative cross section $B(z,\sigma)$ defining the collision operator (\ref{xiaoeq1}) depends only on $|z|$ and the scalar product 
$\frac{z}{|z|}\cdot \sigma.$
Furthermore, the cross section is assumed to be supported in the set where $\frac{z}{|z|} \cdot \sigma \geq 0$. This condition is not a restriction when dealing with the quadratic Boltzmann operator $Q(f,f)$. As noticed in~\cite{alexandre0}, one may indeed reduce to this case after a symmetrization of the cross section since the term $f'f_*'$ appearing in the Boltzmann operator is invariant under the mapping $\sigma \rightarrow -\sigma$.

More specifically, we consider cross sections of  the type
\begin{equation}\label{xiaoeq1.01}
B(v-v_*,\sigma)=\Phi(|v-v_*|)b\Big(\frac{v-v_*}{|v-v_*|} \cdot \sigma\Big), \quad \cos \theta=\frac{v-v_*}{|v-v_*|} \cdot \sigma, \quad 0 \leq \theta \leq \frac{\pi}{2},
\end{equation}
with a kinetic factor
\begin{equation}\label{xiaosa0}
\Phi(|v-v_*|)=|v-v_*|^{\gamma}, \quad  \gamma \in ]-d,+\infty[,
\end{equation}
and a factor related to the collision angle with a singularity 
\begin{equation}\label{xiaosa1}
(\sin \theta)^{d-2} b(\cos \theta) \sim K \theta^{-1-2s},
\end{equation} 
when $\theta \to 0$, $\theta>0$, for some constants $K>0$ and $0 < s <1$. Notice that this singularity is not integrable
$$\int_0^{\frac{\pi}{2}}(\sin \theta)^{d-2}b(\cos \theta)d\theta=+\infty.$$
This non-integrability property plays a major r\^ole regarding the qualitative behavior of the solutions of the Boltzmann equation. Indeed, as first observed by L.~Desvillettes for the Kac equation~\cite{D95}, grazing collisions that account for the non-integrability of the angular factor near $\theta \sim 0^+$, 
do induce smoothing effects for the solutions of the non-cutoff Kac equation, or more generally for the solutions of the non-cutoff Boltzmann equation; whereas those solutions are at most as regular as the initial data, see e.g. \cite{36}, when the collision cross section is assumed to be integrable, or after removing the singularity by using a cutoff function (Grad's angular cutoff assumption).

Being concerned with a close to equilibrium framework, the setting of the problem can be formulated as follows. First of all, without loss of generality, we consider the fluctuation around
$$\mu(v)=(2\pi)^{-\frac{d}{2}}e^{-\frac{|v|^2}{2}},$$
the unique normalized equilibrium with mass 1, momentum 0 and temperature 1, by setting 
$$f=\mu+\sqrt{\mu}g.$$ 
Since $Q(\mu,\mu)=0$ by the conservation of the kinetic energy
$$|v|^2+|v_{\ast}|^2=|v'|^2+|v_{\ast}'|^2,$$
the Boltzmann collision operator can be split into three terms
$$Q(\mu+\sqrt{\mu}g,\mu+\sqrt{\mu}g)=Q(\mu,\sqrt{\mu}g)+Q(\sqrt{\mu}g,\mu)+Q(\sqrt{\mu}g,\sqrt{\mu}g),$$
whose linearized part is 
$$Q(\mu,\sqrt{\mu}g)+Q(\sqrt{\mu}g,\mu).$$
Define
$$\mathcal{L}g=\mathcal{L}_1g+\mathcal{L}_2g,$$
with 
$$\mathcal{L}_1g=-\mu^{-1/2}Q(\mu,\mu^{1/2}g), \quad \mathcal{L}_2g=-\mu^{-1/2}Q(\mu^{1/2}g,\mu),$$
the original Boltzmann equation is reduced to the Cauchy problem for the fluctuation
$$\begin{cases}
\partial_tg+v\cdot\nabla_{x}g+\mathcal{L}g=\mu^{-1/2}Q(\sqrt{\mu}g,\sqrt{\mu}g)\\
g|_{t=0}=g_0.
\end{cases}$$
This linearized Boltzmann operator $\mathcal{L}$ is known (see e.g. \cite{17}) to be an unbounded symmetric operator on $L^2$ (acting in the velocity variable) such that its Dirichlet form satisfies to
$$\mathbf{D}g=(\mathcal{L}g,g)_{L^2} \geq 0,$$
and that 
$$\mathbf{D}g=0 \Leftrightarrow g=\mathbf{P}g,$$ 
where
$$\mathbf{P}g=(a+b \cdot v+c|v|^2)\mu^{1/2},$$ 
with $a,c \in \rr$, $b \in \rr^d$, being the orthogonal projection onto the space of the so-called collisional invariants
$$\textrm{Span}\big\{\mu^{1/2},v_1 \mu^{1/2},...,v_d\mu^{1/2},|v|^2\mu^{1/2}\big\}.$$  
Furthermore, recent works in particular those by P.L.~Lions \cite{lions}, R.~Alexandre, L.~Desvillettes, C.~Villani, B.~Wennberg \cite{alexandre0}, C.~Mouhot~\cite{44}, C.~Mouhot, R.M.~Strain \cite{strain},
R. Alexandre, Y. Morimoto, C.-J. Xu, S. Ukai, T. Yang \cite{amuxy5-1}
 or P.T.~Gressman, R.M. Strain \cite{gress1,gress2} have highlighted that the non-cutoff Boltzmann operator enjoys remarkable coercive properties. Indeed, the results proved in~\cite{44,strain} show that the linearized Boltzmann operator enjoys the following coercive estimate 
\begin{equation}\label{xiaosa2}
(\mathcal{L}g,g)_{L^2(\rr^d)} \gtrsim \|\langle v \rangle^{s+\frac{\gamma}{2}}(g-\mathbf{P}g)\|_{L^2(\rr^d)}^2+\|g-\mathbf{P}g\|_{H^s_{\textrm{loc}}(\rr^d)}^2,
\end{equation}
with $\langle v \rangle=(1+|v|^2)^{1/2}$. This coercive estimate was improved in \cite{amuxy5-1} by showing that the latter term $\|g-\mathbf{P}g\|_{H^s_{\textrm{loc}}(\rr^d)}^2$ in (\ref{xiaosa2}) may be substituted by the global improvement 
$$\|\langle v \rangle^{\frac{\gamma}{2}}(g-\mathbf{P}g)\|_{H^s(\rr^d)}^2.$$
See~\cite{amuxy5-1,gress1,gress2} for optimal coercive estimates in appropriate weighted anisotropic Sobolev spaces. Regarding those coercive estimates,
the studies \cite{alexandre0,lions} have decisively made clear that the smoothing
effect of the non-cutoff Boltzmann equation is produced by the term
$$\int_{\rr^{2d}}\int_{S^{d-1}}B(v-v_{*},\sigma)g'_*(f- f')^2d\sigma dvdv_*,$$
thanks to the celebrated cancellation lemma, see~\cite{alexandre0} (Lemma~1 and Proposition~2). The discovery of these special features of the non-cutoff Boltzmann operator have led those authors to conjecture that this collision operator behaves and induces smoothing effects as a fractional Laplacian $-(-\Delta)^s$. See \cite{alexandre0}, p.~331. This conjecture accounts for the choice of the operator (\ref{yo2}) for the linearized non-cutoff Boltzmann equation. Notice that this operator is only a simplified model which does not account for the actual anisotropic features of the linearized non-cutoff Boltzmann equation and that current investigations are led to determine the exact microlocal structure of this operator. 
See the first works in this direction by R.~Alexandre~\cite{a1,a3,alexandre1}.

\subsection{Wick calculus}\label{appendix}

The purpose of this appendix is to recall few facts and basic properties of the Wick quantization with parameter that are used in this paper. We refer the reader to the works of N.~Lerner \cite{lerner}, \cite{cubo} and \cite{birkhauser} for thorough and extensive presentations of this quantization and some of its applications.

The main property of the Wick quantization is its property of positivity, i.e., that non-negative Hamiltonians define non-negative operators
$$a \geq 0 \Rightarrow a^{\textrm{Wick}(\lambda)} \geq 0.$$
We recall that this is not the case for the Weyl quantization nor the standard
quantization and refer to \cite{lerner} for an example of non-negative Hamiltonian defining an operator which is not non-negative. Before defining properly the Wick quantization, we first need to recall the definition of the wave-packets transform with parameter $\lambda>0$ of a function $u \in \mathscr{S}(\rr^n)$, 
$$W_{\lambda}u(y,\eta)=(u,\varphi_{y,\eta}^{\lambda})_{L^2(\rr^n)}, \ (y,\eta) \in \rr^{2n}.$$
where 
$$\varphi_{y,\eta}^{\lambda}(x)=(2\lambda)^{n/4}e^{- \pi \lambda |x-y|^2}e^{2i \pi (x-y) \cdot \eta}, \ x \in \mathbb{R}^n,$$
with $|\cdot|$ standing for the Euclidean norm and $x \cdot y$ the standard dot product on $\rr^n$. With this definition, one can check (Lemma 1.3 in \cite{cubo}) that 
the mapping $u \mapsto W_{\lambda}u$ is continuous from $\mathscr{S}(\rr^n)$ to $\mathscr{S}(\rr^{2n})$, isometric from $L^{2}(\rr^n)$ to $L^2(\rr^{2n})$ and that we have the
reconstruction formula
\begin{equation}\label{lay0.1}
\forall u \in \mathscr{S}(\rr^n), \forall x \in \rr^n, \ u(x)=\int_{\rr^{2n}}{W_{\lambda}u(y,\eta)\varphi_{y,\eta}^{\lambda}(x)dyd\eta}.
\end{equation}
By denoting by $\Sigma_Y^{\lambda}$ the operator defined by the Weyl quantization of the symbol 
$$p_Y(X)=2^n e^{-2\pi\Gamma_{\lambda}(X-Y)}, \ X=(x,\xi) \in \rr^{2n}, \ Y=(y,\eta) \in \rr^{2n},$$
with 
$$\Gamma_{\lambda}(X)=\lambda |x|^2+\frac{|\xi|^2}{\lambda}, \ X=(x,\xi) \in \rr^{2n},$$
which is a rank-one orthogonal projection,
$$\big{(}\Sigma_Y^{\lambda} u\big{)}(x)=W_{\lambda}u(Y)\varphi_Y^{\lambda}(x)=(u,\varphi_Y^{\lambda})_{L^2(\rr^n)}\varphi_Y^{\lambda}(x),$$
we define the Wick quantization with parameter $\lambda>0$ of any $L^{\infty}(\rr^{2n})$  symbol $a$ as
\begin{equation}\label{lay0.2}
a^{\textrm{Wick}(\lambda)}=\int_{\rr^{2n}}{a(Y)\Sigma_Y^{\lambda} dY}.
\end{equation}
More generally, one can extend this definition when the symbol $a$ belongs to $\mathscr{S}'(\rr^{2n})$ by defining the operator $a^{\textrm{Wick}(\lambda)}$ for any $u$ and $v$ in $\mathscr{S}(\rr^{n})$ by
$$\langle a^{\textrm{Wick}(\lambda)}u,\overline{v}\rangle_{\mathscr{S}'(\rr^{n}), \mathscr{S}(\rr^{n})}=\langle a(Y),(\Sigma_Y^{\lambda}u,v)_{L^2(\rr^n)}\rangle_{\mathscr{S}'(\rr^{2n}),  \mathscr{S}(\rr^{2n})},$$
where $\langle \textrm{\textperiodcentered},\textrm{\textperiodcentered}\rangle_{\mathscr{S}'(\rr^n),\mathscr{S}(\rr^n)}$ denotes the duality bracket between the
spaces $\mathscr{S}'(\rr^n)$ and $\mathscr{S}(\rr^n)$. 
Notice from Proposition~3.1 in~\cite{cubo} that 
\begin{equation}\label{lay10}
a^{\textrm{Wick}(\lambda)}=W_{\lambda}^*a^{\mu}W_{\lambda}, \ 1^{\textrm{Wick}(\lambda)}=\textrm{Id}_{L^2(\rr^n)},
\end{equation}
where $a^{\mu}$ stands for the multiplication operator by $a$ on $L^2(\rr^{2n})$.
The Wick quantization with parameter $\lambda>0$ is a positive quantization
\begin{equation}\label{lay0.5}
a \geq 0 \Rightarrow a^{\textrm{Wick}(\lambda)} \geq 0. 
\end{equation}
In particular, real Hamiltonians get quantized in this quantization by formally self-adjoint operators and one has (Proposition 3.1 in \cite{cubo}) that $L^{\infty}(\rr^{2n})$ symbols define bounded operators on $L^2(\rr^n)$ such that 
\begin{equation}\label{lay0}
\|a^{\textrm{Wick}(\lambda)}\|_{\mathcal{L}(L^2(\rr^n))} \leq \|a\|_{L^{\infty}(\rr^{2n})}.
\end{equation}
According to Proposition~3.2 in~\cite{cubo} (see (51) in~\cite{cubo}), the Wick and Weyl quantizations of a symbol $a$ are linked by the following identities
\begin{equation}\label{lay1bis}
a^{\textrm{Wick}(\lambda)}=\tilde{a}^w,
\end{equation}
with
\begin{equation}\label{lay2bis}
\tilde{a}(X)=\int_{\rr^{2n}}{a(X+Y)e^{-2\pi \Gamma_{\lambda}(Y)}2^ndY}, \ X \in \rr^{2n},
\end{equation}
and
\begin{equation}\label{lay1}
a^{\textrm{Wick}(\lambda)}=a^w+r_{\lambda}(a)^w,
\end{equation}
where $r_{\lambda}(a)$ stands for the symbol
\begin{equation}\label{lay2}
r_{\lambda}(a)(X)=\int_0^1\int_{\rr^{2n}}{(1-\theta)a''(X+\theta Y)Y^2e^{-2\pi \Gamma_{\lambda}(Y)}2^ndYd\theta}, \ X \in \rr^{2n}.
\end{equation}
We recall that we use here the following normalization for the Weyl quantization
\begin{equation}\label{lay3}
(a^wu)(x)=\int_{\rr^{2n}}{e^{2i\pi(x-y)\cdot\xi}a\big(\frac{x+y}{2},\xi\big)u(y)dyd\xi}.
\end{equation}
Let us finally mention that the operator $\pi_{\lambda}=W_{\lambda}W_{\lambda}^*$ is an orthogonal projection on a closed proper subspace of $L^2(\rr^{2n})$ whose kernel is given by
\begin{equation}\label{yo1}
K_{\lambda}(X,Y)=e^{-\frac{\pi}{2}\Gamma_{\lambda}(X-Y)}e^{i\pi (x-y)\cdot(\xi+ \eta)}, \ X=(x,\xi) \in \rr^{2n}, \ Y=(y,\eta) \in \rr^{2n}. 
\end{equation}
Indeed, for any $u \in \mathscr{S}(\rr_x^n)$ and $v \in \mathscr{S}(\rr_Y^{2n})$, we may write
\begin{multline*}
(W_{\lambda}u,v)_{L^2(\rr_Y^{2n})}=\int_{\rr_Y^{2n}}\Big(\int_{\rr_x^n}u(x)\overline{\varphi_{Y}^{\lambda}}(x)dx\Big)\overline{v}(Y)dY\\ =
\int_{\rr_x^{n}}u(x)\Big(\overline{\int_{\rr_Y^{2n}}\varphi_{Y}^{\lambda}(x)v(Y)dY}\Big)dx
=(u,W_{\lambda}^*v)_{L^2(\rr_x^n)}.
\end{multline*}
It follows that 
$$(W_{\lambda}^*v)(x)=\int_{\rr_Y^{2n}}\varphi_{Y}^{\lambda}(x)v(Y)dY.$$
Writing 
\begin{multline*}
(\pi_{\lambda}v)(X)=(W_{\lambda}W_{\lambda}^*v)(X)=\int_{\rr_t^n}\Big(\int_{\rr_Y^{2n}}\varphi_{Y}^{\lambda}(t)v(Y)dY\Big)\overline{\varphi_{X}^{\lambda}}(t)dt\\
=\int_{\rr_Y^{2n}}\Big(\int_{\rr_t^{n}}\varphi_{Y}^{\lambda}(t)\overline{\varphi_{X}^{\lambda}}(t)dt\Big)v(Y)dY=\int_{\rr_Y^{2n}}K_{\lambda}(X,Y)v(Y)dY,
\end{multline*}
with 
$$K_{\lambda}(X,Y)=\int_{\rr_t^{n}}\varphi_{Y}^{\lambda}(t)\overline{\varphi_{X}^{\lambda}}(t)dt.$$
By using the identity
$$2|t-y|^2+2|t-x|^2=|x-y|^2+|x+y-2t|^2,$$
an explicit computation gives that 
\begin{align*}
K_{\lambda}(X,Y)=& \ (2\lambda)^{\frac{n}{2}}\int_{\rr_t^{n}}e^{-\pi \lambda |t-y|^2}e^{-\pi \lambda |t-x|^2}e^{2i\pi(t-y)\cdot \eta}e^{-2i\pi(t-x)\cdot \xi}dt\\
= & \ (2\lambda)^{\frac{n}{2}}e^{-\frac{\pi}{2}\lambda |x-y|^2}e^{2i\pi (x\cdot \xi-y\cdot \eta)}\int_{\rr_t^{n}}e^{-\frac{\pi}{2} \lambda |x+y-2t|^2}e^{2i\pi t \cdot (\eta-\xi)}dt.
\end{align*}
Since 
\begin{multline*}
\int_{\rr_t^{n}}e^{-\frac{\pi}{2} \lambda |x+y-2t|^2}e^{2i\pi t \cdot (\eta-\xi)}dt=\int_{\rr_t^{n}}e^{-2\pi \lambda |t-\frac{x+y}{2}|^2}e^{2i\pi t \cdot (\eta-\xi)}dt
\\ =e^{i\pi (x+y) \cdot (\eta-\xi)}\int_{\rr_t^{n}}e^{-2\pi \lambda |t|^2}e^{2i\pi t \cdot (\eta-\xi)}dt=(2\lambda)^{-\frac{n}{2}}e^{i\pi (x+y) \cdot (\eta-\xi)}e^{-\frac{\pi}{2\lambda} |\eta-\xi|^2},
\end{multline*}
it follows that 
\begin{align*}
K_{\lambda}(X,Y)= & \ e^{-\frac{\pi}{2}\lambda |x-y|^2}e^{-\frac{\pi}{2\lambda} |\eta-\xi|^2}e^{2i\pi (x\cdot \xi-y\cdot \eta)}e^{i\pi (x+y) \cdot (\eta-\xi)}\\
= & \ e^{-\frac{\pi}{2}\Gamma_{\lambda}(X-Y)}e^{i\pi (x-y)\cdot(\xi+ \eta)}.
\end{align*}

\bigskip
\bigskip

\noindent
\textbf{Acknowledgements.} 
The research of the second author was
supported by  Grant-in-Aid for Scientific Research No.22540187,
Japan Society of the Promotion of Science. 
The third author is very grateful to the Japan Society for the Promotion of Science and JSPS London for supporting his stay at Kyoto University during the summer 2010.  He would like to thank Kyoto University and JSPS very warmly for their very kind hospitality and the exceptional working surroundings in which this work was done. Finally, all the three authors are very grateful to the referees for enriching comments and relevant suggestions which have helped to improve significantly the presentation of this article.


\begin{thebibliography}{aa}

\bibitem{a1}
R.~Alexandre, \textit{Remarks on 3D Boltzmann linear equation without cutoff}, Transport Theory Statist. Phys. 28 (5) (1999) 433-473
\bibitem{a3}
R.~Alexandre, \textit{Around 3D Boltzmann nonlinear operator without angular cutoff, a new formulation}, Math. Model. Numer. Anal. 34 (3) (2000) 575-590

\bibitem{alexandre0} R. Alexandre, L. Desvillettes, C. Villani, B. Wennberg, \textit{Entropy dissipation and long-range interactions}, Arch. Ration. Mech. Anal. 152 (2000) 327-355
\bibitem{alexandre1} R. Alexandre, M. Safadi, \textit{Littlewood Paley decomposition and regularity issues in Boltzmann homogeneous equations I. Non cutoff and Maxwell cases}, Math. Models Methods Appl. Sci. 15 (6) (2005) 907-920
\bibitem{al-saf-1}R. Alexandre, M. Safadi, \textit{Littlewood-Paley theory and regularity issues in
Boltzmann homogeneous equations II. Non cutoff case and non
Maxwellian molecules}, {Discrete Contin. Dyn. Syst}.  24 (2009)
1-11
\bibitem{amuxy2} R. Alexandre, Y. Morimoto, S. Ukai, C.-J. Xu, T. Yang, \textit{Uncertainty principle and regularity for Boltzmann equation}, J. Funct. Anal. 255 (8) (2008) 2013-2066
\bibitem{amuxy3} R. Alexandre, Y. Morimoto, S. Ukai, C.-J. Xu, T. Yang, \textit{Regularizing effect and local existence
for non-cutoff Boltzmann equation}, { Arch. Ration. Mech. Anal.} 198 (2010) 39-123

\bibitem{amuxy4} R. Alexandre, Y. Morimoto, S. Ukai, C.-J. Xu, T. Yang,  \textit{Global existence and full regularity of the Boltzmann
equation without angular cutoff,} to appear in {Comm. Math. Phys.}
 http://hal.archives-ouvertes.fr/hal-00439227/fr/

\bibitem{amuxy5-1} R. Alexandre, Y. Morimoto, S. Ukai, C.-J. Xu, T. Yang,
\textit{The Boltzmann equation without angular cutoff in the whole space I: Global existence for soft potential}, {to appear in J. Funct. Anal.} http://hal.archives-ouvertes.fr/hal-00496950/fr/


\bibitem{amuxy5-3} R. Alexandre, Y. Morimoto, S. Ukai, C.-J. Xu, T. Yang,
\emph{Boltzmann equation without angular cutoff in the whole space III: Qualitative properties of solutions},
preprint,
http://hal.archives-ouvertes.fr/hal-00510633/fr/

\bibitem{amuxy7} R. Alexandre, Y. Morimoto, S. Ukai, C.-J. Xu, T. Yang, \textit{Smoothing effect for spatially homogeneous Boltzmann equation}, work in progress

\bibitem{bouchut} F. Bouchut, \textit{Hypoelliptic regularity in kinetic equations}, J. Math. Pures Appl. (9) 81  (2002), no. 11,
1135-1159

\bibitem{17} C.~Cercignani, \textit{The Boltzmann Equation and its Applications}, Applied Mathematical Sciences, vol. 67 (1988), Springer-Verlag, New York 

\bibitem{19} C.~Cercignani, R.~Illner, M.~Pulvirenti, \textit{The mathematical theory of dilute gases}, (1994) Springer-Verlag, New York 



\bibitem{chen1}
H.~Chen, W.-X.~Li, C.-J. Xu, \textit{Gevrey hypoellipticity for linear and non-linear Fokker-Planck equations}, J. Differential Equations 246 (2009), no. 1, 320-339
\bibitem{chen}
H.~Chen, W.-X.~Li, C.-J.~Xu, \textit{Gevrey hypoellipticity for a class of kinetic equations}, Comm. Partial Differential Equations 36 (2011), no. 4,  693-728
\bibitem{desv1} L. Desvillettes, B. Wennberg, \textit{Smoothness of the solution of the spatially homogeneous Boltzmann equation without cutoff}, Comm. Partial Differential Equations, 29 (1-2) (2004) 133-155


\bibitem{D95} L.~Desvillettes, \textit{About the regularization properties
of the non cut-off Kac equation}, Comm. Math. Phys.  168 (1995)
417-440

\bibitem{gress1} P.T.~Gressman, R.M. Strain, \textit{Global Classical Solutions of the Boltzmann Equation without Angular Cut-off}, J. Amer. Math. Soc. 
24 (2011), no. 3, 771-847

\bibitem{gress2} P.T.~Gressman, R.M. Strain, \textit{Sharp anisotropic estimates for the Boltzmann collision operator and its entropy production}, preprint (2010)

\bibitem{landau} F.~H\'erau, K.~Pravda-Starov, \textit{Anisotropic hypoelliptic estimates for Landau-type operators}, J. Math. Pures Appl. 95 (2011) 513-552
\bibitem{hormander} L. H\"{o}rmander, {\it The analysis of linear partial differential operators} (vol. I--IV), Springer Verlag~(1985)
\bibitem{HMUY} Z.H. Huo, Y. Morimoto, S. Ukai, T. Yang, \textit{Regularity of
solutions for spatially homogeneous Boltzmann equation without
Angular cutoff}, {Kinetic and Related Models},  1 (2008)
453-489
\bibitem{L-X} N. Lekrine, C.-J. Xu, \textit{Gevrey regularizing effect of Cauchy problem
for non-cutoff homogeneous Kac's equation},
 Kinetic and Related Models,  2 (2009) 647-666
\bibitem{lerner}
N. Lerner, \textit{The Wick calculus of pseudodifferential operators and some of its applications}, Cubo Mat. Educ. 5 (2003) 213-236 
\bibitem{cubo} N. Lerner, \textit{Some Facts About the Wick Calculus}, Pseudo-differential operators. Quantization and signals. Lecture Notes in Mathematics, 1949, 135--174. Springer-Verlag, Berlin (2008) 
\bibitem{birkhauser} N. Lerner, \textit{Metrics on the phase space and non-selfadjoint pseudo-differential operators}, Pseudo-Differential Operators. Theory and Applications, 3. Birkh\"auser Verlag, Basel (2010)

\bibitem{lions}  P.L. Lions, \textit{Regularity and compactness for Boltzmann collision operator without angular
cut-off}, C. R. Acad. Sci. Paris Series I, 326 (1998) 37-41 


\bibitem{MU} Y. Morimoto, S. Ukai,
\textit{Gevrey smoothing effect of solutions for  spatially
homogeneous nonlinear Boltzmann  equation without angular cutoff},
 J. Pseudo-Differ. Oper. Appl. 1 (2010) 139-159 

\bibitem{Mo2} Y. Morimoto, S. Ukai, C-J. Xu, T. Yang, \textit{Regularity of solutions to the spatially homogeneous Boltzmann equation without angular cutoff}, Discrete Contin. Dyn. Syst. 24 (2009), no. 1, 187-212



\bibitem{MoXu}   Y. Morimoto, C-J. Xu, \textit{Hypoellipticity for a class of kinetic equations}, J. Math. Kyoto Univ. 47, no. 1 (2007) 129-152

\bibitem{Mo1} Y. Morimoto, C-J. Xu, \textit{Ultra-analytic effect of Cauchy problem for a class of kinetic equations}, J. Differential Equations, 247 (2) (2009) 596-617


\bibitem{44}
C. Mouhot, \textit{Explicit coercivity estimates for the Boltzmann and Landau operators}, Comm. Partial Differential Equations, 31 (2006) 1321-1348

\bibitem{strain}
C. Mouhot, R.M. Strain, \textit{Spectral gap and coercivity estimates for linearized Boltzmann collision operators without angular cutoff}, 
J. Math. Pures Appl. (9) 87 (2007), no. 5, 515-535 

\bibitem{perthame} B. Perthame, \textit{Higher Moments for Kinetic Equations: the Vlasov-Poisson and Fokker-Planck cases}, Math. Methods in the Applied Sciences, 13 (1990) 441-452 

\bibitem{lms} K. Pravda-Starov, \textit{A complete study of the pseudo-spectrum for the rotated harmonic oscillator}, J. London Math. Soc. (2) 73 (2006) 745-761 
\bibitem{Ukai} S. Ukai, \textit{Local solutions in Gevrey classes to the nonlinear Boltzmann equation without cutoff}, Japan J. Appl. Math. 1 (1) (1984) 141-156
\bibitem{Villani} C. Villani, \textit{On a new class of weak solutions to the spatially homogeneous Boltzmann and Landau equations}, Arch. Ration. Mech. Anal. 143 (1998) 273-307
\bibitem{villani2}
C. Villani, {\it A review of mathematical topics in collisional kinetic theory}, Handbook of Mathematical Fluid Dynamics, vol. I, North-Holland, Amsterdam (2002) 71-305

\bibitem{36}
B. Wennberg, \textit{Regularity in the Boltzmann equation and the Radon transform}, Comm. Partial Differential Equations, 19 (1994) 2057-2074


\end{thebibliography}
\end{document}